\documentclass[12pt]{article}
\usepackage[a4paper,vmargin={3cm,3cm},hmargin={2.5cm,2.5cm},bottom=30mm]{geometry}
\pdfoutput=1
\usepackage{dsfont}
\usepackage{mathrsfs}
\usepackage{multirow}
\usepackage{slashbox,multirow}
\usepackage{threeparttable}
\usepackage{lscape}
\usepackage[sectionbib]{natbib}
\usepackage{fancyhdr,latexsym,amsmath,amsfonts,amssymb,amsbsy,amsthm,url}
\usepackage{amsmath}
\usepackage[utf8]{inputenc}
\usepackage{tikz}
\date{}

 \usepackage{color}
\usepackage[pdfpagemode=UseNone,bookmarksopen=false,colorlinks=blue]{hyperref}
\hypersetup{ colorlinks=true, linkcolor=blue, citecolor=blue,
  filecolor=blue, urlcolor=red}

\allowdisplaybreaks
\bibfont{\footnotesize}

\newtheorem{assumption}{Assumption}

\newtheorem{lemma}{Lemma}

\newtheorem{proposition}{Proposition}[section]
\newtheorem{theorem}{Theorem}[section]
\theoremstyle{remark}

\begin{document}

%%%%%%%%%%%%%%%%%%%%%%%%%%%%%%%%%%%%%%%%%%%%%%%%%%%%%%%%%%%%%%%%%%%%%%%%%%%%%%%%%%%%%%%%%%%%%%%%%%%%%%%%%%%%%%%%%%%%%%%%%%%%
%%%%%%%%%%%%%%%%%%%%%%%%%%%%%%%%%%%%%%%%%%%%%%%%%%%%%%%%%%%%%%%%%%%%%%%%%%%%%%%%%%%%%%%%%%%%%%%%%%%%%%%%%%%%%%%%%%%%%%%%%%%%

%\renewcommand{\baselinestretch}{2}

\title{\textbf{On estimators of the mean of infinite dimensional data in finite populations}}
%\runtitle{Quantile processes in finite populations}
\author{Anurag Dey and Probal Chaudhuri\\
\textit{Indian Statistical Institute, Kolkata}}
\maketitle

%%%%%%%%%%%%%%%%%%%%%%%%%%%%%%%%%%%%%%%%%%%%%%%%%%%%%%%%%%%%%%%%%%%%%%%%%%%%%%%%%%%%%%%%%%%%%%%%%%%%%%%%%%%%%%%%%%%%%%%%%%%%

\begin{abstract}
The Horvitz-Thompson (HT), the Rao-Hartley-Cochran (RHC) and the generalized regression (GREG) estimators of the finite population mean are considered, when the observations are from an infinite dimensional space. We compare these estimators based on their asymptotic distributions under some commonly used sampling designs and some superpopulations satisfying linear regression models. We show that the GREG estimator is asymptotically at least as efficient as any of the other two estimators under different sampling designs considered in this paper. Further, we show that the use of some well known sampling designs utilizing auxiliary information may have an adverse effect on the performance of the GREG estimator, when the degree of heteroscedasticity present in linear regression models is not very large. On the other hand, the use of those sampling designs improves the performance of this estimator, when the degree of heteroscedasticity present in linear regression models is large.  We develop methods for determining the degree of heteroscedasticity, which in turn determines the choice of appropriate sampling design to be used with the GREG estimator.  We also investigate the consistency of the covariance operators of the above estimators. We carry out some numerical studies using real and synthetic data, and our theoretical results are supported by the results obtained from those numerical studies.\\ 
\end{abstract}
\newpage


\section{Introduction}
In the recent past, \cite{cardot2011horvitz}, \cite{cardot2013confidence}, \cite{cardot2014variance}, etc. considered the HT estimator (see \cite{horvitz1952generalization}) of the finite population mean, when population observations are from some functional space. \cite{cardot2013comparison} and \cite{cardot2013uniform} also constructed a model assisted estimator for finite population mean function based on some  homoscedastic  linear regression models. This model assisted estimator can be related to the GREG estimator considered  earlier  in \cite{MR1707846} for finite dimensional data.  All these authors  investigated different asymptotic properties of the HT and the model assisted estimators in $\mathcal{C}[0,T]$, the space of continuous functions defined on $[0, T]$, under sampling designs, which satisfy some regularity conditions. These sampling designs include the simple random sampling without replacement (SRSWOR), stratified sampling design with SRSWOR  and  the rejective sampling designs (see \cite{MR0178555}), etc. 
\par
\vspace{.1cm}

Suppose that $\mathcal{P}$ is a finite population consisting of $N$ units $U_1,\ldots,U_N$, and $s\subset\mathcal{P}$ is a sample of size $n$ $(< N)$. Let $\mathcal{S}$ be the collection of all subsets of $\mathcal{P}$ having size $n$. Then, a  without replacement  sampling design $P(s)$ is a probability distribution on $\mathcal{S}$ such that $0 \leq P(s)\leq 1$ for all $s\in \mathcal{S}$ and $\sum_{s\in\mathcal{S}}P(s)$=$1$. In this article, we consider SRSWOR, Lahiri-Midzuno-Sen (LMS) sampling design (see \cite{Lahiri1951sampling}, \cite{Midzuno1952sampling} and \cite{Sen1953sampling}), RHC sampling design (see \cite{rao1962simple}) and high entropy $\pi$PS (HE$\pi$PS) sampling designs (see  Section  \ref{sec 4}). It is to be noted that some size variable is used for implementing all of above the sampling designs except SRSWOR. In this article, we consider the extensions of the HT and the RHC estimators for the population mean of a study variable $y$ in an infinite dimensional separable Hilbert space $\mathcal{H}$ because these estimators are widely used design unbiased estimators of the population mean for finite dimensional data. We also consider the extension of the GREG estimator for the population mean of $y$, which is not a design unbiased estimator but known to be asymptotically  often more efficient than other estimators  under SRSWOR for finite dimensional data (see \cite{MR0474575}). We compare the HT, the RHC and the GREG estimators using their asymptotic distributions under above sampling designs and some superpopulations satisfying linear regression models. The main results obtained from this comparison are the following. 
\begin{itemize}
\item The GREG estimator is asymptotically at least as good as the HT estimator under each of SRSWOR, LMS sampling design and any HE$\pi$PS sampling design. Also, the GREG estimator turns out to be asymptotically at least as good as the RHC estimator under RHC sampling design.

\item If the degree of heteroscedasticity present in linear regression models is not very large, then the use of the well known sampling designs like RHC and any HE$\pi$PS sampling designs instead of SRSWOR may have an adverse effect on the performance of the GREG estimator. In other words, the use of the auxiliary information in the design stage of sampling may have an adverse effect on the performance of the GREG estimator. On the other hand, if the degree of heteroscedasticity present in linear regression models is sufficiently large, then the sampling designs like RHC and any HE$\pi$PS sampling designs lead to an improvement in the performance of the GREG estimator. 
\end{itemize}
In  Section  \ref{sec 2}, we discuss infinite dimensional extensions of the HT, the RHC and  the GREG estimators of the population mean. In  Section  \ref{sec 4}, we compare these estimators using their asymptotic distributions under the sampling designs mentioned above and some superpopulations satisfying linear regression models. In this section, we also discuss the estimation of asymptotic covariance operators of several estimators and show that these estimators of asymptotic covariance operators are consistent. Some numerical results based on both synthetic and real data are presented in  Section  \ref{sec 6}.  Several methods of determining the degree of heteroscedasticity present in linear regression models are provided in  Section  \ref{sec 5}. Description of several sampling designs considered in this article are given in Appendix A.  Proofs of various results are given in Appendix B.
\section{Estimators based on infinite dimensional data}\label{sec 2} 
Suppose that $\mathcal{H}$ is an infinite dimensional separable Hilbert space with associated inner product $\langle\cdot,\cdot\rangle$, and $Y_1,\ldots, Y_N$ are population values of a study variable $y$ in $\mathcal{H}$. The HT estimator of the finite population mean of $y$, $\overline{Y}$=$\sum_{i=1}^N Y_i/N$, is defined as 
\begin{equation}\label{eq 4}
\hat{\overline{Y}}_{HT}=\sum_{i\in s}(N\pi_i)^{-1}Y_i,
\end{equation}
where  $\pi_i$=$\sum_{s\ni U_i}P(s)$ is the inclusion probability of the $i^{th}$ population unit $U_i$ for $i$=$1,\ldots,N$.  For notational convenience, we shall write  $i\in s$ (or $s\ni i$) instead of $U_i\in s$ (or $s\ni U_i$). 
\par
\vspace{.1cm}

 Before we write the expression of the RHC estimator, we describe the RHC sampling design briefly. Suppose that $X_1,\ldots,X_N$ are known population values on a size variable $x$ in $(0,\infty)$. Then, in the RHC sampling design (see \cite{rao1962simple}, \cite{chaudhuri2006feasibility}, etc.), $\mathcal{P}$ is first divided randomly into $n$ disjoint groups, say $\mathcal{P}_1,\ldots, \mathcal{P}_n$ of sizes $N_1,\ldots,N_n$ such that $\sum_{i=1}^n N_i$=$N$, and then one unit is selected from each group independently.  For more details on RHC sampling design, see Appendix A.  Now, suppose that $Q_i$ is the total of the $x$ values of that randomly formed group from which the unit $U_i$ is selected in the sample $s$. Then, the RHC estimator of $\overline{Y}$ can be expressed as 
\begin{equation}\label{eq 5}
 \hat{\overline{Y}}_{RHC}=\sum_{i\in s}(NX_i)^{-1} Q_iY_i.
\end{equation} 
\par
\vspace{.1cm}

\cite{rao1962simple} considered the RHC estimator for a real valued study variable. The RHC estimator is easily computable than other unbiased estimators under other unequal probability sampling designs without replacement (e.g., the HT or the Des Raj estimator (see \cite{MR0474578}) under probability proportional to size sampling without replacement). Moreover, the RHC estimator has smaller variance than the usual unbiased estimator under the probability proportional to size sampling with replacement. Also, its variance can be estimated by a non negative unbiased estimator. These results  continue to hold,  when we consider the RHC estimator for a $\mathcal{H}$ valued study variable. 
\par
\vspace{.1cm}

\cite{robinson1983asymptotic}, \cite{sarndal1992model}, \cite{deville1992calibration}, \cite{MR1707846}, etc. considered the GREG estimator for finite dimensional data. Suppose that $z$=$(z_1,\ldots,z_d)$ is a $\mathbb{R}^d$-valued ($d\geq 1$) covariate with population values $Z_1,\ldots,Z_N$ and known population total $\sum_{i=1}^N Z_i$. It will be appropriate to note that the size variable $x$ may be one of the real valued components of $z$ in some cases.  All vectors in Euclidean spaces will be taken as row vectors and superscript $T$ will be used to denote their transpose.  Further, suppose that $\mathcal{G}$ is any separable Hilbert space with inner product $\langle\cdot,\cdot\rangle$, and $\mathcal{B}(\mathcal{G},\mathcal{H})$ is the class of all bounded linear operators from $\mathcal{G}$ to $\mathcal{H}$.  It is to be noted that  $\mathcal{B}(\mathcal{G},\mathcal{H})$ is an infinite dimensional Hilbert space associated with the Hilbert-Schmidt (HS) inner product (see \cite{hsing2015theoretical}).  For any $a\in\mathcal{G}$ and $b\in\mathcal{H}$, let us consider the tensor product $a\otimes b\in \mathcal{B}(\mathcal{G},\mathcal{H})$, which is defined as $(a\otimes b)e$=$\langle a,e\rangle b$, $e\in\mathcal{G}$. Suppose that  $\hat{\overline{Z}}$=$\big(\sum_{i\in s}\pi_i^{-1}\big)^{-1}\sum_{i\in s}\pi_i^{-1}Z_i$.  Let us also suppose that the inverse of  $\hat{S}_{zz}$=$\big(\sum_{i\in s}\pi_i^{-1}\big)^{-1}\sum_{i\in s}\pi_i^{-1}(Z_i-\hat{\overline{Z}})^T(Z_i-\hat{\overline{Z}})$  exists. Then, an infinite dimensional version of the GREG estimator for the population mean is defined as 
\begin{equation}\label{eq 6}
\hat{\overline{Y}}_{GREG}=\hat{\overline{Y}}+\hat{S}_{zy}((\overline{Z}-\hat{\overline{Z}})\hat{S}_{zz}^{-1}),
\end{equation} 
where  $\overline{Z}$=$N^{-1}\sum_{i=1}^N Z_i$, $\hat{\overline{Y}}$=$\big(\sum_{i\in s}\pi_i^{-1}\big)^{-1}\sum_{i\in s}\pi_i^{-1}Y_i$ and $\hat{S}_{zy}$=$\big(\sum_{i\in s}\pi_i^{-1}\big)^{-1}\sum_{i\in s}\pi_i^{-1}(Z_i-\hat{\overline{Z}})\otimes(Y_i-\hat{\overline{Y}})$.  Under RHC sampling design, we consider the GREG estimator $\hat{\overline{Y}}_{GREG}$ after replacing $\pi_i^{-1}$ by  $Q_iX_i^{-1}$  (cf. \cite{MR1707846}).
\section{Comparison of estimators under superpopulation models}\label{sec 4} 
In this section, we compare among the HT and the GREG estimators under SRSWOR, LMS sampling design and HE$\pi$PS sampling designs, and the RHC and the GREG estimators under RHC sampling design. For this, we first assume that the observations $\{(Y_i,Z_i,X_i):1\leq i\leq N\}$ are independent and identically distributed (i.i.d.) $\mathcal{H}\times\mathbb{R}^{d+1}$ valued random variables on a probability space $(\Omega,\mathcal{F},\mathbb{P})$. Earlier,  \cite{sarndal1992model},  \cite{rao2003small}, \cite{fuller2011sampling}, \cite{chaudhuri2014modern} (see chap. $5$), \cite{MR3670194}, etc. considered superpopulation models for constructing different estimators and studying their behaviour, when the observations on $y$ are from some finite dimensional space. Note that the population values $X_1,\ldots,X_N$, which are now random variables on $\Omega$, are used to implement sampling designs like LMS, RHC and HE$\pi$PS sampling designs. In such cases, a function $P(s,\omega)$ is defined on $\mathcal{S}\times\Omega$ (see \cite{MR3670194}) so that for each $s\in \mathcal{S}$, $P(s,\omega)$ is a random variable on $\Omega$, and for each $\omega\in \Omega$, $P(s,\omega)$ is a probability distribution on $\mathcal{S}$. It is to be noted that $P(s,\omega)$ is a sampling design for each $\omega\in \Omega$. If a sampling design does not depend on $X_1,\ldots,X_N$ (e.g., SRSWOR), then for each $s\in \mathcal{S}$, the corresponding $P(s,\omega)$ is a degenerate random variable on $\Omega$. We define our asymptotic framework as follows. Let $\{\mathcal{P}_{\nu}\}$ be a sequence of nested populations with $N_\nu$, $n_\nu \rightarrow \infty$ as $\nu \rightarrow \infty$ (see \cite{MR648029}), where $N_\nu$ and $n_{\nu}$ are, respectively, the population size and the sample size corresponding to the $\nu^{th}$ population. We shall frequently drop the limiting index $\nu$ for the sake of notational simplicity. 
\par
\vspace{.1cm}

We now discuss HE$\pi$PS sampling designs. Suppose that a sampling design $P(s,\omega)$ is such that the Kullback–Leibler divergence  $D(P||R)$=$\sum_{s \in \mathcal{S}}P(s,\omega)\log\big(P(s,\omega)(R(s,\omega))^{-1}\big)$  converges to $0$ as $\nu\rightarrow\infty$ \textit{a.s.} $[\mathbb{P}]$ for some rejective sampling design $R(s,\omega)$ ( for the description of rejective sampling design, see Appendix A).  Such a sampling design is known as the high entropy sampling design (cf. \cite{MR1624693}, \cite{cardot2014variance}, \cite{MR3670194}, etc.). We call a sampling design $P(s,\omega)$ a HE$\pi$PS sampling design if it is a high entropy sampling design as well as a $\pi$PS ( i.e., inclusion probability $\pi$ proportional to size\color{black}) sampling design (see  \cite{bondesson2006pareto} and the references therein). An example of the HE$\pi$PS sampling design is the Rao-Sampford (RS) sampling design (see \cite{MR1624693}  and Appendix A). 
\par
\vspace{.1cm}

 Before we state our main results, let us consider some assumptions on distributions of $\{Y_i,Z_i,X_i\}_{i=1}^N$. Suppose that $E_{\mathbb{P}}$ denotes that expectation with respect to the probability measure $\mathbb{P}$  given above. The expectations of $\mathcal{H}$ valued random variables are defined using Bochner integrals (see \cite{hsing2015theoretical}).  Note that in any finite dimensional Euclidean space, we consider the Euclidean norm and denote it by $||\cdot||$. On the other hand, in $\mathcal{H}$, we consider the norm induced by the inner product associated with $\mathcal{H}$ and denote it by $||\cdot||_{\mathcal{H}}$. 
\begin{assumption}\label{ass C1}
 $n N^{-1}\rightarrow \lambda$ as $\nu\rightarrow\infty$,  where $0\leq \lambda<1$.
\end{assumption}
\begin{assumption}\label{ass C2}
 $X_i\leq b$ \textit{a.s.} $[\mathbb{P}]$ for some $0<b<\infty$, $E_{\mathbb{P}}(X_i)^{-2}<\infty$, and $\max_{1\leq i\leq N} X_i\times $\\ $\big(\min_{1\leq i\leq N} X_i\big)^{-1}$=$O(1)$ as $\nu\rightarrow\infty$ \textit{a.s.} $[\mathbb{P}]$.
\end{assumption}
\begin{assumption}\label{ass C3}
 $E_{\mathbb{P}}||Y_i||_{\mathcal{H}}^4 <\infty$,  $E_{\mathbb{P}}||Z_i||^4<\infty$,  and $E_{\mathbb{P}}(Z_i-E_{\mathbb{P}}(Z_i))^T(Z_i-E_{\mathbb{P}}(Z_i))$ is positive definite (p.d.) matrix.
\end{assumption}
Assumption \ref{ass C1} was considered earlier by several authors (cf. \cite{cardot2011horvitz}, \cite{conti2014estimation}, etc.).  The condition, $X_i\leq b$ \textit{a.s.} $[\mathbb{P}]$ for some $0<b<\infty$, in Assumption \ref{ass C2} and Assumption \ref{ass C1} along with $0\leq \lambda<E_{\mathbb{P}}(X_i) b^{-1}$ ensure that $n\max_{1\leq i\leq N} X_i \big(\sum_{i=1}^N X_i\big)^{-1}$ $<1$ for all sufficiently large $\nu$ \textit{a.s.} $[\mathbb{P}]$, which is required for implementing a $\pi$PS sampling design. In the case of any $\pi$PS sampling design, the condition that $\max_{1\leq i\leq N} X_i\big(\min_{1\leq i\leq N} X_i\big)^{-1}$ =$O(1)$ as $\nu\rightarrow\infty$ \textit{a.s.} $[\mathbb{P}]$, which appears in Assumption \ref{ass C2}, is equivalent to the condition that $M\leq n^{-1} N \pi_{i}  \leq M^{\prime}$ for some constants $M,M^{\prime}>0$, any $i$=$1,\ldots,N$ and all sufficiently large $\nu\geq 1$ \textit{a.s.} $[\mathbb{P}]$. This latter condition was considered earlier in \cite{MR3670194}, \cite{wang2011asymptotic} and references therein. Assumption \ref{ass C2} is used to prove some technical results (see Lemmas \ref{lem 5} and \ref{lem 6}--\ref{lem 2}) under LMS, HE$\pi$PS and RHC sampling designs, which will be required to show weak convergence of $\sqrt{n}(\hat{\overline{Y}}_{HT}-\overline{Y})$, $\sqrt{n}(\hat{\overline{Y}}_{RHC}-\overline{Y})$ and $\sqrt{n}(\hat{\overline{Y}}_{GREG}-\overline{Y})$ via uniform approximation (see \cite{kundu2000central}).  Assumption \ref{ass C3} implies that the fourth order raw moments of $Y_i$ and $Z_i$ exist. Now, we state the following proposition.
\begin{proposition}\label{prop 1}
Suppose that Assumptions \ref{ass C1}, \ref{ass C2} and \ref{ass C3} hold. Then, \textit{a.s.} $[\mathbb{P}]$, under SRSWOR and LMS sampling design, $\sqrt{n}(\hat{\overline{Y}}_{HT}-\overline{Y})\xrightarrow{\mathcal{L}}\mathcal{N}$ as $\nu\rightarrow\infty$, where $\mathcal{N}$ is a Gaussian distribution in $\mathcal{H}$ with mean $0$ and some covariance operator. Moreover, if Assumption \ref{ass C1} holds with  $0\leq \lambda<E_{\mathbb{P}}(X_i)b^{-1}$,  and Assumptions \ref{ass C2} and \ref{ass C3} hold, then the same result holds under any HE$\pi$PS sampling design.
\end{proposition}
 
 Next, recall $\{N_i\}_{i=1}^n$ in the context of RHC sampling design from  Section  \ref{sec 2}.  We choose $\{N_i\}_{i=1}^n$ to be as follows. 
\begin{equation}\label{eq 150}
\begin{split}
N_i=
\begin{cases}
Nn^{-1},\text{ for } i=1,\cdots,n, \text{ when } Nn^{-1} \text{ is an integer},\\
\lfloor Nn^{-1}\rfloor,\text{ for } i=1,\cdots,k, \text{ and }\\
\lfloor Nn^{-1} \rfloor+1,\text{ for } i=k+1,\cdots,n, \text{ when }Nn^{-1} \text{ is non-integer},\\
\end{cases}
\end{split}
\end{equation}
where $k$ is such that $\sum_{i=1}^n N_i$=$N$. Here, $\lfloor Nn^{-1} \rfloor$ is the integer part of $Nn^{-1}$. \cite{rao1962simple} showed that this choice of $\{N_i\}_{i=1}^n$ minimizes the variance of the RHC estimator. Now, we state the following propositions.  
\begin{proposition}\label{prop 3}
Suppose that Assumptions \ref{ass C1}, \ref{ass C2} and \ref{ass C3} hold. Then, \textit{a.s.} $[\mathbb{P}]$, under RHC sampling design, $\sqrt{n}(\hat{\overline{Y}}_{RHC}-\overline{Y})\xrightarrow{\mathcal{L}}\mathcal{N}$ as $\nu\rightarrow\infty$, where $\mathcal{N}$ is a Gaussian distribution in $\mathcal{H}$ with mean $0$ and some covariance operator.
\end{proposition}
\begin{proposition}\label{prop 2}
Suppose that Assumptions \ref{ass C1}, \ref{ass C2} and \ref{ass C3} hold. Then, \textit{a.s.} $[\mathbb{P}]$, under SRSWOR and LMS sampling design, $\sqrt{n}(\hat{\overline{Y}}_{GREG}-\overline{Y})\xrightarrow{\mathcal{L}}\mathcal{N}$ as $\nu\rightarrow\infty$, where $\mathcal{N}$ is a Gaussian distribution in $\mathcal{H}$ with mean $0$ and some covariance operator. Further, if Assumption \ref{ass C1} holds with  $0 \leq\lambda<E_{\mathbb{P}}(X_i)b^{-1}$,  and Assumptions \ref{ass C2} and \ref{ass C3} hold, then the same result holds under any HE$\pi$PS sampling design. Moreover, if Assumptions \ref{ass C1}, \ref{ass C2} and \ref{ass C3} hold, then the above result holds under RHC sampling design.
\end{proposition}
 The technique used to prove Propositions \ref{prop 1}--\ref{prop 2} is based on the idea of convergence in distribution via uniform approximation considered in \cite{kundu2000central}. This idea was used in \cite{kundu2000central} to extend central limit theorem for independent random variables from finite dimensional Euclidean space to an infinite dimensional separable Hilbert space (see Proposition $2.1$ in \cite{kundu2000central}). Any infinite dimensional separable Hilbert space (e.g., the space of square integrable functions equipped with $L^2$-inner product) is isometrically isomorphic to the space of square summable sequences $l^2$ because a separable Hilbert space always has a complete orthonormal basis. Further, the $l^2$ space can be conveniently viewed as an infinite dimensional extension of a finite dimensional Euclidean space. Thus it is relatively easy to extend the results from multivariate data setup to the functional data setup using the sequence structure of the $l^2$ space.
\par

 \cite{cardot2011horvitz} and \cite{cardot2013uniform} showed the weak convergence of the HT and the model assisted estimators, respectively, in $\mathcal{C}[0,T]$ under some conditions on sampling designs (see pp. $110$-$111$ in \cite{cardot2011horvitz} and  pp. $569$-$573$ in \cite{cardot2013uniform}). These conditions hold under usual sampling designs like SRSWOR, stratified sampling design with SRSWOR, rejective sampling design, etc. We are able to dispense with these conditions, and show the weak convergence of the HT and the GREG estimators in a separable Hilbert space under SRSWOR, LMS sampling design and any HE$\pi$PS sampling design. We are also able to show the weak convergence of the RHC and the GREG estimators in a separable Hilbert space under RHC sampling design.
\par
\vspace{.1cm}

  We develop our results in a separable Hilbert space framework rather than in a space of continuous functions equipped with supremum norm because we are able to prove Propositions \ref{prop 1}--\ref{prop 2} in the case of a separable Hilbert space framework. The space of continuous functions is a subset of the space of square integrable functions, which is a separable Hilbert space equipped with $L^2$ inner product. Random functions from the space of continuous functions can be expressed as linear combinations of orthonormal basis functions in the space of square integrable functions through the Karhunen-Lo\`eve expansion. 
\par
\vspace{.1cm}

 Next, we carry out the comparison of the estimators mentioned earlier based on the above results. We say that an estimator $\hat{\overline{Y}}_1$ with asymptotic covariance operator $\Sigma$ is asymptotically at least as efficient as another estimator $\hat{\overline{Y}}_2$ with asymptotic covariance operator $\Delta$ if $\Delta-\Sigma$ is non negative definite (n.n.d.),  i.e., if $\langle(\Delta-\Sigma)a,a\rangle \geq 0$ for any $a\in\mathcal{H}$.  We also say that $\hat{\overline{Y}}_1$ is asymptotically more efficient than $\hat{\overline{Y}}_2$ if $\Delta-\Sigma$ is p.d,  i.e., if $\langle(\Delta-\Sigma)a,a\rangle > 0$ for any $a\in\mathcal{H}$ and $a\neq 0$.  We now state the following theorems.
\begin{theorem}\label{thm 2} 
Suppose that Assumptions \ref{ass C1}, \ref{ass C2} and \ref{ass C3} hold. Then, \textit{a.s.} $[\mathbb{P}]$, the GREG estimator is asymptotically at least as efficient as the HT estimator under SRSWOR as well as LMS sampling design. Moreover, \textit{a.s.} $[\mathbb{P}]$, both the HT and the GREG estimators have the same asymptotic distribution under SRSWOR and LMS sampling design. 
\end{theorem}
Before we state the next theorem, let us consider superpopulations satisfying linear regression model 
\begin{equation}\label{eq 2}
Y_i=\beta_0+\sum_{j=1}^d Z_{ji}\beta_j+\epsilon_i X_i^\eta,
\end{equation}
where $i$=$1,\ldots,N$, $\{\epsilon_i\}_{i=1}^N$ are i.i.d. $\mathcal{H}$ valued random variables independent of $\{Z_i,X_i\}_{i=1}^N$ with mean $0$. Here, $Z_i$=$(Z_{1i},\ldots,Z_{di})$, $\beta_j\in\mathcal{H}$ for $j$=$0,\ldots,d$, and $\eta\geq 0$  is the degree of heteroscedasticity present in the linear model given above. For any given $\eta>0$, the conditional total variance of $Y_i$ given $(Z_i,X_i)$, the trace of the conditional covariance operator of $Y_i$ given $(Z_i,X_i)$, increases as the value of $X_i$ increases (cf. \cite{sarndal1992model}). In essence, the parameter $\eta$ determines the rate at which this conditional total variance increases with $X_i$. Similar types of linear model as in \eqref{eq 2} were used for constructing several estimators by earlier authors, when the observations on $y$ are from some finite dimensional Euclidean space (see \cite{cassel1976some}, \cite{sarndal1980pi}, \cite{sarndal1992model} and references therein). A homoscedastic (i.e., when $\eta$=$0$) version of the above linear regression model was considered earlier in \cite{cardot2013comparison} and \cite{cardot2013uniform}   for constructing the model assisted estimator of $\overline{Y}$, when the observations on $y$ are from some functional space. Now, we state the following theorem.
\begin{theorem}\label{thm 1} 
Suppose that \eqref{eq 2} and Assumptions \ref{ass C1}, \ref{ass C2} and \ref{ass C3} hold. Then, \textit{a.s.} $[\mathbb{P}]$, the GREG estimator is asymptotically at least as efficient as the RHC estimator under RHC sampling design. Further, if \eqref{eq 2} holds, Assumption \ref{ass C1} holds with  $0 \leq\lambda<E_{\mathbb{P}}(X_i)b^{-1}$,  and Assumptions \ref{ass C2} and \ref{ass C3} hold, then \textit{a.s.} $[\mathbb{P}]$, the GREG estimator is asymptotically at least as efficient as the HT estimator under any HE$\pi$PS sampling design.   
\end{theorem}
 It follows from the preceding results that the GREG estimator is asymptotically at least as efficient as the HT and RHC estimators under each of the sampling designs considered in this paper. Also, both the HT and the GREG estimators have the same asymptotic distribution under SRSWOR and LMS sampling design. Now, we compare the performance of the GREG estimator under SRSWOR, RHC sampling design and HE$\pi$PS sampling designs based on the degree of heteroscedasticity $\eta$.  
\begin{theorem}\label{thm 3}
Suppose that \eqref{eq 2} holds, and $\epsilon_i$ has a p.d. covariance operator. Further, suppose that Assumption \ref{ass C1} holds with $0 \leq\lambda<E_{\mathbb{P}}(X_i)b^{-1}$, and Assumptions \ref{ass C2} and \ref{ass C3} hold. Then, the sampling designs among SRSWOR, HE$\pi$PS and RHC sampling designs under which the GREG estimator becomes the most efficient estimator \textit{a.s.} $[\mathbb{P}]$ are as mentioned in Table \ref{table 1} below. Further, if Assumption \ref{ass C1} holds with $\lambda$=$0$, and Assumptions \ref{ass C2} and \ref{ass C3} hold then the GREG estimator has the same asymptotic distribution under RHC and any HE$\pi$PS sampling designs.  
\begin{table}[h]

\centering
\caption{Sampling designs for which the GREG estimator becomes the most efficient estimator}
\label{table 1}
\begin{threeparttable}[b]
\begin{tabular}{cccc} 
\hline
& \multirow{2}{*}{$\lambda$=$0$}& $\lambda>0$ \& & $\lambda>0$ \&  \\ 
&&$\lambda^{-1}$ is an integer& $\lambda^{-1}$ is a non-integer\\
\hline 
$\eta<0.5$ & SRSWOR& SRSWOR& SRSWOR\\
\hline 
$\eta=0.5$& \tnote{*} \ SRSWOR, HE$\pi$PS \& RHC& \tnote{*} \ SRSWOR, HE$\pi$PS \& RHC& \tnote{**}\ \ \ SRSWOR \& HE$\pi$PS\\
\hline 
$\eta>0.5$& \tnote{***}\ \ \ \ \ HE$\pi$PS \& RHC& HE$\pi$PS& HE$\pi$PS\\
\hline
\end{tabular}
\begin{tablenotes}
\item[*] GREG estimator has the same asymptotic distribution under SRSWOR,   RHC sampling design and HE$\pi$PS sampling designs for $\eta$=$0.5$, $\lambda>0$ and $\lambda^{-1}$ being an integer.
\vspace{.2cm}
 
\item[**] GREG estimator has the same asymptotic distribution under SRSWOR and HE$\pi$PS sampling designs for $\eta$=$0.5$, $\lambda>0$ and $\lambda^{-1}$ being a non-integer.
\vspace{.2cm}
 
\item[***] GREG estimator has the same asymptotic distribution under HE$\pi$PS and RHC sampling designs for $\eta>0.5$ and $\lambda$=$0$.
\end{tablenotes}
\end{threeparttable}
\end{table}
\end{theorem}
\par
\vspace{.2cm}

Proofs of Theorems \ref{thm 2}--\ref{thm 3} involve some results related to operator theory, which are available in \cite{hsing2015theoretical}. It follows from \eqref{eq 52} in the proof of Theorem \ref{thm 3} that $cov_{\mathbb{P}}(X_i^{2\eta-1},X_i)$, the covariance between $X_i^{2\eta-1}$ and $X_i$, determines the sampling design among SRSWOR, HE$\pi$PS and RHC sampling designs under which the GREG estimator becomes the most efficient estimator. The GREG estimator performs more efficiently under SRSWOR than under RHC and any HE$\pi$PS sampling designs, whenever $cov_{\mathbb{P}}(X_i^{2\eta-1},X_i)<0$. On the other hand, the GREG estimator under RHC as well as any HE$\pi$PS sampling design becomes more efficient than the GREG estimator under SRSWOR in the case of $\lambda$=$0$, whenever $cov_{\mathbb{P}}(X_i^{2\eta-1},X_i)>0$, and the GREG estimator under any HE$\pi$PS sampling design becomes more efficient than the GREG estimator under both SRSWOR and RHC sampling design in the case of $\lambda>0$, whenever $cov_{\mathbb{P}}(X_i^{2\eta-1},X_i)>0$. Now, $x^{2\eta-1}$ is a decreasing function of $x$ for $\eta<0.5$ and an increasing function of $x$ for $\eta>0.5$. Therefore, $cov_{\mathbb{P}}(X_i^{2\eta-1},X_i)<0$ for $\eta<0.5$ and $cov_{\mathbb{P}}(X_i^{2\eta-1},X_i)>0$ for $\eta>0.5$. Thus the use of the auxiliary information in HE$\pi$PS and RHC sampling designs has an adverse effect on the performance of the GREG estimator, when $\eta<0.5$. On the other hand, for the case of $\eta>0.5$, the use of HE$\pi$PS and RHC sampling designs improves the performance of the GREG estimator. 
\par
\vspace{.2cm}

Note that if we consider a generalized version of the linear regression model in \eqref{eq 2} as $Y_i$=$\beta_0+\sum_{j=1}^d Z_{ji}\beta_j+\epsilon_i g(X_i)$ for $i$=$1,\ldots,N$ and some non-negative real valued function $g$, then it can be shown in the same way as in the proof of Theorem \ref{thm 1} that the conclusion of Theorem \ref{thm 1} holds under the above linear model. It can also be shown in the same way as in the proof of Theorem \ref{thm 3} that the results in $2^{nd}$, $3^{rd}$ and $4^{th}$ rows in Table \ref{table 1} related to Theorem \ref{thm 3} hold, whenever $cov_{\mathbb{P}}( g^2(X_i)X_i^{-1},X_i)<0$, $cov_{\mathbb{P}}( g^2(X_i)X_i^{-1},X_i)$=$0$ and $cov_{\mathbb{P}}( g^2(X_i)X_i^{-1},X_i)>0$, respectively. In particular, the results in $2^{nd}$, $3^{rd}$ and $4^{th}$ rows in Table \ref{table 1} hold if $g^2(x)x^{-1}$ is decreasing, constant and increasing function of $x$, respectively. 
\par
\vspace{.2cm}

 Let us denote the asymptotic covariance operator of $\sqrt{n}(\hat{\overline{Y}}-\overline{Y})$ by $\Sigma$, where $\hat{\overline{Y}}$ denotes one of $\hat{\overline{Y}}_{HT}$, $\hat{\overline{Y}}_{RHC}$ and $\hat{\overline{Y}}_{GREG}$. Next, suppose that $\hat{\overline{Y}}$ is either $\hat{\overline{Y}}_{HT}$ or $\hat{\overline{Y}}_{GREG}$ under one of SRSWOR, LMS sampling design and any HE$\pi$PS sampling design. Then, it follows from the proofs of Propositions \ref{prop 1} and \ref{prop 2} that $\Sigma$=$\lim_{\nu\rightarrow\infty} nN^{-2}\sum_{i=1}^N(V_i-T\pi_i)\otimes(V_i-T\pi_i)(\pi_i^{-1}-1)$ \textit{a.s.} $[\mathbb{P}]$, where  $T$=$\sum_{i=1}^N V_i(1-\pi_i)\big(\sum_{i=1}^N \pi_i(1-\pi_i)\big)^{-1}$  and $\{\pi_i\}_{i=1}^N$ are inclusion probabilities. Further, $V_i$ is $Y_i$ for $\hat{\overline{Y}}$ being $\hat{\overline{Y}}_{HT}$. Also, $V_i$ is $Y_i-\overline{Y}-S_{zy}((Z_i-\overline{Z})S_{zz}^{-1})$ for $\hat{\overline{Y}}$ being $\hat{\overline{Y}}_{GREG}$. Here,  $S_{zy}$=$N^{-1}\sum_{i=1}^N(Z_i-\overline{Z})\otimes(Y_i-\overline{Y})$ and $S_{zz}$=$N^{-1}\sum_{i=1}^N(Z_i-\overline{Z})^T(Z_i-\overline{Z})$. We estimate $\Sigma$ by 
\begin{equation}\label{eq 33}
\hat{\Sigma}=(nN^{-2})\sum_{i\in s}(\hat{V}_i-\hat{T}\pi_i)\otimes(\hat{V}_i-\hat{T}\pi_i)(\pi_i^{-1}-1)\pi_i^{-1},
\end{equation}
 where $\hat{V}_i$ is $Y_i$ or $Y_i-\hat{\overline{Y}}_{HT}-\hat{S}_{zy}((Z_i-\hat{\overline{Z}}_{HT})\hat{S}_{zz}^{-1})$ for $\hat{\overline{Y}}$ being $\hat{\overline{Y}}_{HT}$ or $\hat{\overline{Y}}_{GREG}$, respectively. Also,  $\hat{T}$=$\sum_{i\in s}\hat{V}_i(\pi_i^{-1}-1)\big(\sum_{i\in s}(1-\pi_i)\big)^{-1}$, $\hat{S}_{zz}$=$\big(\sum_{i\in s}\pi_i^{-1}\big)^{-1}\sum_{i\in s}\pi_i^{-1}(Z_i-\hat{\overline{Z}})^T(Z_i-\hat{\overline{Z}})$, and $\hat{S}_{zy}$=$\big(\sum_{i\in s}\pi_i^{-1}\big)^{-1}\sum_{i\in s}\pi_i^{-1}(Z_i-\hat{\overline{Z}})\otimes(Y_i-\hat{\overline{Y}})$.
\par

 Next, suppose that $\hat{\overline{Y}}$ is either $\hat{\overline{Y}}_{RHC}$ or $\hat{\overline{Y}}_{GREG}$ under RHC sampling design. Then, it can be shown from the proofs of Propositions \ref{prop 3} and \ref{prop 2} that  $\Sigma$=$\lim_{\nu\rightarrow\infty} n\gamma\overline{X}N^{-1}\sum_{i=1}^N\big(V_i-X_i\overline{V} \overline{X}^{-1})\otimes(V_i-X_i\overline{V} \overline{X}^{-1}\big) X_i^{-1}$ \textit{a.s.} $[\mathbb{P}]$, where $\gamma$=$\sum_{i=1}^n N_i(N_i-1)(N(N-1))^{-1}$  with $N_i$ being the size of the $i^{th}$ group formed randomly in the first step of the RHC sampling design (see  Section  \ref{sec 2}), $i$=$1,\ldots,n$. Further, $V_i$ is $Y_i$ for $\hat{\overline{Y}}$ being $\hat{\overline{Y}}_{RHC}$. Also, $V_i$ is $Y_i-\overline{Y}-S_{zy}((Z_i-\overline{Z})S_{zz}^{-1})$ for $\hat{\overline{Y}}$ being $\hat{\overline{Y}}_{GREG}$. In this case, we estimate $\Sigma$ by  
\begin{equation}\label{eq 34}
\hat{\Sigma}=n\gamma(\overline{X}N^{-1})\sum_{i\in s}\bigg(\hat{V}_i-X_i\hat{\overline{V}}_{RHC}\overline{X}^{-1}\bigg)\otimes\bigg(\hat{V}_i-X_i\hat{\overline{V}}_{RHC} \overline{X}^{-1}\bigg)(Q_iX_i^{-2}),
\end{equation}
 where $\hat{V}_i$ is $Y_i$ or $Y_i-\hat{\overline{Y}}_{RHC}-\hat{S}_{zy}((Z_i-\hat{\overline{Z}}_{RHC})\hat{S}_{zz}^{-1})$ for $\hat{\overline{Y}}$ being $\hat{\overline{Y}}_{RHC}$ or $\hat{\overline{Y}}_{GREG}$, respectively. Further,  $\hat{\overline{V}}_{RHC}$=$\sum_{i\in s}(NX_i)^{-1} \hat{V}_i Q_i$, $\hat{\overline{Z}}_{RHC}$=$\sum_{i\in s}(NX_i)^{-1}Z_i Q_i $ and $\hat{S}_{zy}$ and $\hat{S}_{zz}$ are the same as above with $\pi_i^{-1}$ replaced by $Q_i X_i^{-1}$. Also, recall $b$ from Assumption \ref{ass C2}.  Now, we state the following theorem concerning the consistency of $\hat{\Sigma}$ as an estimator of $\Sigma$ with respect to the HS norm (see \cite{hsing2015theoretical}). 
\begin{theorem}\label{prop 4}
Let us consider $\Sigma$, the asymptotic covariance operator of $\sqrt{n}(\hat{\overline{Y}}-\overline{Y})$, and its estimator $\hat{\Sigma}$ from the preceding  discussion.  Suppose that Assumptions \ref{ass C1}, \ref{ass C2} and \ref{ass C3} hold. Then, a.s. $[\mathbb{P}]$, under SRSWOR and LMS sampling design, $\hat{\Sigma}\xrightarrow{p}\Sigma$ as $\nu\rightarrow\infty$.  Here, the convergence in probability holds with respect to the HS norm.  Further, if Assumption \ref{ass C1} holds with  $0\leq \lambda<E_{\mathbb{P}}(X_i)b^{-1}$,  and Assumptions \ref{ass C2} and \ref{ass C3} hold, then the same result holds under any HE$\pi$PS sampling design. Moreover, if Assumptions \ref{ass C1}, \ref{ass C2} and \ref{ass C3} hold, then the above result holds under RHC sampling design. 
\end{theorem}

\section{Data analysis}\label{sec 6}
\subsection{Analysis based on synthetic data}\label{subsec 6.1} In this section, we consider a finite population of size $N$=$1000$ generated as follows. We first generate the observations $X_1,\ldots,X_N$ on the size variable $x$ from a gamma distribution with mean $500$ and standard deviation (s.d.) $100$. Here, we assume that the covariate $z$ and the size variable $x$ are same. Then, we generate the population observations on $y$ from $L^2[0,1]$ using linear regression models $Y_{i}(t)$=$1000+ \beta(t) X_i+\epsilon_{i}(t)X_i^\eta$, where $\beta(t)$=$1$, $t$ and $1-(t-0.5)^2$, $\eta$=$0.1 k$ for $k$=$0,1,\ldots,10$, and $\{\epsilon_i(t)\}_{t\in[0,1]}$'s are i.i.d. copies of standard Brownian motion with mean $0$ and covariance kernel $\sigma(s,t)$=$s\wedge t$. The population observations on $y$ are generated at $t_1,\ldots,  t_r $, where $ r $=$100$ and $t_j$=$j  r^{-1} $ for $j$=$1,\ldots,  r $. We now consider the estimation of the mean of $y$. We compare the HT and the GREG estimators under SRSWOR and RS sampling design, and the RHC and the GREG estimators under RHC sampling design in terms of relative efficiencies as defined in the following paragraph. The RS sampling design is chosen as a HE$\pi$PS sampling design since it is easier to implement than any other HE$\pi$PS sampling design. We shall not report the results under LMS sampling design because these results are very close to the results under SRSWOR as expected from our theoretical results. 
\par

Suppose that each curve in a population of $N$ curves from $L^2[0,1]$ is observed at $t_1,\ldots, t_r  \in [0,1]$ for some $ r  >1$. Let us consider $I$ samples each of size $n$ from this population. Then, the mean squared error of an estimator of $\overline{Y}$, say $\hat{\overline{Y}}$, under sampling design $P(s)$ is computed as  $MSE(\hat{\overline{Y}},P)$=$(rI)^{-1}\sum_{l=1}^I\sum_{j=1}^r (\hat{\overline{Y}}_{l}(t_j)-\overline{Y}(t_j))^2$  (see \cite{cardot2011horvitz}, \cite{cardot2013comparison}, etc.), where $\hat{\overline{Y}}_{l}$ is an estimate of $\overline{Y}$ based on the $l^{th}$ sample, $l$=$1,\ldots,I$. Further, we define the relative efficiency of an estimator $\hat{\overline{Y}}_1$ under sampling design $P_1(s)$ compared to another estimator $\hat{\overline{Y}}_2$ under sampling design $P_2(s)$ by 
$$
RE(\hat{\overline{Y}}_1,P_1|\hat{\overline{Y}}_2, P_2)=MSE(\hat{\overline{Y}}_2,P_2)\big(MSE(\hat{\overline{Y}}_1,P_1)\big)^{-1}.
$$
 We say that $\hat{\overline{Y}}_1$ under $P_1(s)$ is more efficient than $\hat{\overline{Y}}_2$ under $P_2(s)$ if RE($\hat{\overline{Y}}_1$, $P_1$ $|$ $\hat{\overline{Y}}_2$, $P_2$)$>1$. We compute relative efficiencies of the estimators mentioned in the preceding paragraph based on $I$=$1000$ samples each of size $n$=$100$. We plot the relative efficiency of the HT estimator compared to the GREG estimator under each of SRSWOR and RS sampling design as well as the relative efficiency of the RHC estimator compared to the GREG estimator under RHC sampling design for different $\eta$. We also plot the relative efficiency of the GREG estimator under SRSWOR compared to the GREG estimator under each of RS and RHC sampling designs. We use the $R$ software for drawing samples as well as computing estimators. For RS sampling design, we use the `pps' package in $R$.  The $R$ codes for this simulation study are  available at \href{https://github.com/deyanurag/Simulation-code.git}{{https://github.com/deyanurag/Simulation-code.git}}.  The results obtained from the above analyses are summarized as follows. 
\begin{figure}[h]
\includegraphics[height=8cm,width=13cm]{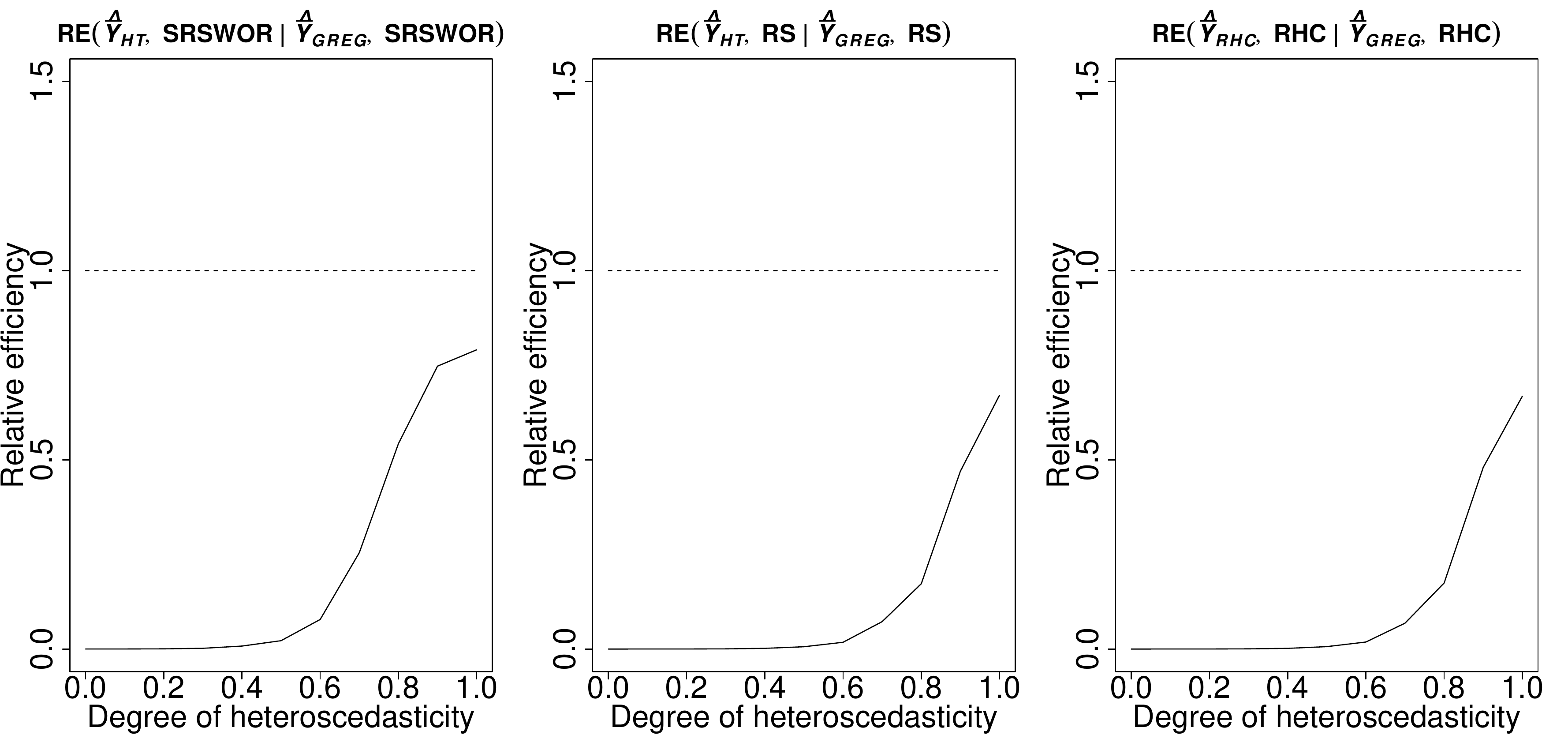}
\centering
\caption{Comparison of HT, GREG and RHC estimators under different sampling designs for $\beta(t)$=$1$}
\label{Fig 5}
\end{figure}

\begin{figure}[h!]

\includegraphics[height=8cm,width=13cm]{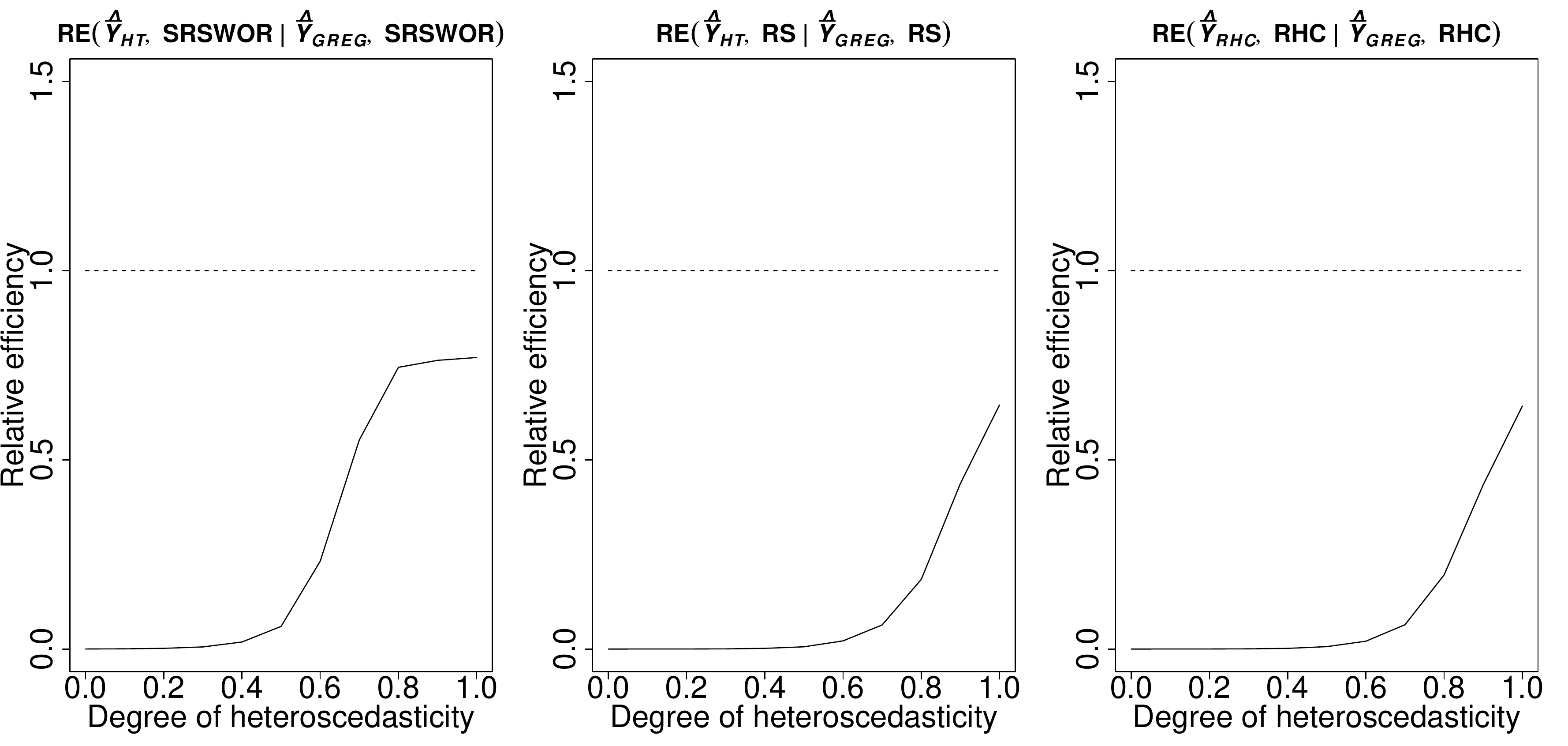}
\centering
\caption{Comparison of HT, GREG and RHC estimators under different sampling designs for $\beta(t)$=$t$.}
\label{Fig 1}
\end{figure}

\begin{figure}[h!]

\includegraphics[height=8cm,width=13cm]{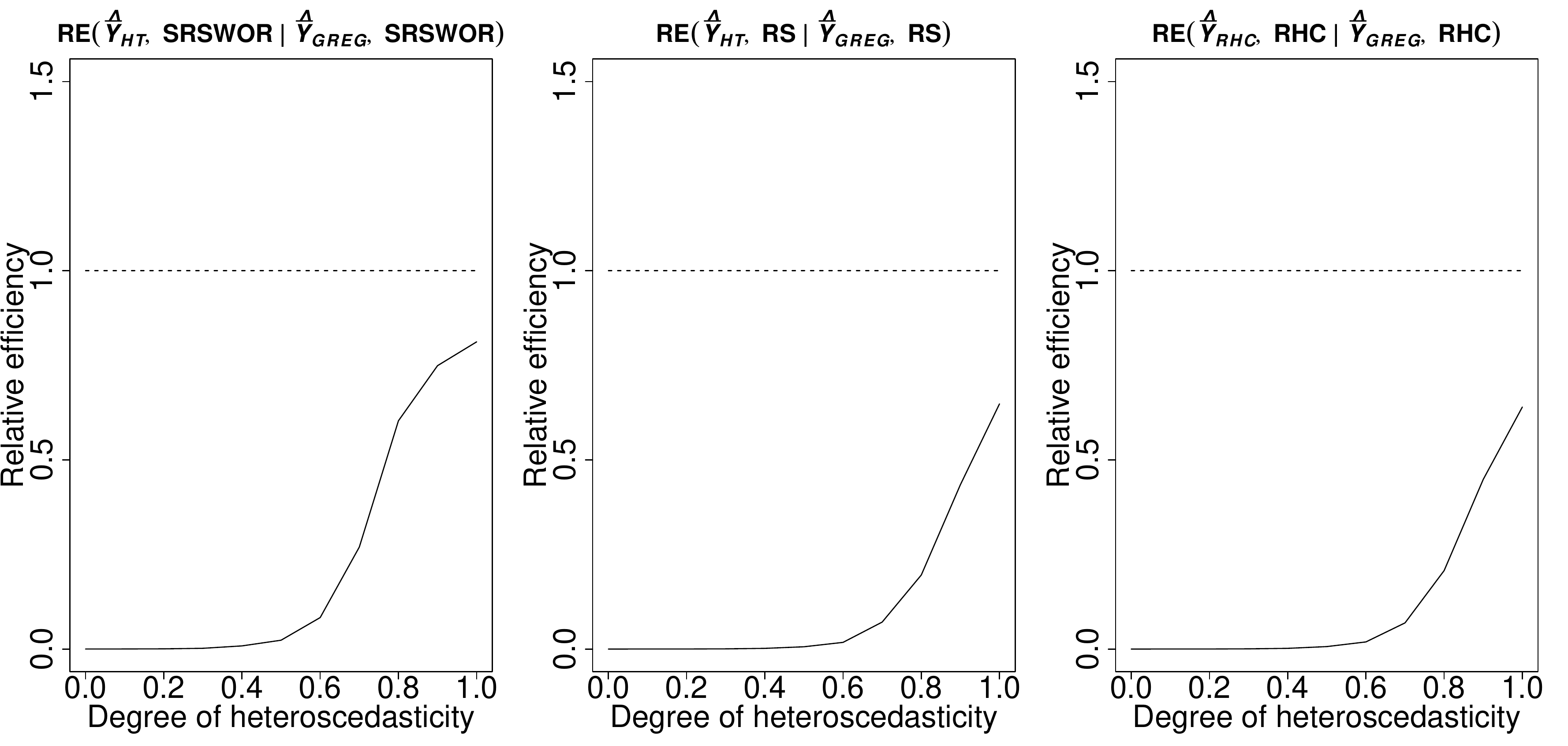}
\centering
\caption{Comparison of HT, GREG and RHC estimators under different sampling designs for $\beta(t)$=$1-(t-0.5)^2$}
\label{Fig 3}
\end{figure}

\begin{figure}[h!]
\includegraphics[height=8cm,width=13cm]{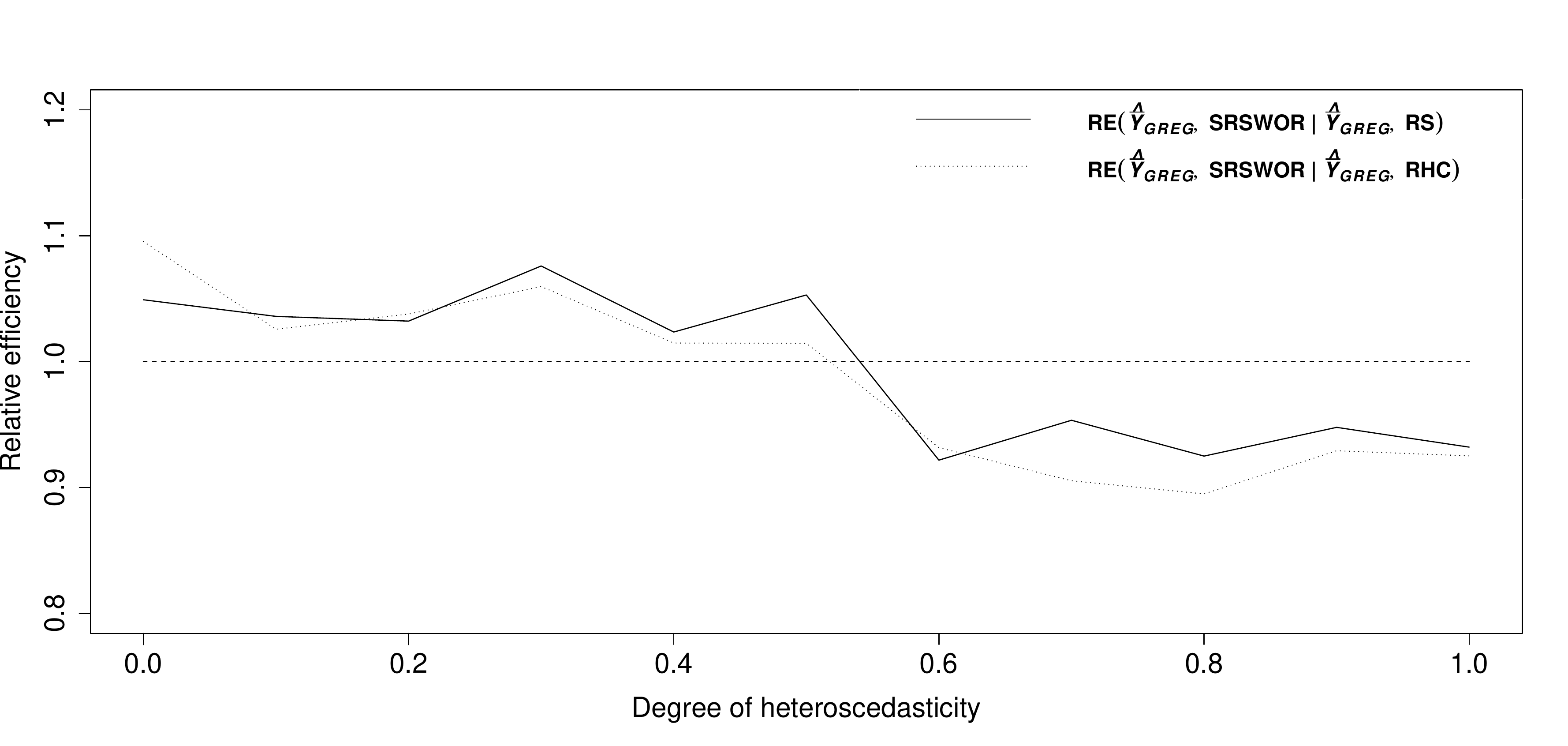}
\centering
\caption{Comparison of GREG estimators under different sampling designs for $\beta(t)$=$1$}
\label{Fig 6}
\end{figure}

\begin{figure}[h!]
\includegraphics[height=8cm,width=13cm]{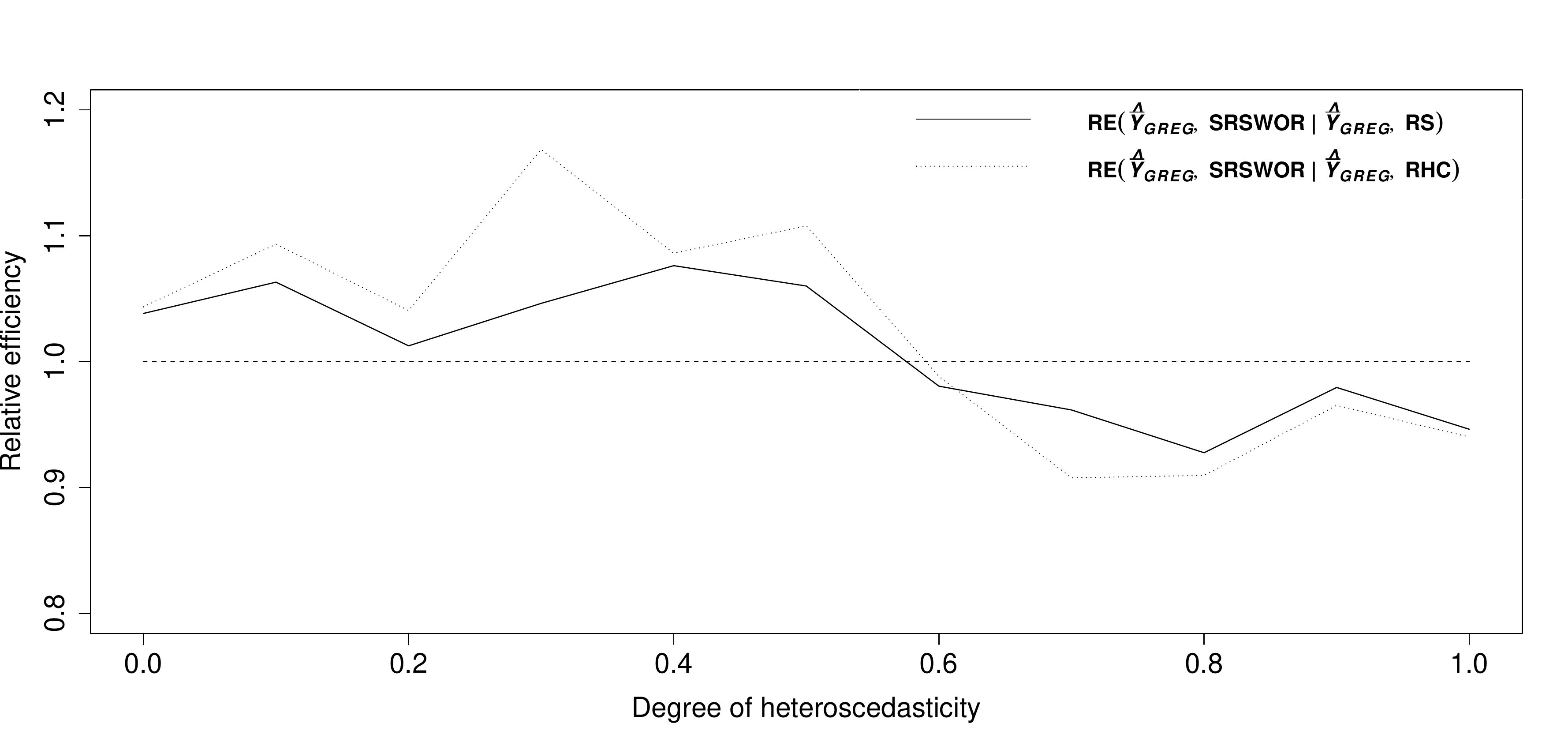}
\centering
\caption{Comparison of GREG estimators under different sampling designs for $\beta(t)$=$t$}
\label{Fig 2}
\end{figure}

\begin{figure}[h!]
\includegraphics[height=8cm,width=13cm]{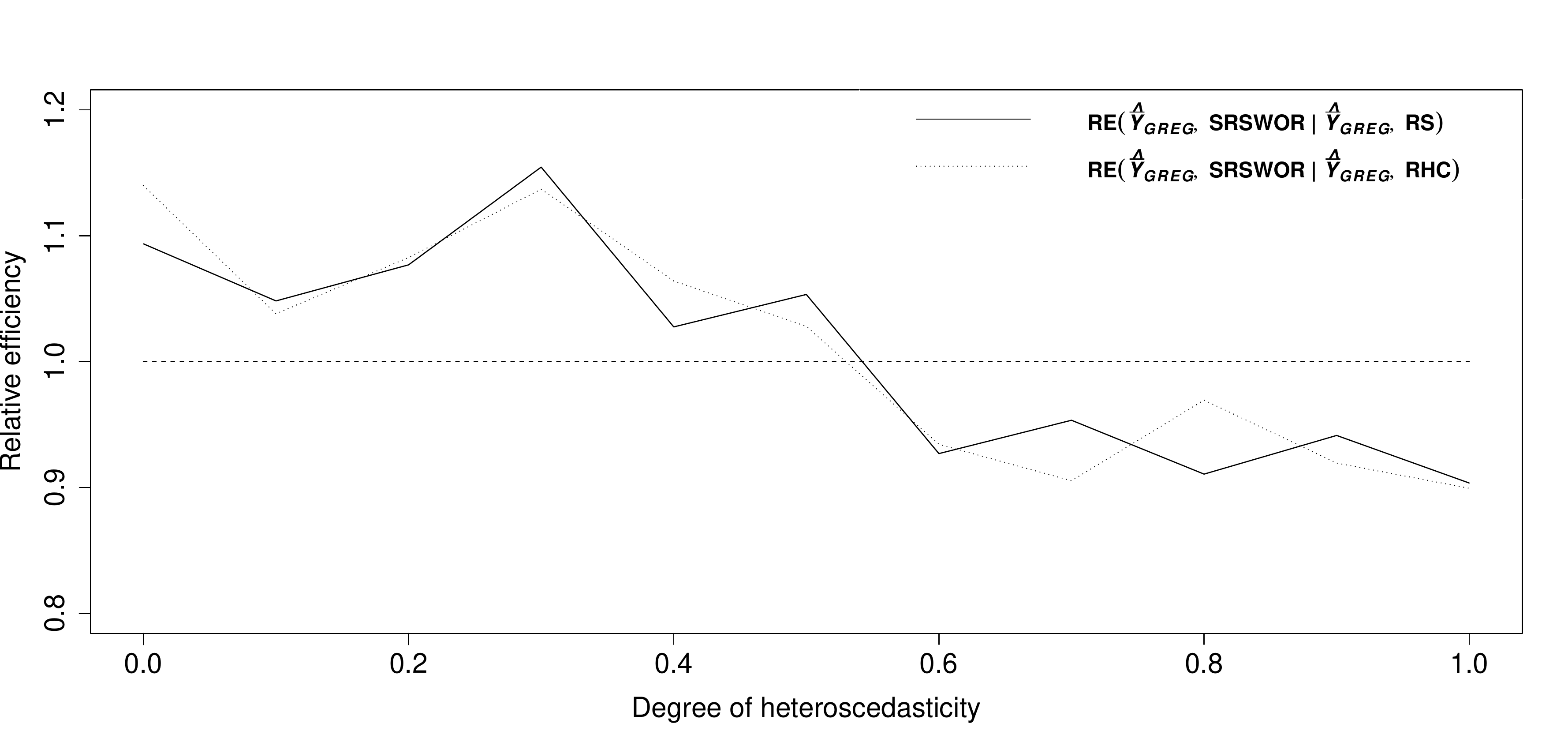}
\centering
\caption{Comparison of GREG estimators under different sampling designs for $\beta(t)$=$1-(t-0.5)^2$}
\label{Fig 4}
\end{figure}

\begin{itemize}
\item[(i)] It follows from Figures \ref{Fig 5}, \ref{Fig 1} and \ref{Fig 3} that the relative efficiency curve of the HT estimator compared to the GREG estimator under each of SRSWOR and RS sampling design and that of the RHC estimator compared to the GREG estimator under RHC sampling design always lie below the $y=1$ line (dashed line), when $\beta(t)$=$1,t$ or $1-(t-0.5)^2$. This implies that the GREG estimator is more efficient than the HT estimator under SRSWOR and RS sampling design, and the GREG estimator is more efficient than the RHC estimator under RHC sampling design for different $\eta$. The above results are in conformity with Theorems \ref{thm 2} and \ref{thm 1}. 
\item[(ii)] We see from Figures \ref{Fig 6}, \ref{Fig 2} and \ref{Fig 4} that the relative efficiency curve of the GREG estimator under SRSWOR compared to that under each of RS and RHC sampling designs lies above $y=1$ line, when $\eta<0.5$ and $\beta(t)$=$1,t$ or $1-(t-0.5)^2$. However, these lines lie below $y=1$ line, when $\eta>0.5$. This means that the use of the sampling designs like RS and RHC have an adverse effect on the performance of the GREG estimator, when $\eta< 0.5$. However, the use of the above sampling designs improves the performance of the GREG estimator, when $\eta>0.5$. Thus the above empirical results corroborate the theoretical results stated in Theorem \ref{thm 3}.
\end{itemize}

\subsection{Analysis based on real data}\label{subsec 6.3}
In this section, we consider Electricity Customer Behaviour Trial data available in Irish Social Science Data Archive (ISSDA, \href{https://www.ucd.ie/issda/}{{https://www.ucd.ie/issda/}}). In this dataset, we have electricity consumption of Irish households measured (in kWh) at the end of every half an hour during the period, $14^{th}$ July in $2009$ to $31^{st}$ December in $2010$. We are interested in the estimation of the mean electricity consumption curve in the summer months, viz. June, July and August in $2010$ and in the winter month of December in $2010$. It is to be noted that we consider the estimation of the mean electricity consumption curve only in the winter month of December in $2010$ because the data for the other two months in the winter of $2010$, viz. January and February in $2011$ are unavailable. In this dataset, we have $N$=$5372$ households for which electricity consumption data are available during July and August of $2009$ and all the summer months of $2010$. We also have $N$=$5092$ households for which electricity consumption data are available during December of both $2009$ and $2010$. Further, for each unit, there are $4416$ and $1488$ measurement points in summer months and December of $2010$, respectively. Electricity consumption in summer months and December of $2010$ can be viewed as electricity consumption curves in $L^2[0,T_1]$ and $L^2[0,T_2]$, respectively, where $T_1$=$30\times 4416$=$132480$ and $T_2$=$30\times 1488$=$44640$. For estimating the mean electricity consumption curve in the summer months of $2010$, we choose the mean electricity consumption in July and August of $2009$ as the size variable $x$, the mean electricity consumption in July of $2009$ as the first covariate $z_1$ and the mean electricity consumption in August of $2009$ as the second covariate $z_2$. On the other hand, for estimating the mean electricity consumption curve in December of $2010$, we choose the mean electricity consumption in December of $2009$ as both the size variable $x$ and the covariate $z$. In case of the above estimation problems, we compare the estimators considered in the preceding section in terms of relative efficiencies (see  Section  \ref{subsec 6.1}).  We compute relative efficiencies of these estimators based on $I$=$1000$ samples each of size $n$=$100$, where these samples are selected from the two datasets consisting of $5372$ and $5092$ observations, respectively.  The results obtained from the above analyses are summarized as follows. 
\par
\vspace{.1cm}

\begin{figure}[h]
\includegraphics[height=8cm,width=13cm]{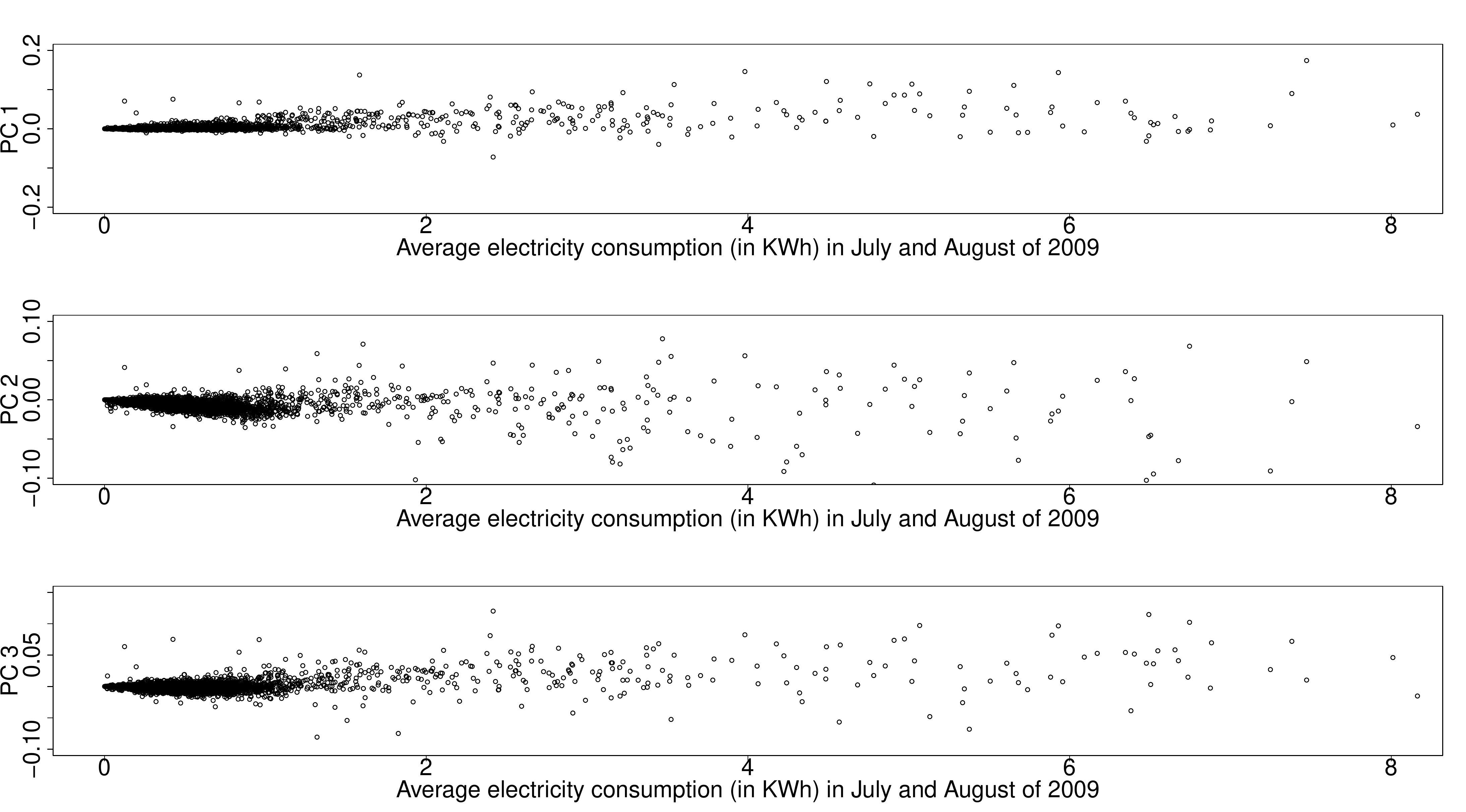}
\centering
\caption{Scatter plots  of the first three principal components  of electricity consumption data versus the size variable}
\label{Fig 7}
\end{figure}

\begin{figure}[h]
\includegraphics[height=8cm,width=13cm]{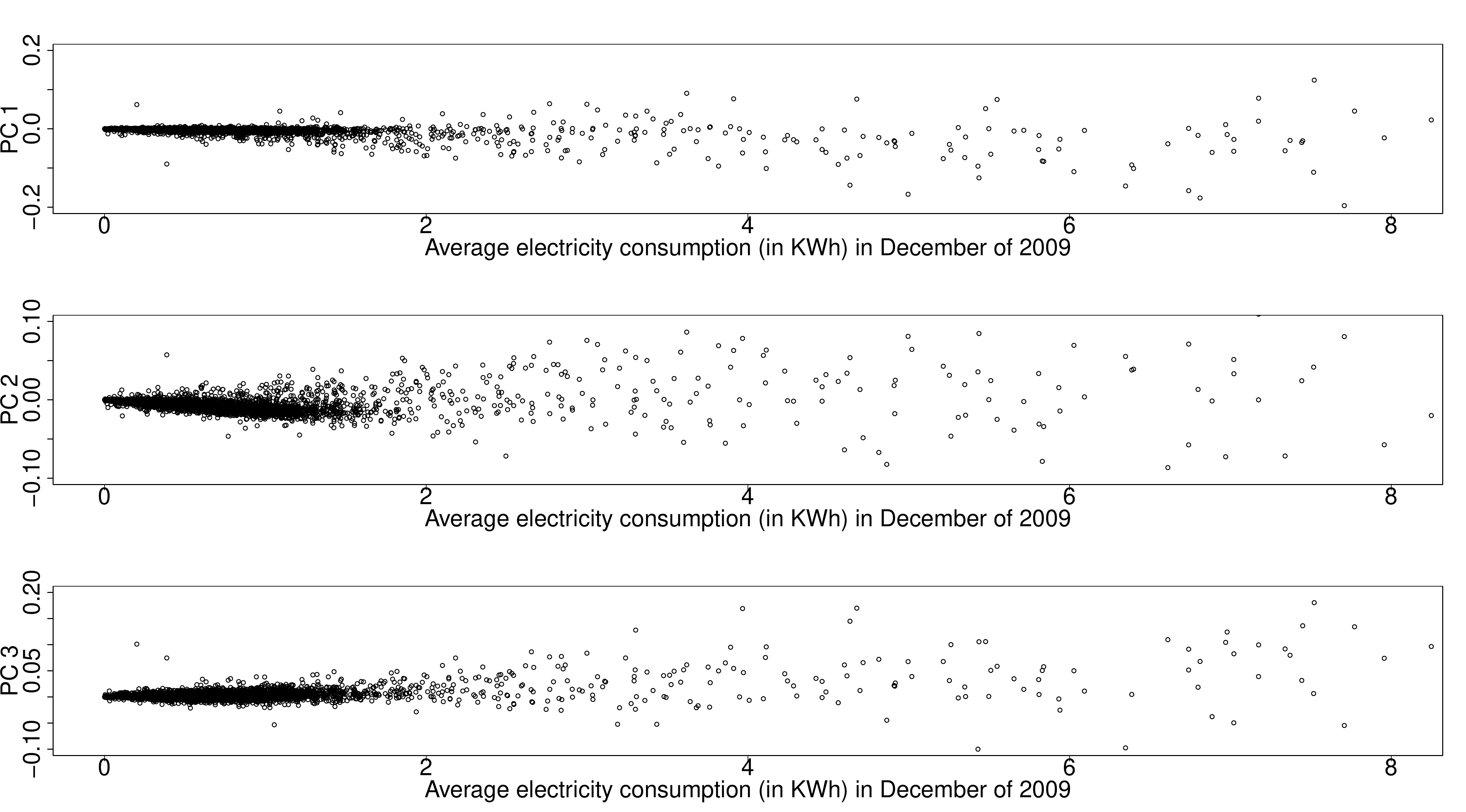}
\centering
\caption{Scatter plots of  the first three principal components  of electricity consumption data versus the size variable}
\label{Fig 8}
\end{figure}

\begin{itemize}
\item[(i)] We see from Table \ref{table 9} that the GREG estimator is more efficient than the HT estimator under SRSWOR and RS sampling design in both  the datasets.  Also, the GREG estimator is more efficient than the RHC estimator under RHC sampling design in both  the datasets.  Therefore, these results support the results stated in Theorems \ref{thm 2} and \ref{thm 1}.

\item[(ii)] In the case of both the  datasets,  we observe the presence of substantial heteroscedasticity in electricity consumption data, when we plot each of the first three principal components (PC) of electricity consumption data against the size variable (see Figures \ref{Fig 7} and \ref{Fig 8}). Further, it follows from Table \ref{table 8} that the GREG estimator under RS sampling design is more efficient than any other estimator under any other sampling design for both the  datasets.  Thus the empirical results stated here are in conformity with the theoretical results stated in Theorem \ref{thm 3}. 
\end{itemize}
\begin{table}[h]
\centering
\caption{Relative efficiencies of the HT, the GREG and the RHC estimators under various sampling designs}
\label{table 9}
\centering
\begin{tabular}{ccc} 
\hline
\multirow{2}{*}{Relative efficiency}& Jun, July and August & December\\
& in $2010$& in $2010$\\
\hline
RE($\hat{\overline{Y}}_{GREG}$, SRSWOR $|$ $\hat{\overline{Y}}_{HT}$, SRSWOR)& $1.529$& $1.805$\\
\hline
RE($\hat{\overline{Y}}_{GREG}$, RS $|$ $\hat{\overline{Y}}_{HT}$, RS)& $1.427$ & $1.263$\\
\hline
RE($\hat{\overline{Y}}_{GREG}$, RHC $|$ $\hat{\overline{Y}}_{RHC}$, RHC)& $1.531$ & $1.251$\\
\hline
\end{tabular}
\end{table}
\begin{table}[h]
\centering
\caption{Relative efficiencies of the GREG estimator under various sampling designs}
\label{table 8}
\centering
\begin{tabular}{ccc} 
\hline
\multirow{2}{*}{Relative efficiency}& Jun, July and August & December\\
& in $2010$& in $2010$\\
\hline
RE($\hat{\overline{Y}}_{GREG}$, RS $|$ $\hat{\overline{Y}}_{GREG}$, SRSWOR)&$2.32$& $1.76$\\
\hline
RE($\hat{\overline{Y}}_{GREG}$, RS $|$ $\hat{\overline{Y}}_{GREG}$, RHC)& $1.018$& $1.012$\\
\hline
\end{tabular}
\end{table}

\section{Determining the degree of heteroscedasticity $\eta$}\label{sec 5}
 In this section, we provide two methods for checking whether the degree of heteroscedasticity $\eta$ in the linear regression model in \eqref{eq 2} in  Section  \ref{sec 4} is bigger than $0.5$ or smaller than $0.5$ based on a pilot survey using SRSWOR. In the first method, we estimate $\eta$ based on some non-parametric estimation methods. In the second method, we choose $\eta$ based on statistical tests of heteroscedasticity.
\subsection{Estimation of $\eta$}\label{subsec 5.1}
Under the linear regression model in \eqref{eq 2}, we have the conditional total variance $tr(cov_{\mathbb{P}}(Y_i|Z_i,X_i))$ =$tr(cov_{\mathbb{P}}(\epsilon_i)) X_i^{2\eta}$, where $tr$ denotes the trace of an operator, and $cov_{\mathbb{P}}(Y_i|Z_i,X_i)$ is the conditional covariance operator of $Y_i$ given $(Z_i,X_i)$. Thus according to   the linear model \eqref{eq 2},  $\log\big(tr(cov_{\mathbb{P}}(Y_i|Z_i,X_i))\big)$ and $log(X_i)$ are linearly related with the slope $2\eta$. Now, in the case of $\mathcal{H}$=$L^2[0,T]$, we have $tr(cov_{\mathbb{P}}(Y_i|Z_i,X_i))$=$\int_{[0,T]}var_{\mathbb{P}}(Y_i(t)|Z_i,X_i)dt$, where $var_{\mathbb{P}}(Y_i(t)|Z_i,X_i)$ is the conditional variance of $Y_i(t)$ given $(Z_i,X_i)$. Suppose that the observations $\{(Y_i,Z_i,X_i):1\leq i\leq N\}$ in the population are generated from the linear model in \eqref{eq 2} and the observations on the study variable $y$ are obtained at $t_1,\ldots,t_r$ in $[0,T]$. Further, suppose that $s$ is a sample of size $n$ drawn based on a pilot survey using SRSWOR. Then, we estimate $tr(cov_{\mathbb{P}}(Y_i|Z_i,X_i))$ based on $\big\{\big(Y_i(t_l),Z_i,X_i\big): i\in s, l=1,\ldots,r \big\}$ as follows. For any $i\in s$ and $l$=$1,\ldots, r$, we first construct  the local average estimator of $E_{\mathbb{P}}(Y_i(t_l)|Z_i,X_i)$, the conditional mean of $Y_i(t_l)$ given $(Z_i,X_i)$, as 
\begin{align}
\begin{split}
&\hat{E}_{\mathbb{P}}(Y_i(t_l)|Z_i,X_i)=\sum_{k\in s} \prod_{j=1}^d \mathds{1}_{[|Z_{ji}-Z_{jk}|\leq h_{1l}]} \mathds{1}_{[|X_i-X_k|\leq h_{1l}]}Y_k(t_l)\times\\
&\bigg(\sum_{k\in s}\prod_{j=1}^d \mathds{1}_{[|Z_{ji}-Z_{jk}|\leq h_{1l}]} \mathds{1}_{[|X_i-X_k|\leq h_{1l}]}\bigg)^{-1}.
\end{split}
\end{align}
 Here, $Z_{ji}$ is $j^{th}$ component of $Z_i$. For any given $l$=$1,\ldots, r$, we compute the bandwidth $h_{1l}$  using  leave one out cross validation based on $\{(Y_i(t_l),Z_i,X_i): i\in s\}$. Now, using $\{\hat{E}_{\mathbb{P}}(Y_i(t_l)|Z_i,X_i): i\in s\}$, we estimate $var_{\mathbb{P}}(Y_i(t_l)|Z_i, X_i)$ by  local sample variance
\begin{align}
\begin{split}
&\widehat{var}_{\mathbb{P}}(Y_i(t_l)|Z_i, X_i)=\sum_{k\in s} \prod_{j=1}^d \mathds{1}_{[|Z_{ji}-Z_{jk}|\leq h_{2l}]} \mathds{1}_{[|X_i-X_k|\leq h_{2l}]}\times \\
&\big(Y_k(t_l)-\hat{E}_{\mathbb{P}}(Y_k(t_l)|Z_k,X_k)\big)^2\bigg(\sum_{k\in s}\prod_{j=1}^d \mathds{1}_{[|Z_{ji}-Z_{jk}|\leq h_{2l}]} \mathds{1}_{[|X_i-X_k|\leq h_{2l}]}\bigg)^{-1}
\end{split}
\end{align}
 for any $i\in s$ and $l$=$1,\ldots, r$. We compute the bandwidth $h_{2l}$ based on $\big\{\big(\big(Y_i(t_l)-\hat{E}_{\mathbb{P}}(Y_i(t_l)|Z_i,X_i)\big)^2,Z_i,$ $X_i\big): i\in s \big\}$  using leave one out cross validation  in the same way as we compute the bandwidth $h_{1l}$. Now, given $\big\{\widehat{var}_{\mathbb{P}}(Y_i(t_l)|Z_i, X_i): i\in s, l=1,\ldots,r\big\}$, we estimate $tr(cov_{\mathbb{P}}(Y_i|Z_i,X_i))$ by $T r^{-1}\times$ $\sum_{l=1}^r \widehat{var}_{\mathbb{P}}(Y_i(t_l)|Z_i,X_i)$ for any $i\in s$. Then, we fit a least square regression line to the data $\big\{\big(log\big(T r^{-1}$ $\times\sum_{l=1}^r \widehat{var}_{\mathbb{P}}(Y_i(t_l)|Z_i,X_i)\big), $ $log(X_i)\big): i\in s\big\}$, and compute the slope of this line. The slope, say $\hat{\theta}$, is expected to be close to $2\eta$  if the linear model in \eqref{eq 2} holds.  Thus $\hat{\eta}$=$0.5\hat{\theta}$ can be considered as an estimator of $\eta$. We demonstrate this method based on real and synthetic data as follows. 
\begin{itemize}
\item[(i)]Let us first consider the datasets from  Section  \ref{subsec 6.3}. Recall from  Section  \ref{subsec 6.3} that in the case of the estimation of the mean electricity consumption curve in June, July and August of $2010$, we have $r$=$4416$. On the other hand, in the case of the estimation of the mean electricity consumption curve in the December of $2010$, we have $r$=$1488$. Also, recall that $T$=$30 r$  in the cases of the estimation problems for both the datasets.  We draw $I$=$100$ samples each of size $n$=$500$ from these populations using SRSWOR and estimate $\eta$ as above based on these samples. Then, we compute the proportion of cases, when $\hat{\eta}>0.5$. It follows that this proportion is $0.72$ in the case of the estimation of the mean electricity consumption curve in June, July and August of $2010$ and $0.76$ in the case of the estimation of the mean electricity consumption curve in the December of $2010$. Recall from  Section  \ref{subsec 6.3} that  in the cases of the estimation problems for both the datasets,  the GREG estimator under RS sampling design is more efficient than any other estimator under any other sampling design when compared in terms of relative efficiencies. These corroborate the results stated in Theorem \ref{thm 3}.
\vspace{.2cm}
  
\item[(ii)] Next, suppose that finite populations each of size $N$=$5000$ are generated from linear models in the same way as in  Section  \ref{subsec 6.1}. Recall from  Section  \ref{subsec 6.1} that $r$=$100$, $T$=$1$ and $\eta$=$0.1 k$ for $k$=$0,\ldots,10$ in this case. We draw $I$=$100$ samples each of size $n$=$500$ from these populations using SRSWOR. Based on each sample $s$, we estimate $\eta$. Now, suppose that $ \hat{\eta}_{lk}$ is the estimate of $0.1 k$ based on the $l^{th}$ sample for $k$=$0,\ldots,10$ and $l$=$1,\ldots,I$. Then, we compute the proportion $l^{-1}\#\{l: \hat{\eta}_{lk}\leq 0.5\}$ for different $\eta$'s and $\beta(t)$'s (see  Section  \ref{subsec 6.1}) in Table \ref{table 3}. It follows from Table \ref{table 3} that for the values of $\eta$ smaller than $0.5$, the proportions are close to $1$. On the other hand, these proportions gradually decrease and become $0$, when $\eta$ becomes larger than $0.5$.  Once again, these corroborate the results stated in Theorem \ref{thm 3}.  
\end{itemize}

\begin{table}[h]
\centering
\caption{Proportion of cases when $\hat{\eta}\leq 0.5$ for different $\eta$'s and $\beta(t)$'s in the case of synthetic data}
\label{table 3}
\begin{center}
\begin{tabular}{cccc} 
\hline
$\eta$& $\beta(t)$=$1$& $\beta(t)$=$t$& $\beta(t)$=$1-(t-0.5)^2$\\
\hline
$0$& $1$ &$1$ & $1$\\
\hline
$0.1$& $1$ & $1$ & $1$\\
\hline
$0.2$& $1$ & $1$ & $1$\\
\hline
$0.3$& $1$ & $0.99$ & $0.98$\\
\hline
$0.4$& $0.99$ & $0.95$ & $0.96$\\
\hline
$0.5$& $0.9$ & $0.92$ & $0.94$\\
\hline
$0.6$& $0.52$ & $0.56$  & $0.59$\\
\hline
$0.7$& $0.2$ & $0.24$ & $0.2$\\
\hline
$0.8$& $0.01$  & $0.02$& $0$\\
\hline
$0.9$& $0$ & $0$ & $0$\\
\hline
$1$&  $0$& $0$& $0$\\
\hline
\end{tabular}
\end{center}
\end{table}
\subsection{Tests for $\eta$}
Under the linear regression model in \eqref{eq 2}, in the case of $\mathcal{H}$=$L^2[0,T]$, we have $X_i^{-\eta}\int_{[0,T]} Y_i(t)dt$= $X_i^{-\eta} \int_{[0,T]} \beta_0(t)dt+\sum_{j=1}^d\big( \int_{[0,T]} \beta_j(t)dt\big) Z_{ji} X_i^{-\eta}+\int_{[0,T]}\epsilon_i(t)dt$ for $i$=$1,\ldots,N$. As in the preceding section, suppose that observations on the study variable $y$ are obtained at $t_1,\ldots,t_r$ in $[0,T]$, and $s$ is a sample of size $n$ drawn based on a pilot survey using SRSWOR. Then we can say that $\{ (\tilde{Y}_i X_i^{-\eta},(1,Z_i)X_i^{-\eta}) : i\in s\}$ are generated from a homoscedastic linear model. Here, $\tilde{Y}_i$=$\int_{[0,T]} Y_i(t)dt$ for $i\in s$. We approximate $\tilde{Y}_i$ by $\hat{Y}_i$=$T r^{-1} \sum_{l=1}^r Y_i(t_l)$. Next, for every $\eta$ in $\{0.1 k: k=0,\ldots,10\}$, we test the null hypothesis $H_{0,\eta}:$ the data $\{(\hat{Y}_i X_i^{-\eta},(1,Z_i) X_i^{-\eta}) :i\in s\}$ are generated from a homoscedastic linear model against the alternative hypothesis $H_{1,\eta}:$ heteroscedasticity is present in the data $\{(\hat{Y}_i X_i^{-\eta},(1,Z_i) X_i^{-\eta}) :i\in s\}$. For this purpose, we use the Breusch-Pagan (BP, see \cite{breusch1979simple}), the White (see \cite{white1980heteroskedasticity}) and the Glejser (see \cite{glejser1969new}) tests  because these are some well known tests for heteroscedasticity.  In these tests, the residuals obtained from the ordinary least square regression between the response and the explanatory variables are expressed in terms of explanatory variables by means of different parametric models, and it is checked whether the explanatory variables have any influence on these residuals. Large $P$-values of the BP, the White and the Glejser tests  are indicative of  substantial evidence in favour of $H_{0,\eta}$. Thus, we select the $\eta$ from $\{0.1 k: k=0,\ldots,10\}$ for which we have the highest $P$-value. We denote this $\eta$ by $\hat{\eta}$. Now, we demonstrate this method based on real and synthetic data as follows.

\begin{itemize}
\item[(i)] As in the preceding section, let us first consider the datasets used in  Section  \ref{subsec 6.3}. We draw $I$=$100$ samples each of size $n$=$500$ from these datasets using SRSWOR and compute $\hat{\eta}$ as above based on each of these samples. Then, for each of the three tests and each of the datasets, we compute the proportion of cases, when $\hat{\eta}>0.5$ (see Table \ref{table 5}). As mentioned in the preceding section, in the cases of both the estimation problems, the GREG estimator under RS sampling design becomes the most efficient estimator when compared in terms of relative efficiencies. These corroborate the results stated in Theorem \ref{thm 3}.
\item[(ii)] Next, we determine $\eta$ as above based on the synthetic data considered in  Section  \ref{subsec 5.1}. We draw $I$=$100$ samples each of size $n$=$500$ from these datasets using SRSWOR and compute $\hat{\eta}$ based on each of these samples. Then, for each of the three tests, every $\eta$ in $\{0.1k:k=1,\ldots,10\}$ and each $\beta(t)$ (see  Section  \ref{subsec 6.1}), we compute the proportion of cases, $\hat{\eta}\leq 0.5$ (see Table \ref{table 4}). As in the previous section, it follows from Table \ref{table 4} that for the values of $\eta$ smaller than $0.5$, these proportions are close to $1$. On the other hand, these proportions gradually decrease and become $0$, when $\eta$ becomes larger than $0.5$.  Once again, these corroborate the results stated in Theorem \ref{thm 3}.  
\end{itemize}
\begin{table}[h]
\caption{Proportion of cases when $\hat{\eta}> 0.5$ for different tests and datasets in the case of electricity consumption data}
\label{table 5}
\centering
\begin{tabular}{ccc} 
\hline
\multirow{2}{*}{Test}& Jun, July and August & December\\
& in $2010$& in $2010$\\
\hline
BP & $0.79$& $0.83$\\
\hline
White & $0.76$& $0.78$\\
\hline
Glejser & $0.84$& $0.8$\\
\hline
\end{tabular}
\end{table}
\begin{table}[h] 
\caption{Proportion of cases when $\hat{\eta}\leq 0.5$ for different $\eta$'s and $\beta(t)$'s in the case of synthetic data}
\label{table 4}
\begin{center}
\renewcommand{\arraystretch}{0.8}
\begin{tabular}{cccccccccc} 
\hline
$\eta$& \multicolumn{3}{c}{$\beta(t)$=$1$}& \multicolumn{3}{c}{$\beta(t)$=$t$}& \multicolumn{3}{c}{$\beta(t)$=$1-(t-0.5)^2$}\\
& BP &White &Glejser & BP & White& Glejser & BP& White& Glejser\\
\hline
$0$& $1$ & $1$ &$1$ &$1$ &$1$& $1$& $1$& $1$& $1$\\
\hline
$0.1$& $1$ &$1$& $1$ & $1$& $1$ & $1$& $1$& $1$& $1$\\
\hline
$0.2$& $0.96$& $0.99$ & $0.99$ & $0.99$& $1$ & $1$& $1$& $0.97$& $0.93$\\
\hline
$0.3$& $0.95$& $0.77$& $0.98$  & $0.9$& $0.91$& $0.93$& $0.98$ & $0.9$& $0.88$\\
\hline
$0.4$& $0.85$& $0.75$& $0.84$ & $0.83$& $0.72$ & $0.88$& $0.87$& $0.9$& $0.79$\\
\hline
$0.5$& $0.6$& $0.65$& $0.68$ & $0.67$& $0.58$ & $0.69$& $0.58$& $0.69$& $0.72$\\
\hline
$0.6$& $0.47$& $0.29$& $0.36$ & $0.43$& $0.46$  & $0.47$& $0.39$& $0.45$& $0.36$\\
\hline
$0.7$& $0.17$& $0.22$& $0.15$ & $0.16$& $0.21$ & $0.14$& $0.13$& $0.25$& $0.17$\\
\hline
$0.8$& $0.09$& $0.07$& $0.06$ & $0.07$& $0.06$& $0.05$& $0.03$& $0.04$& $0.09$\\
\hline
$0.9$& $0.01$& $0.02$ & $0$ & $0.01$& $0.01$ & $0.01$& $0.01$& $0.01$& $0$\\
\hline
$1$&  $0$& $0$& $0$ & $0$& $0$& $0$& $0$& $0$& $0$\\
\hline
\end{tabular}
\end{center}
\end{table}

%%%%%%%%%%%%%%%%%%%%%%%%%%%%%%%%%%%%%%%%%%%%%%
%% Example with single Appendix:            %%
%%%%%%%%%%%%%%%%%%%%%%%%%%%%%%%%%%%%%%%%%%%%%%
\begin{appendix}

\section*{Appendix A: Sampling designs without replacement  considered in the paper}\label{appn 1}
\textbf{Simple random sampling without replacement (SRSWOR):} In SRSWOR, $n$ units are selected from the population $\mathcal{P}$ such that any subset of $n$ units has the same probability =$\big( ^N C_n\big)^{-1}$ of being selected.
\vspace{.2cm}

\noindent{\textbf{Lahiri-Midzuno-Sen (LMS) sampling design (\cite{Lahiri1951sampling}, \cite{Midzuno1952sampling} and \cite{Sen1953sampling}):} In LMS sampling design, the first unit is selected from $\mathcal{P}$, where the $i^{th}$ population unit $U_i$ has the probability =$X_i\big(\sum_{i=1}^N X_i\big)^{-1}$ of being selected for $i$=$1,\ldots,N$. Following the first draw, $n-1$ units are selected from the remaining $N-1$ units in $\mathcal{P}$ using SRSWOR. One can show that in this sampling design, the selection probability of a sample is proportional to the total of the values of the size variable $x$ for the sampled units.}
\vspace{.2cm}

\noindent{\textbf{Rao-Hartley-Cochran (RHC) sampling design (\cite{rao1962simple}):} In RHC sampling design, $\mathcal{P}$ is first divided randomly into $n$ disjoint groups, say $\mathcal{P}_1,\ldots, \mathcal{P}_n$ of sizes $N_1,\ldots,N_n$, respectively, by taking a sample of $N_1$ units from $N$ units using SRSWOR, then a sample of $N_2$ units from the remaining $N-N_1$ units using SRSWOR, then a sample of $N_3$ units from the remaining $N-N_1-N_2$ units using SRSWOR and so on. Following this random split, one unit is selected from each group independently. For each $r$=$1,\ldots,n$, the $q^{th}$ unit from $\mathcal{P}_r$ is selected with probability =$X^{\prime}_{qr}\big(\sum_{q=1}^{N_r} X^{\prime}_{qr}\big)^{-1}$, where $X^{\prime}_{qr}$ is the $x$ value of the $q^{th}$ unit in $\mathcal{P}_r$.}
\vspace{.2cm}

\noindent{\textbf{Rejective sampling design (\cite{MR0178555}):} Suppose that $\alpha_1,\ldots,\alpha_N$ are such that $\alpha_i>0$ for any $i$=$1,\ldots,N$ and $\sum_{i=1}^N \alpha_i$=$1$. Then, in the rejective sampling design, $n$ units are first drawn with replacement, where $U_i$ is selected with probability =$\alpha_i$, $i$=$1,\ldots,N$. If any population unit is selected in the sample more than once, the sample is rejected and the entire procedure is repeated until $n$ distinct units are selected in the sample.}
\vspace{.2cm}

\noindent{\textbf{High entropy sampling design (\cite{MR1624693}):} A sampling design $P(s,\omega)$ is called high entropy sampling design if $D(P||R)$=$\sum_{s \in \mathcal{S}}P(s,\omega)\log\big(P(s,\omega)(R(s,\omega))^{-1}\big)\rightarrow 0$ as $\nu\rightarrow\infty$ \textit{a.s.} $[\mathbb{P}]$ for some rejective sampling design $R(s,\omega)$.}
\vspace{.2cm}

\noindent{\textbf{$\pi$PS sampling design (\cite{MR1624693} and \cite{bondesson2006pareto}):} A sampling design is called $\pi$PS sampling design if its inclusion probabilities $\{\pi_i\}_{i=1}^N$ satisfy the condition $\pi_i$=$nX_i\big(\sum_{i=1}^N X_i\big)^{-1}$ for $i$=$1,\ldots,N$.}
\vspace{.2cm}

\noindent{\textbf{High entropy $\pi$PS (HE$\pi$PS) sampling design:} A sampling design is called a HE$\pi$PS sampling design if it is a high entropy sampling design as well as a $\pi$PS sampling design. }
\vspace{.2cm}

\noindent{\textbf{Rao-Sampford (RS) sampling design (\cite{MR1624693}):} In RS sampling design, a population unit is first selected in such a way that $U_i$ has the probability =$X_i\big(\sum_{i=1}^N X_i\big)^{-1}$ of being selected for $i$=$1,\ldots,N$. After replacing this unit back into the population, $n-1$ units are drawn with replacement, where $U_i$ is selected with probability =$\lambda_i(1-\lambda_i)^{-1}\big(\sum_{i=1}^N \lambda_i(1-\lambda_i)^{-1}\big)^{-1}$ for $\lambda_i$=$n X_i \big(\sum_{i=1}^N X_i\big)^{-1}$. If any population unit is selected in the sample more than once, the sample is rejected and the entire procedure is repeated until $n$ distinct units are selected in the sample. It was shown by \cite{MR1624693} that RS sampling design is a HE$\pi$PS sampling design.}
 
\section*{Appendix B: Lemmas and proofs of propositions and theorems}\label{appn 2}
 In this section, we give the proofs of different Propositions and Theorems.  For technical details, which are related to operator theory and used in the proofs of Lemmas, Propositions and Theorems, the reader is referred to \cite{hsing2015theoretical}. We shall first state some technical lemmas, which will be required to prove our main results.  Proofs of these lemmas are given in the supplement (\cite{deychaudhuri2023}). 
\begin{lemma}\label{lem 5} 
Suppose that Assumption \ref{ass C2} holds. Then, LMS sampling design is a high entropy sampling design. Moreover, under each of SRSWOR, LMS sampling design and any HE$\pi$PS sampling design, we have, for all sufficiently large $\nu$,
\begin{equation}\label{eq 1}
 M\leq  n^{-1} N \pi_{i} \leq M^{\prime}\text{ for some constants }M,M^{\prime}>0\text{ and all }1\leq i\leq N\text{ \textit{a.s.} }[\mathbb{P}].
\end{equation}
\end{lemma}
The condition \eqref{eq 1} was stated earlier in \cite{wang2011asymptotic}, \cite{MR3670194}, etc. However, it was not shown that \eqref{eq 1} holds under LMS and HE$\pi$PS sampling designs. Next, recall from the paragraph containing \eqref{eq 34} in  Section  \ref{sec 4} that $\gamma$=$\sum_{i=1}^n N_i(N_i-1)(N(N-1))^{-1}$ with $N_i$ being the size of the $i^{th}$ group formed randomly in the first step of the RHC sampling design (see Appendix A) for $i$=$1,\ldots,n$. We now state the following lemma.
\begin{lemma}\label{lem 4}
Suppose that Assumption \ref{ass C1} holds. Then, $n\gamma\rightarrow c$ for some $c\geq 1-\lambda>0$ as $\nu\rightarrow\infty$, where $\lambda$ is as in Assumption \ref{ass C1}.
\end{lemma}
Before, we state the next lemma, let us introduce some notations. Let $\{e_j\}_{j=1}^{\infty}$ be an orthonormal basis of the separable Hilbert space $\mathcal{H}$. Suppose that $V_i$ is either $Y_i$ or $Y_i-\overline{Y}-S_{zy}((Z_i-\overline{Z})S_{zz}^{-1})$, where   $S_{zy}$=$N^{-1}\sum_{i=1}^N(Z_i-\overline{Z})\otimes(Y_i-\overline{Y})$ and $S_{zz}$=$N^{-1}\sum_{i=1}^N(Z_i-\overline{Z})^T(Z_i-\overline{Z})$.  Further, suppose that $\Gamma_1$=$nN^{-2}\sum_{i=1}^N(V_i-T\pi_i)\otimes(V_i-T\pi_i)(\pi_i^{-1}-1)$, where  $T$=$\sum_{i=1}^N V_i(1-\pi_i)\big(\sum_{i=1}^N \pi_i(1-\pi_i)\big)^{-1}$,  and $\pi_i$ is the inclusion probability of the $i^{th}$ population unit. Moreover, in the case of RHC sampling design, suppose that  $\Gamma_{2}$=$n\gamma\overline{X}N^{-1}\sum_{i=1}^N\big(V_i-X_i\overline{V} \overline{X}^{-1}\big)\otimes\big(V_i-X_i\overline{V}\overline{X}^{-1}\big)X_i^{-1}$, where $\overline{V}$=$N^{-1}\sum_{i=1}^N V_i$, $\overline{X}$=$N^{-1}$ $\times\sum_{i=1}^N X_i$,  and $\gamma$ is as in the paragraph preceding Lemma \ref{lem 4}.  Let us also recall $b$ from Assumption \ref{ass C2}.  We now state the following lemma.
\begin{lemma}\label{lem 6}
Suppose that Assumptions \ref{ass C1}, \ref{ass C2} and \ref{ass C3} hold. Then, under SRSWOR and LMS sampling design, $\Gamma_1\rightarrow \Sigma_1$ with respect to the HS norm as $\nu\rightarrow\infty$ \textit{a.s.} $[\mathbb{P}]$ for some n.n.d. HS operator $\Sigma_1$. Also, $\sum_{j=1}^{\infty}\langle\Sigma_1 e_j,e_j\rangle<\infty$, and $\sum_{j=1}^{\infty}\langle \Gamma_1e_j,e_j\rangle \rightarrow \sum_{j=1}^{\infty}\langle\Sigma_1 e_j,e_j\rangle$ under the above sampling designs as $\nu\rightarrow\infty$ \textit{a.s.} $[\mathbb{P}]$. Further, if Assumption \ref{ass C1} holds with  $0\leq \lambda<E_{\mathbb{P}}(X_i)b^{-1}$,  and Assumptions \ref{ass C2} and \ref{ass C3} hold, then, the above results hold under any HE$\pi$PS sampling design. Moreover, if Assumptions \ref{ass C1}, \ref{ass C2} and \ref{ass C3} hold, then in the case of RHC sampling design, $\Gamma_2\rightarrow \Sigma_2$ with respect to the HS norm as $\nu\rightarrow\infty$ \textit{a.s.} $[\mathbb{P}]$ for some n.n.d. HS operator $\Sigma_2$. Also, $\sum_{j=1}^{\infty}\langle\Sigma_2 e_j,e_j\rangle<\infty$, and $\sum_{j=1}^{\infty}\langle \Gamma_2 e_j,e_j\rangle \rightarrow \sum_{j=1}^{\infty}\langle\Sigma_2 e_j,e_j\rangle$ as $\nu\rightarrow\infty$ \textit{a.s.} $[\mathbb{P}]$. 
\end{lemma}
Recall $\{V_i\}_{i=1}^N$ and $\{e_j\}_{j=1}^{\infty}$ from the paragraph preceding Lemma \ref{lem 6} and define $W_i$=$(\langle V_i,e_1\rangle,$ $\ldots,$ $\langle V_i,e_r\rangle)$ for $i$=$1,\ldots,N$ and $r\geq 1$. Suppose that $\hat{\overline{V}}_1$=$\sum_{i\in s}(N\pi_i)^{-1}V_i$,  $\hat{\overline{W}}_1$=$\sum_{i\in s} (N\pi_i)^{-1} W_i$ and $\overline{W}$=$N^{-1}\sum_{i=1}^N W_i$.  Moreover, suppose that  $\hat{\overline{V}}_2$=$\sum_{i\in s}(NX_i)^{-1}Q_iV_i$ and $\hat{\overline{W}}_2$=$\sum_{i\in s}(NX_i)^{-1}$ $\times  Q_i W_i$,  where $Q_i$ is the total of the $x$ values of that randomly formed group from which the $i^{th}$ population unit is selected in the sample by RHC sampling design (see Appendix A). Let us also assume that $S_k$=$\sqrt{n}(\hat{\overline{V}}_k-\overline{V})$, and $\Gamma_{k,r}$ is a $r\times r$ matrix such that $((\Gamma_{k,r}))_{jl}$=$\langle \Gamma_k e_j, e_l\rangle$ for $j,l=1,\ldots,r$, $k$=$1,2$ and $r\geq 1$. We now state the following lemma. 
 \begin{lemma}\label{lem 1}
Fix $r\geq 1$. Suppose that Assumptions \ref{ass C1}, \ref{ass C2} and \ref{ass C3} hold. Then, under SRSWOR and LMS sampling design, $(\langle S_1, e_1\rangle ,\ldots,\langle S_1, e_r\rangle )\xrightarrow{\mathcal{L}}N_r(0,\Sigma_{1,r})$ as $\nu\rightarrow\infty$ \textit{a.s.} $[\mathbb{P}]$, where $\Sigma_{1,r}$ is a $r\times r$ matrix such that $((\Sigma_{1,r}))_{jl}$=$\langle \Sigma_1 e_j, e_l\rangle$ for $j,l=1,\ldots,r$, and $\Sigma_1$ is as in the statement of Lemma \ref{lem 6}. Further, if Assumption \ref{ass C1} holds with  $0\leq \lambda<E_{\mathbb{P}}(X_i)b^{-1}$,  and Assumptions \ref{ass C2} and \ref{ass C3} hold, then, the above result holds under any HE$\pi$PS sampling design. Moreover, if Assumptions \ref{ass C1}, \ref{ass C2} and \ref{ass C3} hold, then $(\langle S_2, e_1\rangle ,\ldots,\langle S_2, e_r\rangle )\xrightarrow{\mathcal{L}}N_r(0,\Sigma_{2,r})$ as $\nu\rightarrow\infty$ under RHC sampling design \textit{a.s.} $[\mathbb{P}]$. Here, $\Sigma_{2,r}$ is a $r\times r$ matrix such that $((\Sigma_{2,r}))_{jl}$=$\langle \Sigma_2 e_j, e_l\rangle$ for $j,l=1,\ldots,r$, and $\Sigma_2$ is as in the statement of Lemma \ref{lem 6}. 
\end{lemma} 
Suppose that $\Pi_r$ denotes the orthogonal projection onto the linear span of $\{e_1,\ldots, e_r\}$, i.e., $\Pi_r(a)$=$\sum_{j=1}^r \langle a, e_j\rangle e_j$ for any $r\geq 1$ and $a\in\mathcal{H}$. Further, suppose that  $B_{1,r}$=$\{s\in\mathcal{S}: ||S_1-\Pi_r(S_1)||_{\mathcal{H}} >\epsilon\}$ and $B_{2,r}$=$\{s\in\mathcal{S}: ||S_2-\Pi_r(S_2)||_{\mathcal{H}} >\epsilon\}$ for any $\epsilon>0$.  Now, we state the following lemma.
\begin{lemma}\label{lem 2}
Suppose that Assumptions \ref{ass C1}, \ref{ass C2} and \ref{ass C3} hold, and $P(s,\omega)$ denotes one of SRSWOR and LMS sampling design. Then, for any $\epsilon>0$, $\lim_{r\rightarrow\infty}\overline{\lim}_{\nu\rightarrow\infty}\sum_{s\in B_{1,r}}P(s,\omega)$=$0$ \textit{a.s.} $[\mathbb{P}]$. Further, if Assumption \ref{ass C1} holds with  $0\leq \lambda<E_{\mathbb{P}}(X_i)b^{-1}$,  and Assumptions \ref{ass C2} and \ref{ass C3} hold, then the above result holds under any HE$\pi$PS sampling design. Moreover, suppose that Assumptions \ref{ass C1}, \ref{ass C2} and \ref{ass C3} hold, and $P(s,\omega)$ denotes RHC sampling design. Then, for any $\epsilon>0$, $\lim_{r\rightarrow\infty}\overline{\lim}_{\nu\rightarrow\infty}\sum_{s\in B_{2,r}}P(s,\omega)$=$0$ \textit{a.s.} $[\mathbb{P}]$.
\end{lemma} 
\begin{proof}[Proof of Proposition \ref{prop 1}]
Recall the expression of $\hat{\overline{Y}}_{HT}$ from \eqref{eq 4} in  Section  \ref{sec 2} and note that $S_1$= $\sqrt{n}(\hat{\overline{Y}}_{HT}-\overline{Y})$ if we substitute $V_i$=$Y_i$ in $S_1$. It follows from Lemma \ref{lem 1} that $(\langle S_1, e_1\rangle,\ldots,\langle S_1, e_r\rangle)\xrightarrow{\mathcal{L}}N_r(0,\Sigma_{1,r})$ as $\nu\rightarrow\infty$ for any $r\geq 1$ under SRSWOR, LMS sampling design and any HE$\pi$PS sampling design \textit{a.s.} $[\mathbb{P}]$. Here, $\Sigma_{1,r}$ is a $r\times r$ matrix such that $((\Sigma_{1,r}))_{jl}$=$\langle\Sigma_1 e_j,e_l\rangle$, and $\Sigma_1$=$\lim_{\nu\rightarrow\infty}\Gamma_1$ \textit{a.s.} $[\mathbb{P}]$. Further, it follows from the $1^{st}$ paragraph in the proof of Lemma \ref{lem 6}  in the supplement (\cite{deychaudhuri2023})  that $\Sigma_1$=$\Delta_1$ for SRSWOR and LMS sampling design, and $\Sigma_1$=$\Delta_2$ for any HE$\pi$PS sampling design. Therefore, by continuous mapping theorem, $\Pi_r(S_1)$=$\sum_{j=1}^r\langle S_1,e_j\rangle e_j\xrightarrow{\mathcal{L}}\mathcal{N}_1\circ \Pi_r^{-1}$ as $\nu\rightarrow\infty$ under the above sampling designs for any $r\geq 1$ \textit{a.s.} $[\mathbb{P}]$, where $\mathcal{N}_1$ is the Gaussian distribution in $\mathcal{H}$ with mean $0$ and covariance operator $\Sigma_1$. Moreover, in view of Lemma \ref{lem 2}, we have $\lim_{r\rightarrow\infty}\overline{\lim}_{\nu\rightarrow\infty} \sum_{s\in B_{1,r}} P(s,\omega)$=$0$ \textit{a.s.} $[\mathbb{P}]$, where $P(s,\omega)$ denotes one of the above sampling designs. Then, by Proposition $2.1$ in \cite{kundu2000central}, $\sqrt{n}(\hat{\overline{Y}}_{HT}-\overline{Y})\xrightarrow{\mathcal{L}}\mathcal{N}_1$ as $\nu\rightarrow\infty$ under the above sampling designs \textit{a.s.} $[\mathbb{P}]$. 
\end{proof}
\begin{proof}[Proof of Proposition \ref{prop 3}]
Recall the expression of $\hat{\overline{Y}}_{RHC}$ from \eqref{eq 5} in  Section  \ref{sec 2} and note that $S_2$= $\sqrt{n}(\hat{\overline{Y}}_{RHC}-\overline{Y})$ if we substitute $V_i$=$Y_i$ in $S_2$. It follows in view of Lemma \ref{lem 1} that under RHC sampling design, $(\langle S_2, e_1\rangle,\ldots,\langle S_2, e_r\rangle)\xrightarrow{\mathcal{L}}N_r(0,\Sigma_{2,r})$ as $\nu\rightarrow\infty$ for any $r\geq 1$ \textit{a.s.} $[\mathbb{P}]$. Here, $\Sigma_{2,r}$ is a $r\times r$ matrix such that $((\Sigma_{2,r}))_{jl}$=$\langle\Sigma_2 e_j,e_l\rangle$, and $\Sigma_2$=$\lim_{\nu\rightarrow\infty}\Gamma_2$ \textit{a.s.} $[\mathbb{P}]$. Further, it follows from the $2^{nd}$ paragraph in the proof of Lemma \ref{lem 6}  in the supplement (\cite{deychaudhuri2023})  that $\Sigma_2$=$ \Delta_3$.  Therefore, by continuous mapping theorem, $\Pi_r(S_2)$=$\sum_{j=1}^r\langle S_2,e_j\rangle e_j\xrightarrow{\mathcal{L}}\mathcal{N}_2\circ \Pi_r^{-1}$ as $\nu\rightarrow\infty$ under RHC sampling design for any $r\geq 1$ \textit{a.s.} $[\mathbb{P}]$, where $\mathcal{N}_2$ is the Gaussian distribution in $\mathcal{H}$ with mean $0$ and covariance operator $\Sigma_2$. Next, it follows from Lemma \ref{lem 2} that $\lim_{r\rightarrow\infty}\overline{\lim}_{\nu\rightarrow\infty} \sum_{s\in B_{2,r}} P(s,\omega)$=$0$ \textit{a.s.} $[\mathbb{P}]$, where $P(s,\omega)$ denotes RHC sampling design. Then, by Proposition $2.1$ in \cite{kundu2000central}, $\sqrt{n}(\hat{\overline{Y}}_{RHC}-\overline{Y})\xrightarrow{\mathcal{L}}\mathcal{N}_2$ as $\nu\rightarrow\infty$ under RHC sampling design \textit{a.s.} $[\mathbb{P}]$.
\end{proof}
Before we write the proof of Proposition \ref{prop 2}, we state the following lemma.  The proof of this lemma is given in the supplement (\cite{deychaudhuri2023}).  Let $V_i^{\sharp}$ be one of $Y_i\otimes Z_i, Z_i^T Z_i, Z_i$ and $1$ for $i$=$1,\ldots,N$. Also, let  $\overline{V}^{\sharp}$=$N^{-1}\sum_{i=1}^N V_i^{\sharp}$,  $S_1^{\sharp}$=$\sqrt{n}(\sum_{i\in s}(N\pi_i)^{-1} V_i^{\sharp}-\overline{V}^{\sharp})$ and $S_2^{\sharp}$=$\sqrt{n}(\sum_{i\in s}(NX_i)^{-1}Q_i V_i^{\sharp}-\overline{V}^{\sharp})$.  In this case, we denote the associated norm by $||\cdot||_{\mathcal{G}}$. Note that $||\cdot||_{\mathcal{G}}$=the Euclidean norm, when $V_i^{\sharp}$ is one of $Z_i^T Z_i, Z_i$ and $1$, and $||\cdot||_{\mathcal{G}}$=the HS norm, when $V_i^{\sharp}$=$Y_i\otimes Z_i$.  
\begin{lemma}\label{lem 3}
Suppose that Assumptions \ref{ass C1}, \ref{ass C2} and \ref{ass C3} hold. Then,  $||S_1^{\sharp}||_{\mathcal{G}}$=$O_p(1)$  under SRSWOR, LMS sampling design and any HE$\pi$PS sampling design, and  $||S_2^{\sharp}||_{\mathcal{G}}$=$O_p(1)$  under RHC sampling design  as $\nu\rightarrow\infty$ \textit{a.s.} $[\mathbb{P}]$.
\end{lemma}
\begin{proof}[Proof of Proposition \ref{prop 2}] Recall from \eqref{eq 6} in  Section  \ref{sec 2} that $\hat{\overline{Y}}_{GREG}$=$\hat{\overline{Y}}+\hat{S}_{zy}((\overline{Z}-\hat{\overline{Z}})\hat{S}_{zz}^{-1})$, where  $\hat{\overline{Y}}$=$\big(\sum_{i\in s}\pi_i^{-1}\big)^{-1}\sum_{i\in s}\pi_i^{-1}Y_i$, $\hat{\overline{Z}}$=$\big(\sum_{i\in s}\pi_i^{-1}\big)^{-1}\sum_{i\in s}\pi_i^{-1}Z_i$, $\hat{S}_{zz}$=$\big(\sum_{i\in s}\pi_i^{-1}\big)^{-1}\times$ $\sum_{i\in s}\pi_i^{-1}(Z_i-\hat{\overline{Z}})^T(Z_i-\hat{\overline{Z}})$, and $\hat{S}_{zy}$=$\big(\sum_{i\in s}\pi_i^{-1}\big)^{-1}\sum_{i\in s}\pi_i^{-1}(Z_i-\hat{\overline{Z}})\otimes(Y_i-\hat{\overline{Y}})$.  Note that 
\begin{equation}\label{eq 27}
\hat{\overline{Y}}_{GREG}-\overline{Y}=\Theta(\hat{\overline{V}}_1-\overline{V})+B, 
\end{equation}
 where $\hat{\overline{V}}_1$=$\sum_{i\in s}(N\pi_i)^{-1}V_i$, $V_i$=$Y_i-\overline{Y}-S_{zy}((Z_i-\overline{Z})S_{zz}^{-1})$, $\Theta$=$(\sum_{i\in s}\pi_i^{-1})^{-1}$, $B$=$S_{zy}((\hat{\overline{Z}}-\overline{Z}) S_{zz}^{-1})-\hat{S}_{zy}((\hat{\overline{Z}}-\overline{Z}) \hat{S}_{zz}^{-1})$, $S_{zy}$=$N^{-1}\sum_{i=1}^N(Z_i-\overline{Z})\otimes(Y_i-\overline{Y})$, and $S_{zz}$=$N^{-1}\sum_{i=1}^N(Z_i-\overline{Z})^T(Z_i-\overline{Z})$. Using Lemmas \ref{lem 1} and \ref{lem 2}, it can be shown in the same way as in the proof of Proposition \ref{prop 1} that as $\nu\rightarrow\infty$, $\sqrt{n}(\hat{\overline{V}}_1-\overline{V})\xrightarrow{\mathcal{L}}\mathcal{N}_3$ under SRSWOR, LMS sampling design and any HE$\pi$PS sampling design \textit{a.s.} $[\mathbb{P}]$, where $\mathcal{N}_3$ is the Gaussian distribution in $\mathcal{H}$ with mean $0$ and covariance operator $\Sigma_1$. Here, $\Sigma_1$=$\lim_{\nu\rightarrow\infty}\Gamma_1$ \textit{a.s.} $[\mathbb{P}]$. It follows from the last paragraph in the proof of Lemma \ref{lem 6}  in the supplement (\cite{deychaudhuri2023})  that  $\Sigma_1$=$\Delta_4$ under SRSWOR and LMS sampling design, and $\Sigma_1$=$\Delta_5$  under any HE$\pi$PS sampling design. Now, to establish the weak convergence of $\sqrt{n}(\hat{\overline{Y}}_{GREG}-\overline{Y})$ under the above sampling designs \textit{a.s.} $[\mathbb{P}]$, it is enough to show that $\Theta\xrightarrow{p} 1$ and $\sqrt{n}B \xrightarrow{p} 0$ under these sampling designs as $\nu\rightarrow\infty$ \textit{a.s.} $[\mathbb{P}]$.
\par
\vspace{.2cm}

Suppose that  $||\cdot||_{op}$  denotes the operator norm. Note that except the operator norm, we use only the HS norm for the operators considered in this article and denote it by  $||\cdot||_{HS}$.  Also, note that 
\begin{equation}\label{eq 28}
||B||_{\mathcal{H}}\leq (||S^{-1}_{zz}||_{op}||S_{zy}-\hat{S}_{zy}||_{op}+||\hat{S}_{zy}||_{op}||S_{zz}^{-1}-\hat{S}_{zz}^{-1}||_{op})||\hat{\overline{Z}}-\overline{Z}||. 
\end{equation}
 It follows in view of Lemma \ref{lem 3} that as $\nu\rightarrow\infty$, 
\begin{align}\label{eq 29}
\begin{split}
&\bigg|\bigg|\sum_{i\in s}(N\pi_i)^{-1}(Y_i\otimes Z_i)-N^{-1}\sum_{i=1}^N (Y_i\otimes Z_i)\bigg|\bigg|_{HS}=o_p(1)\text{, } \bigg|\bigg|\sum_{i\in s}(N\pi_i)^{-1}Z_i^T Z_i-\\
& N^{-1}\sum_{i=1}^N Z_i^T Z_i\bigg|\bigg|=o_p(1), \sqrt{n}\bigg|\bigg|\hat{\overline{Z}}_1-\overline{Z}\bigg|\bigg|=O_p(1),\text{ and }\sum_{i\in s}(N\pi_i)^{-1}-1=o_p(1) 
\end{split}
\end{align}
 under the sampling designs considered in the previous paragraph \textit{a.s.} $[\mathbb{P}]$, where $\hat{\overline{Z}}_1$= $\sum_{i\in s}(N\pi_i)^{-1}Z_i$. Consequently, in view of Assumption \ref{ass C3},  
\begin{align}\label{eq 30}
\begin{split}
&\sqrt{n}\bigg|\bigg|\hat{\overline{Z}}-\overline{Z}\bigg|\bigg|=O_p(1), \bigg|\bigg|\hat{S}_{zz}-S_{zz}\bigg|\bigg|_{op}\leq \bigg|\bigg|\hat{S}_{zz}-S_{zz}\bigg|\bigg|=o_p(1)\text{ and }\\
&\bigg|\bigg|\hat{S}_{zy}-S_{zy}\bigg|\bigg|_{op} \leq \bigg|\bigg|\hat{S}^{*}_{zy}-S^{*}_{zy}\bigg|\bigg|_{HS}=o_p(1) 
\end{split}
\end{align}
 as $\nu\rightarrow\infty$ under these sampling designs \textit{a.s.} $[\mathbb{P}]$. Here,  $\hat{S}^{*}_{zy}$=$\big(\sum_{i\in s} \pi_i^{-1}\big)^{-1}\sum_{i\in s} \pi_i^{-1}(Y_i-\hat{\overline{Y}})\otimes(Z_i-\hat{\overline{Z}})$ and $S^{*}_{zy}$=$N^{-1}\sum_{i=1}^N(Y_i-\overline{Y})\otimes(Z_i-\overline{Z})$  are adjoints of $\hat{S}_{zy}$ and $S_{zy}$, respectively. Now, recall $C_{zz}$ and $C_{zy}$ from the $3^{rd}$ paragraph in the proof of Lemma \ref{lem 6}  in the supplement (\cite{deychaudhuri2023}).  Note that $||S_{zz}-C_{zz}||$=$o(1)$ and $||S_{zy}-C_{zy}||_{HS}$=$o(1)$  as $\nu\rightarrow\infty$ \textit{a.s.} $[\mathbb{P}]$ in view of Assumption \ref{ass C3}. Also, note that $C_{zz}^{-1}$ exists by Assumption \ref{ass C3}. Consequently,  $||S^{-1}_{zz}||_{op}$=$O(1)$, $||\hat{S}^{-1}_{zz}-S^{-1}_{zz}||_{op}$=$o_p(1)$ and $||\hat{S}_{zy}||_{op}$=$O_p(1)$  as $\nu\rightarrow\infty$ \textit{a.s.} $[\mathbb{P}]$. Thus  $\sqrt{n}||B||_{\mathcal{H}}$=$o_p(1)$  and $\Theta-1$=$o_p(1)$ as $\nu\rightarrow\infty$ under the above-mentioned sampling designs \textit{a.s.} $[\mathbb{P}]$. Hence, the weak convergence of $\sqrt{n}(\hat{\overline{Y}}_{GREG}-\overline{Y})$ follows under these sampling designs by using Proposition $2.1$ in \cite{kundu2000central}.
\par
\vspace{.2cm}

Let us next consider the RHC sampling design. Recall from  Section  \ref{sec 2} that we consider $\hat{\overline{Y}}_{GREG}$ under RHC sampling design with $\pi_i^{-1}$ replacing  $Q_i X_i^{-1}$.  Then, under this sampling design, 
\begin{equation}\label{eq 31}
\hat{\overline{Y}}_{GREG}-\overline{Y}=\Theta(\hat{\overline{V}}_2-\overline{V})+B, 
\end{equation}
where $\hat{\overline{V}}_2$=$\sum_{i\in s}(N X_i)^{-1}Q_iV_i$  for $V_i$=$Y_i-\overline{Y}-S_{zy}((Z_i-\overline{Z})S_{zz}^{-1})$, and $\Theta$ and $B$ are the same as defined in the $1^{st}$ paragraph of this proof with $\pi_i^{-1}$ replaced by  $Q_i X_i^{-1}$.  Using Lemmas \ref{lem 1} and \ref{lem 2}, it can be shown in a similar way as in the proof of Proposition \ref{prop 3} that $\sqrt{n}(\hat{\overline{V}}_2-\overline{V})\xrightarrow{\mathcal{L}}\mathcal{N}_4$ as $\nu\rightarrow\infty$ under RHC sampling design \textit{a.s.} $[\mathbb{P}]$, where $\mathcal{N}_4$ is the Gaussian distribution in $\mathcal{H}$ with mean $0$ and covariance operator $\Sigma_2$. It follows from the last paragraph in the proof of Lemma \ref{lem 6}  in the supplement (\cite{deychaudhuri2023})  that  $\Sigma_2$=$\Delta_6$=$\lim_{\nu\rightarrow\infty}\Gamma_2$  \textit{a.s.} $[\mathbb{P}]$. Moreover, using Lemma \ref{lem 3}, it can be shown in the same way as in the preceding paragraph of this proof that $\Theta\xrightarrow{p} 1$ and $\sqrt{n}B \xrightarrow{p} 0$ as $\nu\rightarrow\infty$ under RHC sampling design \textit{a.s.} $[\mathbb{P}]$. Threfore, the weak convergence of $\sqrt{n}(\hat{\overline{Y}}_{GREG}-\overline{Y})$ follows under this sampling design by using Proposition $2.1$ in \cite{kundu2000central}. 
\end{proof} 
\begin{proof}[Proof of Theorem \ref{thm 2}]
Let us recall the expressions of $\Delta_1$ and $\Delta_4$ from the proof of Lemma \ref{lem 6}  in the supplement (\cite{deychaudhuri2023}).  It follows from the proof of Proposition \ref{prop 2} that \textit{a.s.} $[\mathbb{P}]$, $\sqrt{n}(\hat{\overline{Y}}_{GREG}-\overline{Y})$ has the same asymptotic covariance operator $\Delta_4$ under SRSWOR and LMS sampling design. It further follows from the proof of Proposition \ref{prop 1} that \textit{a.s.} $[\mathbb{P}]$, the asymptotic covariance operator of $\sqrt{n}(\hat{\overline{Y}}_{HT}-\overline{Y})$ is $\Delta_1$ under SRSWOR as well as LMS sampling design. Let $A_i$=$\langle Y_i, a\rangle$ for $a\in\mathcal{H}$ and $i$=$1,\ldots,N$. Then, we have 
\begin{align}\label{eq 32}
\begin{split}
&\langle(\Delta_1-\Delta_4)a,a\rangle=(1-\lambda)\big(E_{\mathbb{P}}(A_i-E_{\mathbb{P}}(A_i))^2-E_{\mathbb{P}}(A_i-E_{\mathbb{P}}(A_i)-C_{za}C_{zz}^{-1}(Z_i-E_{\mathbb{P}}(Z_i))^T)^2\big)\\
&=(1-\lambda)C_{za}C_{zz}^{-1}C_{za}^T 
\end{split}
\end{align}
 for $C_{za}$=$E_{\mathbb{P}}(A_i-E_{\mathbb{P}}(A_i))(Z_i-E_{\mathbb{P}}(Z_i))$ and $C_{zz}$=$E_{\mathbb{P}}(Z_i-E_{\mathbb{P}}(Z_i))^T(Z_i-E_{\mathbb{P}}(Z_i))$. Note that we have $C_{za}C_{zz}^{-1} C_{za}^T\geq 0$ for any $a\in\mathcal{H}$ by Assumption \ref{ass C3}. In fact, there exists $a\in \mathcal{H}$ such that $a\neq 0$ and $C_{za}$=$0$. Therefore, $\Delta_1-\Delta_4$ is p.s.d. Hence, \textit{a.s.} $[\mathbb{P}]$, the GREG estimator is asymptotically at least as efficient as the HT estimator under SRSWOR and LMS sampling design. Moreover, \textit{a.s.} $[\mathbb{P}]$, both the HT and the GREG estimators have the same asymptotic distribution under SRSWOR and LMS sampling design.  
\end{proof}
\begin{proof}[Proof of Theorem \ref{thm 1}]
Let us recall the expressions of  $\Delta_2$, $\Delta_3$, $\Delta_5$ and $\Delta_6$  from the proof of Lemma \ref{lem 6}  in the supplement (\cite{deychaudhuri2023}).  It can be shown from the proofs of Propositions \ref{prop 3} and \ref{prop 2} that \textit{a.s.} $[\mathbb{P}]$, asymptotic covariance operators of $\sqrt{n}(\hat{\overline{Y}}_{RHC}-\overline{Y})$ and $\sqrt{n}(\hat{\overline{Y}}_{GREG}-\overline{Y})$ under RHC sampling design are  $\Delta_3$ and $\Delta_6$, respectively.  Now, it follows from the linear regression model in \eqref{eq 2} in  Section  \ref{sec 4} that 
\begin{align}\label{eq 35}
\begin{split}
&\langle \Delta_3 a,a\rangle=c\bigg[\mu_x E_{\mathbb{P}}(\tilde{\epsilon}_i)^2 E_{\mathbb{P}}\big(X_i^{2\eta-1}\big)+\mu_x E_{\mathbb{P}}\bigg(\tilde{\beta}_0+\sum_{j=1}^d\tilde{\beta}_jZ_{ji}\bigg)^2 X_i^{-1}-\bigg(\sum_{j=0}^d\tilde{\beta}_j\mu_j\bigg)^2\bigg]\\
&\text{and }\langle \Delta_6 a,a\rangle=c \mu_x E_{\mathbb{P}}(\tilde{\epsilon}_i)^2 E_{\mathbb{P}}\big(X_i^{2\eta-1}\big), 
\end{split}
\end{align}
 where $c$=$\lim_{\nu\rightarrow\infty}n\gamma>0$, $a\in\mathcal{H}$, $\tilde{\epsilon}_i$=$\langle\epsilon_i,a\rangle$, $\mu_x$=$E_{\mathbb{P}}(X_i)$, $\tilde{\beta}_j$=$\langle\beta_j,a\rangle$ for $j$=$0,\ldots, d$, $\mu_0$=$1$, and $\mu_j$=$E_{\mathbb{P}}(Z_{ji})$ for $j$=$1,\ldots,d$. Therefore, 
\begin{equation}\label{eq 37}
\langle (\Delta_3-\Delta_6) a,a\rangle=c \mu_x E_{\mathbb{P}}\bigg(\tilde{\beta}_0+\sum_{j=1}^d\tilde{\beta}_jZ_{ji}-X_i\sum_{j=0}^d\tilde{\beta}_j\mu_j\mu_x^{-1}\bigg)^2 X_i^{-1} \geq 0
\end{equation}
 for any $a\in\mathcal{H}$. Thus  $\Delta_3-\Delta_6$  is n.n.d. Hence, \textit{a.s.} $[\mathbb{P}]$, the GREG estimator is asymptotically at least as efficient as the RHC estimator under RHC sampling design. Next, it follows from the proofs of Propositions \ref{prop 1} and \ref{prop 2} that \textit{a.s.} $[\mathbb{P}]$, asymptotic covariance operators of $\sqrt{n}(\hat{\overline{Y}}_{HT}-\overline{Y})$ and $\sqrt{n}(\hat{\overline{Y}}_{GREG}-\overline{Y})$ under any HE$\pi$PS sampling design are $\Delta_2$ and  $\Delta_5$,  respectively. Further, it follows from the linear regression model in \eqref{eq 2} in  Section  \ref{sec 4} that 
\begin{align}\label{eq 49}
\begin{split}
&\langle \Delta_2 a,a\rangle=\bigg[ E_{\mathbb{P}}(\tilde{\epsilon}_i)^2 \bigg\{ \mu_x E_{\mathbb{P}}\big(X_i^{2\eta-1}\big)-\lambda E_{\mathbb{P}}\big(X_i^{2\eta}\big)\bigg\}+ E_{\mathbb{P}}\bigg\{\bigg(\tilde{\beta}_0+\sum_{j=1}^d\tilde{\beta}_jZ_{ji}\bigg)^2 \big(X_i^{-1}\mu_x-\lambda\big)\bigg\}\\
&-\chi^{-1}\mu_x^{-1}\bigg\{(1-\lambda)\tilde{\beta}_0 \mu_x +\bigg(\sum_{j=1}^d\tilde{\beta}_j\big(\mu_j\mu_x-\lambda \mu_{jx}\big)\bigg\}^2\bigg]\\
&\text{and }\langle \Delta_5 a,a\rangle=  E_{\mathbb{P}}(\tilde{\epsilon}_i)^2 \bigg(\mu_x E_{\mathbb{P}}\big(X_i^{2\eta-1}\big)-\lambda E_{\mathbb{P}}\big(X_i^{2\eta}\big)\bigg), 
\end{split}
\end{align}
where $\mu_{jx}$=$E_{\mathbb{P}}(Z_{ji}X_i)$ for $j$=$1,\ldots,d$ and  $\chi$=$\mu_x-\lambda \mu_x^{-1} E_{\mathbb{P}}(X_i)^2 $.  Now, since Assumption \ref{ass C2} holds and $0\leq \lambda\leq \mu_x b^{-1}$, we have
\begin{align}\label{eq 51}
\begin{split}
&\langle (\Delta_2-\Delta_5) a,a\rangle=  E_{\mathbb{P}}\bigg[\bigg\{\bigg(\tilde{\beta}_0+\sum_{j=1}^d\tilde{\beta}_jZ_{ji}\bigg)-\chi^{-1} X_i\bigg(\sum_{j=0}^d\tilde{\beta}_j\mu_j-\lambda \tilde{\beta}_0-\\
&\sum_{j=1}^d \lambda\tilde{\beta}_j\mu_{jx}\mu_x^{-1}\bigg)\bigg\}^2\big(X_i^{-1}\mu_x -\lambda\big)\bigg] \geq 0
\end{split}
\end{align}
 Thus using similar arguments as above, we can say that \textit{a.s.} $[\mathbb{P}]$, the GREG estimator is asymptotically at least as efficient as the HT estimator under any HE$\pi$PS sampling design. 
\end{proof}

\begin{proof}[Proof of Theorem \ref{thm 3}]
 Recall from the proofs of Theorems \ref{thm 2} and \ref{thm 1} that \textit{a.s.} $[\mathbb{P}]$, the asymptotic covariance operators of the GREG estimator under SRSWOR, any HE$\pi$PS sampling design and RHC sampling design are $\Delta_4$, $\Delta_5$ and $\Delta_6$, respectively. Also, recall from the proof of Theorem \ref{thm 1} that 
\begin{equation}
\langle \Delta_5 a,a\rangle=E_{\mathbb{P}}(\tilde{\epsilon}_i)^2 \bigg(\mu_x E_{\mathbb{P}}\big(X_i^{2\eta-1}\big)-\lambda E_{\mathbb{P}}\big(X_i^{2\eta}\big)\bigg)\text{ and }\langle \Delta_6 a,a\rangle=c\mu_x E_{\mathbb{P}}(\tilde{\epsilon}_i)^2 E_{\mathbb{P}}\big(X_i^{2\eta-1}\big)
\end{equation}
 for any $a\in\mathcal{H}$ under the linear regression model in \eqref{eq 2} in  Section  \ref{sec 4}.  It can be further shown using \eqref{eq 2} in  Section  \ref{sec 4}  of the main text and ($16$)  in the proof of Lemma \ref{lem 6}  in the supplement (\cite{deychaudhuri2023})  that 
\begin{equation}\label{eq 38}
\langle \Delta_4 a,a \rangle=(1-\lambda)E_{\mathbb{P}}(\tilde{\epsilon}_i)^2 E_{\mathbb{P}}(X_i^{2\eta}) 
\end{equation}
 for any $a\in\mathcal{H}$. Therefore, we have 
\begin{align}\label{eq 52}
\begin{split}
&\langle (\Delta_4-\Delta_5) a,a \rangle=E_{\mathbb{P}}(\tilde{\epsilon}_i)^2 cov_{\mathbb{P}}\bigg(X_i^{2\eta-1},X_i\bigg) \\
&\langle (\Delta_6-\Delta_5) a,a \rangle=E_{\mathbb{P}}(\tilde{\epsilon}_i)^2\bigg(\lambda E_{\mathbb{P}}\big(X_i^{2\eta}\big)-(1-c)E_{\mathbb{P}}\big(X_i^{2\eta-1}\big)\mu_x\bigg)\text{ and }\\
& \langle (\Delta_4-\Delta_6) a,a \rangle=E_{\mathbb{P}}(\tilde{\epsilon}_i)^2\bigg((1-\lambda) E_{\mathbb{P}}\big(X_i^{2\eta}\big)-c E_{\mathbb{P}}\big(X_i^{2\eta-1}\big)\mu_x\bigg)
\end{split}
\end{align}
for any $a\in\mathcal{H}$. Note that $E_{\mathbb{P}}(\tilde{\epsilon}_i)^2$=$\langle E_{\mathbb{P}}(\epsilon_i\otimes \epsilon_i)a,a\rangle>0$ for any $a\in\mathcal{H}$ since $E_{\mathbb{P}}(\epsilon_i\otimes \epsilon_i)$ is p.d. Also, note that $cov_{\mathbb{P}}\big(X_i^{2\eta-1},X_i\big)>0$ for $\eta>0.5$, $cov_{\mathbb{P}}\big(X_i^{2\eta-1},X_i\big)=0$ for $\eta=0.5$ and $cov_{\mathbb{P}}\big(X_i^{2\eta-1},X_i\big)<0$ for $\eta<0.5$. Further, it follows from Lemma \ref{lem 4} that $c$=$1$ for $\lambda$=$0$, $c$=$1-\lambda$ for $\lambda>0$ and $\lambda^{-1}$ being an integer, and $c>1-\lambda$ for $\lambda>0$ and $\lambda^{-1}$ being a non-integer.  Therefore, the results in Table \ref{table 2} below hold, and hence the results stated in Table \ref{table 1} hold. 

\begin{table}[h]

\caption{Relations among $\Delta_4$, $\Delta_5$ and $\Delta_6$}
\label{table 2}
\centering
\begin{tabular}{cccc} 
\hline
& \multirow{2}{*}{$\lambda$=$0$}& $\lambda>0$ \& & $\lambda>0$ \&  \\ 
& & $\lambda^{-1}$ is an integer& $\lambda^{-1}$ is a non-integer\\
\hline 
\multirow{2}{*}{$\eta<0.5$}& $\Delta_5-\Delta_4$ and & $\Delta_5-\Delta_4$ and  & $\Delta_5-\Delta_4$ and \\
& $\Delta_6-\Delta_4$ are p.d. &$\Delta_6-\Delta_4$ are p.d.& $\Delta_6-\Delta_4$ are p.d.\\ 
\hline 
\multirow{2}{*}{$\eta=0.5$}& \multirow{2}{*}{$\Delta_4$=$\Delta_5$=$\Delta_6$}& \multirow{2}{*}{$\Delta_4$=$\Delta_5$=$\Delta_6$} & $\Delta_4$=$\Delta_5$ and \\
&&& $\Delta_6-\Delta_4$ are p.d.\\
\hline 
\multirow{2}{*}{$\eta>0.5$}&$\Delta_5$=$\Delta_6$ and  & $\Delta_4-\Delta_5$ and  & $\Delta_4-\Delta_5$ and \\
& $\Delta_4-\Delta_5$ are p.d.& $\Delta_6-\Delta_5$ are p.d.& $\Delta_6-\Delta_5$ are p.d.\\
\hline
\end{tabular}
\end{table}
\par
\vspace{.2cm}

Next, if we put $\lambda$=$0$ and $c$=$1$, respectively, in the expressions of $\Delta_5$ and $\Delta_6$ in the proof of Lemma \ref{lem 6}  in the supplement (\cite{deychaudhuri2023}),   we have $\Delta_5$=$\Delta_6$. Thus \textit{a.s.} $[\mathbb{P}]$, the GREG estimator has the same asymptotic covariance operator under RHC and any HE$\pi$PS sampling designs. Hence, \textit{a.s.} $[\mathbb{P}]$, the GREG estimator has the same asymptotic distribution under RHC and any HE$\pi$PS sampling designs.  This completes the proof of the theorem.
\end{proof}
\begin{proof}[Proof of Theorem \ref{prop 4}] Recall the expression of $\hat{\Sigma}$ from \eqref{eq 33} in  Section  \ref{sec 4} and note that 
\begin{equation}\label{eq 39}
\hat{\Sigma}=(nN^{-2})\bigg(\sum_{i\in s}(\hat{V}_i\otimes \hat{V}_i)(\pi_i^{-1}-1)\pi_i^{-1}-\sum_{i\in s}(1-\pi_i)\hat{T}\otimes\hat{T}\bigg) 
\end{equation}
 with $\hat{T}$=$\sum_{i\in s}\hat{V}_i(\pi_i^{-1}-1)\big(\sum_{i\in s}(1-\pi_i)\big)^{-1}$. Let us first consider the case, when $\Sigma$ denotes the asymptotic covariance operator of $\sqrt{n}(\hat{\overline{Y}}_{HT}-\overline{Y})$ and $\hat{\Sigma}$ is its estimator. Then, we have $\hat{V}_i$=$Y_i$ in $\hat{\Sigma}$. Now, recall the expression of $\Gamma_1$ from the paragraph preceding Lemma \ref{lem 6} and note that 
\begin{equation}\label{eq 40}
\Gamma_1=(nN^{-2})\bigg(\sum_{i=1}^N(V_i\otimes V_i)(\pi_i^{-1}-1)-\sum_{i=1}^N\pi_i(1-\pi_i)T\otimes T\bigg) 
\end{equation}
 with $T$=$\sum_{i=1}^N V_i(1-\pi_i)\big(\sum_{i=1}^N \pi_i(1-\pi_i)\big)^{-1}$. Let us substitute $V_i$=$Y_i$ in $\Gamma_1$. We shall first show that under SRSWOR, LMS sampling design and any HE$\pi$PS sampling design, $\hat{\Sigma}-\Gamma_1\xrightarrow{p} 0$ with respect to the HS norm as $\nu\rightarrow\infty$ \textit{a.s.} $[\mathbb{P}]$. It follows by Assumption \ref{ass C3} that  $N^{-1}\sum_{i=1}^N||Y_i||_{\mathcal{H}}^2$=$O(1)$   as $\nu\rightarrow\infty$ \textit{a.s.} $[\mathbb{P}]$. It also follows by \eqref{eq 1} in the statement of Lemma \ref{lem 5} that as $\nu\rightarrow\infty$, $N^{-1}\sum_{i=1}^N\big(n^{-1}N\pi_i(1-\pi_i)\big)^2$=$O(1)$ under the above sampling designs \textit{a.s.} $[\mathbb{P}]$. Then, using the same line of arguments as in the proof of Lemma \ref{lem 3}  in the supplement (\cite{deychaudhuri2023}),  it can be shown that 
\begin{align}\label{eq 41}
\begin{split}
&n^{-1}\bigg(\sum_{i\in s}(1-\pi_i)-\sum_{i=1}^N\pi_i(1-\pi_i)\bigg)=o_p(1)\text{ and }\\
&N^{-1}\bigg|\bigg|\sum_{i\in s}\hat{V}_i(\pi_i^{-1}-1)-\sum_{i=1}^N V_i(1-\pi_i)\bigg|\bigg|_{\mathcal{H}}=o_p(1) 
\end{split}
\end{align}
 as $\nu\rightarrow\infty$ \textit{a.s.} $[\mathbb{P}]$. Moreover, $n^{-1}\sum_{i=1}^N\pi_i(1-\pi_i)$ is bounded away from $0$ as $\nu\rightarrow\infty$ \textit{a.s.} $[\mathbb{P}]$ because \eqref{eq 1} in Lemma \ref{lem 5} and Assumption \ref{ass C1} hold. Consequently, under all of the above-mentioned sampling designs, $(nN^{-2})\big(\sum_{i\in s}(1-\pi_i)(\hat{T}\otimes\hat{T})-\sum_{i=1}^N\pi_i(1-\pi_i)(T\otimes T)\big)\xrightarrow{p}0$ with respect to the HS norm as $\nu\rightarrow\infty$ \textit{a.s.} $[\mathbb{P}]$. Similarly, $(nN^{-2})\big(\sum_{i\in s}(\hat{V}_i\otimes \hat{V}_i)(\pi_i^{-1}-1)\pi_i^{-1}-\sum_{i=1}^N(V_i\otimes V_i)(\pi_i^{-1}-1)\big)\xrightarrow{p} 0$ with respect to the HS norm as $\nu\rightarrow\infty$ \textit{a.s.} $[\mathbb{P}]$. Thus under the above sampling designs, $\hat{\Sigma}-\Gamma_1\xrightarrow{p} 0$ with respect to the HS norm as $\nu\rightarrow\infty$ \textit{a.s.} $[\mathbb{P}]$. Recall from  Section  \ref{sec 4} that $\Sigma$=$\lim_{\nu\rightarrow\infty}\Gamma_1$ \textit{a.s.} $[\mathbb{P}]$. Therefore, under the aforesaid sampling designs, $\hat{\Sigma}\xrightarrow{p} \Sigma$ with respect to the HS norm as $\nu\rightarrow\infty$ \textit{a.s.} $[\mathbb{P}]$.
\par
\vspace{.2cm}

Let us next consider the case, when $\Sigma$ denotes the asymptotic covariance operator of $\sqrt{n}(\hat{\overline{Y}}_{GREG}-\overline{Y})$ and $\hat{\Sigma}$ denotes its estimator. Then, $\hat{\Sigma}$ is the same as described in the preceding paragraph with $\hat{V}_i$=$Y_i-\hat{\overline{Y}}_{HT}-\hat{S}_{zy}((Z_i-\hat{\overline{Z}}_{HT})\hat{S}_{zz}^{-1})$. Let us also consider $\Gamma_1$ with $V_i$=$Y_i-\overline{Y}-S_{zy}((Z_i-\overline{Z}) S_{zz}^{-1})$. Note that 
\begin{align}\label{eq 42}
\begin{split}
&N^{-1}\bigg(\sum_{i\in s}\hat{V}_i(\pi_i^{-1}-1)-\sum_{i=1}^N V_i(1-\pi_i)\bigg)=N^{-1}\sum_{i\in s}(\hat{V}_i-V_i)(\pi_i^{-1}-1)+ \\
&N^{-1}\bigg(\sum_{i\in s}V_i(\pi_i^{-1}-1)-\sum_{i=1}^N V_i(1-\pi_i)\bigg).
\end{split}
\end{align} 
 It can be shown in the same way as in the proof of Lemma \ref{lem 3}  in the supplement (\cite{deychaudhuri2023})  that  $||N^{-1}(\sum_{i\in s}V_i\times$ $(\pi_i^{-1}-1)-\sum_{i=1}^N V_i(1-\pi_i))||_{\mathcal{H}}$=$o_p(1)$  under the sampling designs considered in the previous paragraph as $\nu\rightarrow\infty$ \textit{a.s.} $[\mathbb{P}]$. Further, it can be shown that  $||N^{-1}\sum_{i\in s}(\hat{V}_i-V_i)(\pi_i^{-1}-1)||_{\mathcal{H}}$=$o_p(1)$  as $\nu\rightarrow\infty$ \textit{a.s.} $[\mathbb{P}]$ since  $||\hat{\overline{Y}}_{HT}-\overline{Y}||_{\mathcal{H}}$=$o_p(1)$, $||\hat{S}_{zy}-S_{zy}||_{op}$=$o_p(1)$, $||\hat{S}^{-1}_{zz}-S^{-1}_{zz}||_{op}$=$o_p(1)$, $||\hat{S}_{zy}||_{op}$=$O_p(1)$ and $||S_{zz}^{-1}||_{op}$=$O(1)$  as $\nu\rightarrow\infty$ \textit{a.s.} $[\mathbb{P}]$ (see the proof of Proposition \ref{prop 2}). Then, $(nN^{-2})(\sum_{i\in s}(1-\pi_i)(\hat{T}\otimes\hat{T})-\sum_{i=1}^N\pi_i(1-\pi_i)(T\otimes T))\xrightarrow{p} 0$ with respect to the HS norm as $\nu\rightarrow\infty$ \textit{a.s.} $[\mathbb{P}]$. Similarly, $(nN^{-2})(\sum_{i\in s}(\hat{V}_i\otimes \hat{V}_i)(\pi_i^{-1}-1)\pi_i^{-1}-\sum_{i=1}^N(V_i\otimes V_i)(\pi_i^{-1}-1))\xrightarrow{p} 0$ with respect to the HS norm as $\nu\rightarrow\infty$ \textit{a.s.} $[\mathbb{P}]$. Hence, under the above sampling designs, $\hat{\Sigma}-\Gamma_1\xrightarrow{p} 0$, and hence $\hat{\Sigma}\xrightarrow{p} \Sigma$ with respect to the HS norm as $\nu\rightarrow\infty$ \textit{a.s.} $[\mathbb{P}]$. 
\par
\vspace{.2cm}

Next, consider the case, when $\Sigma$ denotes the asymptotic covariance operator of $\sqrt{n}(\hat{\overline{Y}}_{RHC}-\overline{Y})$ or $\sqrt{n}(\hat{\overline{Y}}_{GREG}-\overline{Y})$ under RHC sampling design, and $\hat{\Sigma}$ denotes its estimator. Recall from \eqref{eq 34} in  Section  \ref{sec 4} that in this case, 
\begin{align}\label{eq 43}
\begin{split}
&\hat{\Sigma}=n\gamma(\overline{X}N^{-1})\sum_{i\in s}\bigg(\hat{V}_i-X_i\hat{\overline{V}}_{RHC} \overline{X}^{-1}\bigg)\otimes\bigg(\hat{V}_i-X_i\hat{\overline{V}}_{RHC}\overline{X}^{-1}\bigg)(Q_i X_i^{-2})=\\
&n\gamma\bigg((\overline{X}N^{-1})\sum_{i\in s}(\hat{V}_i\otimes\hat{V}_i)Q_i X_i^{-2}-\hat{\overline{V}}_{RHC}\otimes\hat{\overline{V}}_{RHC}\bigg). 
\end{split}
\end{align}
 Also, recall the expression of $\Gamma_2$ from the paragraph preceding Lemma \ref{lem 6} and note that 
\begin{equation}\label{eq 44}
 \Gamma_2=n\gamma\bigg((\overline{X}N^{-1})\sum_{i=1}^N(V_i\otimes V_i)X_i^{-1}-\overline{V}\otimes \overline{V}\bigg). 
\end{equation} 
 Then, it can be shown in a similar way as in the earlier cases that under RHC sampling design, $\hat{\Sigma}-\Gamma_2\xrightarrow{p} 0$ with respect to the HS norm as $\nu\rightarrow\infty$ \textit{a.s.} $[\mathbb{P}]$. Therefore, under RHC sampling design, $\hat{\Sigma}\xrightarrow{p} \Sigma$ with respect to the HS norm as $\nu\rightarrow\infty$ \textit{a.s.} $[\mathbb{P}]$ because $\Sigma$=$\lim_{\nu\rightarrow\infty}\Gamma_2$ \textit{a.s.} $[\mathbb{P}]$ (see  Section  \ref{sec 4}).
\end{proof}                                                                                                                                                                                                                       
\end{appendix}

%%%%%%%%%%%%%%%%%%%%%%%%%%%%%%%%%%%%%%%%%%%%%%
%% Support information, if any,             %%
%% should be provided in the                %%
%% Acknowledgements section.                %%
%%%%%%%%%%%%%%%%%%%%%%%%%%%%%%%%%%%%%%%%%%%%%%
\section*{Acknowledgments}
The authors would like to thank Irish Social Science Data Archive (ISSDA) and its administrator, Ruby O'Riordan, for making Electricity Customer Behaviour Trial data available to the authors. The authors gratefully acknowledge careful reading of an earlier version of the paper by two anonymous reviewers. Critical comments and constructive suggestions from the reviewers led to significant improvement of the paper. The authors thank the editor for allowing additional space for presenting several technical details in the paper.
%%%%%%%%%%%%%%%%%%%%%%%%%%%%%%%%%%%%%%%%%%%%%%
%% Funding information, if any,             %%
%% should be provided in the                %%
%% funding section.                         %%
%%%%%%%%%%%%%%%%%%%%%%%%%%%%%%%%%%%%%%%%%%%%%%

%%%%%%%%%%%%%%%%%%%%%%%%%%%%%%%%%%%%%%%%%%%%%%
%% Supplementary Material, including data   %%
%% sets and code, should be provided in     %%
%% {supplement} environment with title      %%
%% and short description. It cannot be      %%
%% available exclusively as external link.  %%
%% All Supplementary Material must be       %%
%% available to the reader on Project       %%
%% Euclid with the published article.       %%
%%%%%%%%%%%%%%%%%%%%%%%%%%%%%%%%%%%%%%%%%%%%%%

\section*{Supplementary material}
In this supplement, we give the proofs of Lemmas \ref{lem 5}--\ref{lem 3}.

%%%%%%%%%%%%%%%%%%%%%%%%%%%%%%%%%%%%%%%%%%%%%%%%%%%%%%%%%%%%%
%%                  The Bibliography                       %%
%%                                                         %%
%%  imsart-???.bst  will be used to                        %%
%%  create a .BBL file for submission.                     %%
%%                                                         %%
%%  Note that the displayed Bibliography will not          %%
%%  necessarily be rendered by Latex exactly as specified  %%
%%  in the online Instructions for Authors.                %%
%%                                                         %%
%%  MR numbers will be added by VTeX.                      %%
%%                                                         %%
%%  Use \cite{...} to cite references in text.             %%
%%                                                         %%
%%%%%%%%%%%%%%%%%%%%%%%%%%%%%%%%%%%%%%%%%%%%%%%%%%%%%%%%%%%%%

%% if your bibliography is in bibtex format, uncomment commands:
%\bibliographystyle{imsart-nameyear} % Style BST file (imsart-number.bst or imsart-nameyear.bst)
%\bibliography{bibliography}       % Bibliography file (usually '*.bib')

%% or include bibliography directly:

\newpage
%%%%%%%%%%%%%%%%%%%%%%%%%%%%%%%%%%%%%%%%%%%%%%%%%%%%%%%%%%%%%%%%%%%%%%%%%%%%%%%%%%%%%%%%%%%%%%%%%%%%%%%%%%%%%%%%%%%%%%%%%%%%
\vskip .55cm
\noindent
Anurag Dey\\
{\it Indian Statistical Institute, Kolkata}
\vskip 1.3pt
\noindent
E-mail: deyanuragsaltlake64@gmail.com
\vskip 1.3pt

\noindent
Probal Chaudhuri\\
{\it Indian Statistical Institute, Kolkata}
\vskip 1.3pt
\noindent
E-mail:probalchaudhuri@gmail.com
% \vskip .3cm
%\centerline{(Received ???? 20??; accepted ???? 20??)}\par
%\end{document}
%%%%%%%%%%%%%%%%%%%%%%%%%%%%%%%%%%%%%%%%%%%%%%%%%%%%%%%%%%%%%%%%%%%%%%%%%%%%%%%%%%%%%%%%%%%%%%%%%%%%%%%%%%%%%%%%%%%%%%%%%%%%
%%%%%%%%%%%%%%%%%%%%%%%%%%%%%%%%%%%%%%%%%%%%%%%%%%%%%%%%%%%%%%%%%%%%%%%%%%%%%%%%%%%%%%%%%%%%%%%%%%%%%%%%%%%%%%%%%%%%%%%%%%%%
\end{document}

%%%%%%%%%%%%%%%%%%%%%%%%%%%%%%%%%%%%%%%%%%%%%%%%%%%%%%%%%%%%%%%%%%%%%%%%%%%%%%%%%%%%%%%%%%%%%%%%%%%%%%%%%%%%%%%%%%%%%%%%%%%%
%%%%%%%%%%%%%%%%%%%%%%%%%%%%%%%%%%%%%%%%%%%%%%%%%%%%%%%%%%%%%%%%%%%%%%%%%%%%%%%%%%%%%%%%%%%%%%%%%%%%%%%%%%%%%%%%%%%%%%%%%%%%

%\renewcommand{\baselinestretch}{2}

\title{\textbf{Supplementary material for ``On estimators of the mean of infinite dimensional data in finite populations"}}
%\runtitle{Quantile processes in finite populations}
\author{Anurag Dey and Probal Chaudhuri\\
\textit{Indian Statistical Institute, Kolkata}}
\maketitle

%%%%%%%%%%%%%%%%%%%%%%%%%%%%%%%%%%%%%%%%%%%%%%%%%%%%%%%%%%%%%%%%%%%%%%%%%%%%%%%%%%%%%%%%%%%%%%%%%%%%%%%%%%%%%%%%%%%%%%%%%%%%

\begin{abstract}
In this supplement, we give the proofs of Lemmas $1$--$6$, which are stated in the Appendix B of the main text.\\ 
\end{abstract}

\textbf{Keywords and phrases:} Asymptotic normality, Consistency of estimators, Covariance operator, Heteroscedasticity, High entropy sampling design, Inclusion probability, Relative efficiency, Separable Hilbert space.

\section{Proofs of Lemmas $1$--$6$}
\begin{proof}[\textbf{Proof of Lemma $1$}]
 Suppose that $P(s,\omega)$ and $R(s,\omega)$ denote LMS sampling design and SRSWOR, respectively. Note that SRSWOR is a rejective sampling design. Then, $P(s,\omega)$= $\overline{x}(^NC_n \overline{X})^{-1}$  and $R(s,\omega)$=$(^NC_n)^{-1}$, where  $\overline{x}$=$n^{-1}\sum_{i\in s}X_i$ and $\overline{X}$=$N^{-1}\sum_{i=1}^N X_i$. Moreover, by Cauchy-Schwarz inequality and Assumption  $2$,  we have 
\begin{equation}\label{eq 7}
D(P||R)=E_{R}\bigg\{ \overline{x}\overline{X}^{-1}\log\bigg(\overline{x}\overline{X}^{-1}\bigg)\bigg\}\leq K E_{R}\bigg|\overline{x}\overline{X}^{-1}-1\bigg|\leq K E_{R}\bigg(\overline{x}\overline{X}^{-1}-1\bigg)^2
\end{equation} 
for some constant $K>0$ \textit{a.s.} $[\mathbb{P}]$ since $\log(x)\leq|x-1|$ for any $x>0$.  Here, $E_{R}$ denotes the expectation with respect to SRSWOR.  Moreover, by Assumption   $2$,  we have
\begin{align}\label{eq 8}
\begin{split}
&nE_{R}\bigg(\overline{x} \overline{X}^{-1}-1\bigg)^2=\bigg(1-nN^{-1}\bigg)\bigg(N(N-1)^{-1}\bigg) \bigg(S^2_x \overline{X}^{-2}\bigg)\leq 2 \sum_{i=1}^N X_i^2\bigg(N\overline{X}^2\bigg)^{-1} \\
&\leq 2\bigg\{\bigg(\max_{1\leq i \leq N} X_{i}\bigg)\bigg(\min_{1\leq i \leq N} X_{i}\bigg)^{-1}\bigg\}^2=O(1) 
\end{split}
\end{align}
as $\nu\rightarrow\infty$ \textit{a.s.} $[\mathbb{P}]$, where $S_x^2$=$N^{-1}\sum_{i=1}^N X_i^2-\overline{X}^2$.  Hence, $D(P||R)\rightarrow 0$ as $\nu\rightarrow\infty$ \textit{a.s.} $[\mathbb{P}]$ by \eqref{eq 7} and \eqref{eq 8}. Thus LMS sampling design is a high entropy sampling design
\par
\vspace{.1cm}

Next, note that  ($10$) in Lemma $1$  holds  trivially under SRSWOR. Now, let $\{\pi_i\}_{i=1}^N$ be the inclusion probabilities of LMS sampling design. Then,   $\pi_i$=$(n-1)(N-1)^{-1}+\big(X_i(\sum_{i=1}^N X_i)^{-1}\big)\big((N-n)(N-1)^{-1}\big)$ and $\pi_i-nN^{-1}$=$-(N-n)(N(N-1))^{-1}(X_i \overline{X}^{-1}-1)$. Further, we have 
\begin{align}\label{eq 9}
\begin{split}
&n^{-1}N\bigg|\pi_i-n N^{-1}\bigg|=(N-n)(n(N-1))^{-1}\bigg| X_i \overline{X}^{-1}-1\bigg|\leq(N-n)(n(N-1))^{-1}\times\\
&\bigg[1+\bigg(\max_{1\leq i \leq N} X_{i}\bigg)\bigg(\min_{1\leq i \leq N} X_{i}\bigg)^{-1}\bigg].
\end{split}
\end{align} 
 Therefore,  $\lim_{\nu\rightarrow\infty}\max_{1\leq i \leq N}|n^{-1}N\pi_{i}-1|$=$0$  \textit{a.s.} $[\mathbb{P}]$ by Assumption   $2$  and \eqref{eq 9}. Hence,   ($10$) in Lemma $1$  holds under LMS sampling design. Further,  ($10$) in Lemma $1$  holds under any HE$\pi$PS sampling design since Assumption  $2$  holds. 
\end{proof}
\begin{proof}[\textbf{Proof of Lemma $2$}]
Let us first consider the case of $\lambda$=$0$. Note that 
\begin{equation}\label{eq 55}
n (Nn^{-1}-1)(N-n)(N(N-1))^{-1}\leq n\gamma\leq n (Nn^{-1}+1)(N-n)(N(N-1))^{-1}
\end{equation} 
by  ($4$) in Section $3$ of the main text.   Moreover, $n (Nn^{-1}+1)(N-n)(N(N-1))^{-1}$=$(1+nN^{-1})(N-n)(N-1)^{-1}\rightarrow 1$ and $n (Nn^{-1}-1)(N-n)(N(N-1))^{-1}$=$(1-nN^{-1})(N-n)(N-1)^{-1}\rightarrow 1$ as $\nu\rightarrow\infty$ because Assumption  $1$  holds and $\lambda$=$0$. Thus we have $n\gamma\rightarrow 1$ as $\nu\rightarrow\infty$ in this case.
\par
\vspace{.2cm}

Next, consider the case, when $\lambda>0$ and $\lambda^{-1}$ is an integer. Here, we consider the following sub-cases. Let us first consider the sub-case, when $Nn^{-1}$ is an integer for all sufficiently large $\nu$. Then, by  ($4$) in the main text,  we have $n\gamma$=$(N-n)(N-1)^{-1}$ for all sufficiently large $\nu$. Now, since Assumption  $1$ holds,  we have 
\begin{equation}\label{eq 53}
(N-n)(N-1)^{-1}\rightarrow 1-\lambda\text{ as }\nu\rightarrow\infty.
\end{equation}
\par

Further, consider the sub-case, when $Nn^{-1}$ is a non-integer and $Nn^{-1}-\lambda^{-1}\geq 0$ for all sufficiently large $\nu$. Then by  ($4$) in the main article,  we have 
\begin{equation}\label{eq 57}
n\gamma=N(N-1)^{-1}nN^{-1} \lfloor Nn^{-1}\rfloor \big(2-\big(nN^{-1} \lfloor Nn^{-1} \rfloor \big)-nN^{-1}\big)
\end{equation}
for all sufficiently large $\nu$. Now, since Assumption  $1$  holds, we have $0\leq Nn^{-1}-\lambda^{-1}<1$ for all sufficiently large $\nu$. Then, $\lfloor Nn^{-1} \rfloor$=$\lambda^{-1}$ for all sufficiently large $\nu$, and hence
\begin{equation}\label{eq 54}
N(N-1)^{-1} nN^{-1} \lfloor Nn^{-1} \rfloor \bigg(2-\big(nN^{-1} \lfloor Nn^{-1}\rfloor\big)-nN^{-1}\bigg)\rightarrow 1-\lambda\text{ as }\nu\rightarrow\infty.
\end{equation}
\par

 Next,  consider the sub-case, when $Nn^{-1}$ is a non-integer and $Nn^{-1}-\lambda^{-1}< 0$ for all sufficiently large $\nu$. Then, the result in \eqref{eq 57} above holds by  ($4$) in the main text,  and $-1\leq Nn^{-1}-\lambda^{-1}<0$ for all sufficiently large $\nu$ by Assumption  $1$.  Therefore, $\lfloor Nn^{-1} \rfloor$=$\lambda^{-1}-1$ for all sufficiently large $\nu$, and hence the result in \eqref{eq 54} above holds. Thus, in the case of $\lambda>0$ and $\lambda^{-1}$ being an integer, $n\gamma$ converges to $1-\lambda$ as $\nu\rightarrow\infty$ through all the sub-sequences, and hence $n\gamma\rightarrow 1-\lambda$ as $\nu\rightarrow\infty$. Thus we have $c$=$1-\lambda$ in this case.
\par
\vspace{.2cm}

Finally, consider the case, when $\lambda>0$, and $\lambda^{-1}$ is a non-integer. Then, $Nn^{-1}$ must be a non-integer for all sufficiently large $\nu$, and hence $n\gamma$=$N(N-1)^{-1} nN^{-1} \lfloor Nn^{-1} \rfloor \big(2-\big(nN^{-1} \lfloor Nn^{-1} \rfloor \big)-nN^{-1}\big)$ for all sufficiently large $\nu$ by  ($4$) in the main paper.  Note that in this case, $Nn^{-1}-\lfloor \lambda^{-1} \rfloor\rightarrow \lambda^{-1}-\lfloor \lambda^{-1} \rfloor\in(0,1)$ as $\nu\rightarrow\infty$ by Assumption  $1$.  Therefore, $\lfloor \lambda^{-1} \rfloor<Nn^{-1}<\lfloor \lambda^{-1} \rfloor+1$ for all sufficiently large $\nu$, and hence $\lfloor N n^{-1} \rfloor$=$\lfloor \lambda^{-1} \rfloor$ for all sufficiently large $\nu$. Thus $n\gamma\rightarrow \lambda \lfloor \lambda^{-1} \rfloor (2-\lambda \lfloor \lambda^{-1} \rfloor-\lambda)$ as $\nu\rightarrow\infty$ by Assumption  $1$.  Now, if $m$=$ \lfloor \lambda^{-1} \rfloor$ and $\lambda^{-1}$ is a non-integer, then $(m+1)^{-1}<\lambda<m^{-1}$. Therefore, $\lambda \lfloor \lambda^{-1}\rfloor (2-\lambda \lfloor\lambda^{-1}\rfloor-\lambda)-1+\lambda$=$-\big(1-(2m+1)\lambda+m(m+1)\lambda^2\big)$=$-(1-m\lambda)(1-(m+1)\lambda)>0$. Thus we have $c$=$\lambda \lfloor \lambda^{-1} \rfloor (2-\lambda \lfloor \lambda^{-1} \rfloor-\lambda)>1-\lambda$ in this case. This completes the proof of the Lemma.
\end{proof}
\begin{proof}[\textbf{Proof of Lemma $3$}]
Let us first consider the case $V_i$=$Y_i$ for $i$=$1,\ldots,N$. Then, we have 

\begin{align}\label{eq 10}
\begin{split}
&\Gamma_1=nN^{-2}\sum_{i=1}^N(V_i-T\pi_i)\otimes(V_i-T\pi_i)(\pi_i^{-1}-1)=nN^{-2}\bigg\{\sum_{i=1}^N (\pi_i^{-1}-1) (Y_i\otimes Y_i)-\\
&\bigg(\sum_{i=1}^N Y_i(1-\pi_i)\otimes \sum_{i=1}^N Y_i(1-\pi_i )\bigg)\bigg(\sum_{i=1}^N\pi_i(1-\pi_i)\bigg)^{-1}\bigg\}. 
\end{split}
\end{align}
 Now,  substituting $\pi_i$=$nN^{-1}$  for SRSWOR, we obtain $\Gamma_1$=$(1-nN^{-1})N^{-1}\sum_{i=1}^N (Y_i-\overline{Y})\otimes (Y_i-\overline{Y})$. Note that  $E_{\mathbb{P}}||Y_i||_{\mathcal{H}}^2<\infty$  in view of Assumption  $3$.  Then, under SRSWOR, 
\begin{equation}\label{eq 47}
\Gamma_1\rightarrow \Delta_1=(1-\lambda)E_{\mathbb{P}}(Y_i-E_{\mathbb{P}}(Y_i))\otimes (Y_i-E_{\mathbb{P}}(Y_i))
\end{equation}
 with respect to the HS norm as $\nu\rightarrow\infty$ \textit{a.s.} $[\mathbb{P}]$ by strong law of large numbers (SLLN) and Assumption  $1$.  Now, suppose that $\Gamma_1^{(1)}$ and $\Gamma_1^{(2)}$ denote $\Gamma_1$ under SRSWOR and LMS sampling design, respectively. Further, suppose that $\{\pi_i\}_{i=1}^N$ are the inclusion probabilities of LMS sampling design.  Then, we have
\begin{align}\label{eq 11}
\begin{split}
&\Gamma_1^{(2)}-\Gamma_1^{(1)}=nN^{-2}\bigg\{\sum_{i=1}^N (\pi_i^{-1}-n^{-1}N)(Y_i\otimes Y_i)\bigg\}-nN^{-2}\bigg\{\bigg(\sum_{i=1}^N Y_i(1-\pi_i)\otimes \sum_{i=1}^N Y_i(1-\pi_i)\bigg)\\
&\times\bigg(\sum_{i=1}^N\pi_i (1-\pi_i)\bigg)^{-1}-\bigg(\sum_{i=1}^N Y_i(1-nN^{-1})\otimes \sum_{i=1}^N Y_i(1-nN^{-1})\bigg)\bigg(n(1-nN^{-1})\bigg)^{-1}\bigg\}
\end{split}
\end{align}
by \eqref{eq 10}.  Further, it follows from the proof of Lemma $1$  that as $\nu\rightarrow\infty$, $\max_{1\leq i \leq N} |n^{-1}N\pi_{i}-1|\rightarrow0$ \textit{a.s.} $[\mathbb{P}]$. It also follows from  Assumption  $3$  that  $N^{-1}\sum_{i=1}^N ||Y_i||_{\mathcal{H}}^2$=$O(1)$  as $\nu\rightarrow\infty$ \textit{a.s.} $[\mathbb{P}]$. Therefore, it can be shown that as $\nu\rightarrow\infty$, 
\begin{equation}\label{eq 45}
nN^{-2}\bigg\{\sum_{i=1}^N (\pi_i^{-1}-n^{-1}N)(Y_i\otimes Y_i)\bigg\}\rightarrow 0\text{ and } 
\end{equation} 
\begin{align}\label{eq 46}
\begin{split}
&nN^{-2}\bigg\{\bigg(\sum_{i=1}^N Y_i(1-\pi_i)\otimes \sum_{i=1}^N Y_i(1-\pi_i) \bigg)\bigg(\sum_{i=1}^N \pi_i(1-\pi_i)\bigg)^{-1}-\\
&\bigg(\sum_{i=1}^N Y_i(1-nN^{-1})\otimes\sum_{i=1}^N Y_i(1-nN^{-1})\bigg)\bigg(n\big(1-nN^{-1}\big)\bigg)^{-1}\bigg\}\rightarrow 0,\text{ and hence}
\end{split}
\end{align} 
 $\Gamma_1^{(2)}-\Gamma_1^{(1)}\rightarrow 0$ with respect to the HS norm \textit{a.s.} $[\mathbb{P}]$ by \eqref{eq 11}. Thus $\Gamma_1\rightarrow\Sigma_1$ as $\nu\rightarrow\infty$ under SRSWOR as well as under LMS sampling design \textit{a.s.} $[\mathbb{P}]$ with $\Sigma_1$=$\Delta_1$. Next, under any HE$\pi$PS sampling design (i.e., a sampling design with $\pi_i$=$nX_i(\sum_{i=1}^N X_i)^{-1}$),  
\begin{align}\label{eq 50}
\begin{split}
&\Gamma_1\rightarrow \Delta_2=E_{\mathbb{P}}\bigg[\bigg\{Y_i-\chi^{-1} X_i\bigg(E_{\mathbb{P}}(Y_i)-\lambda E_{\mathbb{P}}(X_i Y_i)(E_{\mathbb{P}}(X_i))^{-1}\bigg)\bigg\}\otimes\\
&\bigg\{Y_i-\chi^{-1} X_i \bigg(E_{\mathbb{P}}(Y_i)-\lambda E_{\mathbb{P}}(X_i Y_i)(E_{\mathbb{P}}(X_i))^{-1}\bigg)\bigg\}\bigg\{X_i^{-1}E_{\mathbb{P}}(X_i)-\lambda\bigg\} \bigg]
\end{split}
\end{align}
 with respect to the HS norm as $\nu\rightarrow\infty$ \textit{a.s.} $[\mathbb{P}]$ by SLLN because  $E_{\mathbb{P}}||Y_i||_{\mathcal{H}}^2<\infty$,  Assumptions  $1$ and $2$  hold.  Here, $\chi$=$E_{\mathbb{P}}(X_i)-\lambda E_{\mathbb{P}}(X_i)^2 (E_{\mathbb{P}}(X_i))^{-1}$. Note that $\Delta_2$ is a n.n.d. HS operator since Assumption  $1$  holds with $0\leq \lambda<E_{\mathbb{P}}(X_i)b^{-1}$.  Thus as $\nu\rightarrow\infty$, $\Gamma_1\rightarrow\Sigma_1$ under any HE$\pi$PS sampling design \textit{a.s.} $[\mathbb{P}]$ with $\Sigma_1$=$\Delta_2$. Next, note that  $\sum_{j=1}^{\infty}\langle \Delta_1 e_j,e_j\rangle$=$E_{\mathbb{P}}||Y_i-E_{\mathbb{P}}(Y_i)||_{\mathcal{H}}^2<\infty$  and  $\sum_{j=1}^{\infty}\langle \Delta_2 e_j,e_j\rangle$=$E_{\mathbb{P}}\big[||Y_i-\chi^{-1} X_i\{E_{\mathbb{P}}(Y_i)-\lambda E_{\mathbb{P}}(X_i Y_i)(E_{\mathbb{P}}(X_i))^{-1}\}||_{\mathcal{H}}^2\{X_i^{-1}E_{\mathbb{P}}(X_i)-\lambda\}\big]<\infty$  since Assumption  $2$  holds, and  $E_{\mathbb{P}}||Y_i||_{\mathcal{H}}^2<\infty$.  Then, it can be shown in the same way as argued above that as $\nu\rightarrow\infty$, $\sum_{j=1}^{\infty}\langle \Gamma_1 e_j,e_j\rangle$= $nN^{-2}\big\{\sum_{i=1}^N (\pi_i^{-1}-1)||Y_i||_{\mathcal{H}}^2-\sum_{i=1}^N ||Y_i(1-\pi_i)||_{\mathcal{H}}^2\big(\sum_{i=1}^N\pi_i(1-\pi_i)\big)^{-1}\big\}\rightarrow\sum_{j=1}^{\infty} \langle\Delta_1 e_j,e_j\rangle$  under SRSWOR and LMS sampling design, and $\sum_{j=1}^{\infty}\langle \Gamma_1 e_j,e_j\rangle\rightarrow\sum_{j=1}^{\infty}\langle \Delta_2 e_j,e_j\rangle$ under any HE$\pi$PS sampling design \textit{a.s.} $[\mathbb{P}]$. 
\par
\vspace{.1cm}

Next, consider the case of RHC sampling design and $\Gamma_2$ with $V_i$=$Y_i$.  Then, we have 
\begin{align}\label{eq 13}
\begin{split}
&\Gamma_{2}=n\gamma\overline{X}N^{-1}\sum_{i=1}^N\bigg(V_i-X_i\overline{V}\overline{X}^{-1}\bigg)\otimes\bigg(V_i-X_i\overline{V}\overline{X}^{-1}\bigg)X_i^{-1}\\
&=n\gamma\bigg\{\overline{X}N^{-1}\sum_{i=1}^N (Y_i\otimes Y_i)X_i^{-1}-\overline{Y}\otimes\overline{Y}\bigg\}.
\end{split} 
\end{align}
 Note that $n\gamma\rightarrow c>0$ as $\nu\rightarrow\infty$ by Lemma $2$.  Then, by SLLN,  
\begin{equation}\label{eq 48}
\Gamma_{2}\rightarrow \Delta_3=c \big\{E_{\mathbb{P}}(X_i)E_{\mathbb{P}}\big((Y_i\otimes Y_i)X_i^{-1}\big)-E_{\mathbb{P}}(Y_i)\otimes E_{\mathbb{P}}(Y_i)\big\}
\end{equation}
  with respect to the HS norm as $\nu\rightarrow \infty$ \textit{a.s.} $[\mathbb{P}]$. Thus $\Sigma_2$=$\Delta_3$ in this case. It follows that  $\sum_{j=1}^{\infty}\langle \Delta_3 e_j,e_j\rangle$ =$c\big\{E_{\mathbb{P}}(X_i)E_{\mathbb{P}}(||Y_i||_{\mathcal{H}}^2 X_i^{-1})-||E_{\mathbb{P}}(Y_i)||_{\mathcal{H}}^2\big\}<\infty$  since Assumption  $2$  holds, and  $E_{\mathbb{P}}||Y_i||_{\mathcal{H}}^2<\infty$.  Further, it can be shown using SLLN that  $\sum_{j=1}^{\infty}\langle \Gamma_2 e_j,e_j\rangle$=$n\gamma\big\{\overline{X}N^{-1}\sum_{i=1}^N ||Y_i||_{\mathcal{H}}^2  X_i^{-1} -||\overline{Y}||_{\mathcal{H}}^2\big\}\rightarrow  \sum_{j=1}^{\infty}\langle \Delta_3 e_j,e_j\rangle$ as $\nu\rightarrow \infty$ \textit{a.s.} $[\mathbb{P}]$.
\par
\vspace{.1cm}

 Let us next consider the case $V_i$=$Y_i-\overline{Y}-S_{zy}((Z_i-\overline{Z})S_{zz}^{-1})$ for $i$=$1,\ldots,N$. It follows from SLLN that  $N^{-1}\sum_{i=1}^N ||V_i||_{\mathcal{H}}^2$=$O(1)$  as $\nu\rightarrow\infty$ \textit{a.s.} $[\mathbb{P}]$ because Assumption  $3$  holds. Then, it can be shown using similar arguments as in the $1^{st}$ paragraph of this proof that as $\nu\rightarrow\infty$, 
\begin{align}\label{eq 56}
\begin{split}
&\Gamma_1\rightarrow\Delta_4=(1-\lambda)E_{\mathbb{P}}\bigg\{\bigg(Y_i-E_{\mathbb{P}}(Y_i)-C_{zy}\big(( Z_i-E_{\mathbb{P}}(Z_i))C_{zz}^{-1}\big)\bigg)\otimes \\
&\bigg(Y_i-E_{\mathbb{P}}(Y_i)-C_{zy}\big(( Z_i-E_{\mathbb{P}}(Z_i)) C_{zz}^{-1}\big)\bigg)\bigg\}
\end{split}
\end{align}
 under SRSWOR and LMS sampling design, and 
\begin{align}
\begin{split}
&\Gamma_1\rightarrow\Delta_5= E_{\mathbb{P}}\bigg[\bigg\{Y_i-E_{\mathbb{P}}(Y_i)-C_{zy}\big(( Z_i-E_{\mathbb{P}}(Z_i)) C_{zz}^{-1}\big)+\chi^{-1}X_i\lambda \bigg(E_{\textbf{P}}(X_iY_i)-\\
&E_{\textbf{P}}(X_i)E_{\textbf{P}}(Y_i)-C_{zy}\big(( E_{\textbf{P}}(X_i Z_i)-E_{\mathbb{P}}(X_i)E_{\mathbb{P}}(Z_i)) C_{zz}^{-1}\big) \bigg)(E_{\textbf{P}}(X_i))^{-1}\bigg\}\otimes \\
& \bigg\{Y_i-E_{\mathbb{P}}(Y_i)-C_{zy}\big(( Z_i-E_{\mathbb{P}}(Z_i)) C_{zz}^{-1}\big)+\chi^{-1}X_i\lambda \bigg(E_{\textbf{P}}(X_iY_i)-E_{\textbf{P}}(X_i)E_{\textbf{P}}(Y_i)\\
&-C_{zy}\big(( E_{\textbf{P}}(X_i Z_i)-E_{\mathbb{P}}(X_i)E_{\mathbb{P}}(Z_i)) C_{zz}^{-1}\big) \bigg)(E_{\textbf{P}}(X_i))^{-1}\bigg\}\bigg\{X_i^{-1}E_{\mathbb{P}}(X_i)-\lambda\bigg\}\bigg]
\end{split}
\end{align}
 under any HE$\pi$PS sampling design with respect to the HS norm \textit{a.s.} $[\mathbb{P}]$. Here, $C_{zy}$=$E_{\mathbb{P}}(Z_i-E_{\mathbb{P}}(Z_i))\otimes (Y_i-E_{\mathbb{P}}(Y_i))$ and $C_{zz}$=$E_{\mathbb{P}}(Z_i-E_{\mathbb{P}}(Z_i))^T(Z_i-E_{\mathbb{P}}(Z_i))$. Thus as $\nu\rightarrow\infty$, $\Gamma_1\rightarrow\Sigma_1$ with $\Sigma_1$=$\Delta_4$ under SRSWOR and LMS sampling design, and $\Gamma_1\rightarrow\Sigma_1$ with $\Sigma_1$=$\Delta_5$ under any HE$\pi$PS sampling design \textit{a.s.} $[\mathbb{P}]$. Note that $\sum_{j=1}^{\infty}\langle \Delta_4 e_j,e_j\rangle$= $(1-\lambda)E_{\mathbb{P}}\big|\big|Y_i-E_{\mathbb{P}}(Y_i)-C_{zy}(( Z_i-E_{\mathbb{P}}(Z_i))C_{zz}^{-1})||_{\mathcal{H}}^2<\infty$, and $\sum_{j=1}^{\infty} \langle \Delta_5 e_j,e_j\rangle$=$E_{\mathbb{P}}\big[\big|\big|Y_i-E_{\mathbb{P}}(Y_i)-C_{zy}\big(( Z_i-E_{\mathbb{P}}(Z_i)) C_{zz}^{-1}\big)+\chi^{-1}X_i\lambda \big\{E_{\textbf{P}}(X_iY_i)-E_{\textbf{P}}(X_i)E_{\textbf{P}}(Y_i)-C_{zy}\big(( E_{\textbf{P}}(X_i Z_i)-E_{\mathbb{P}}(X_i)E_{\mathbb{P}}(Z_i)) C_{zz}^{-1}\big) \big\}(E_{\textbf{P}}(X_i))^{-1}\big|\big|_{\mathcal{H}}^2 \big\{X_i^{-1}E_{\mathbb{P}}(X_i)-\lambda\big\}\big]<\infty$  since Assumptions  $2$ and $3$  hold. Then, it can be shown in a similar way as in the $1^{st}$ paragraph of this proof that $\sum_{j=1}^{\infty}\langle \Gamma_1 e_j,e_j\rangle\rightarrow \sum_{j=1}^{\infty}\langle \Delta_4 e_j,e_j\rangle$ under SRSWOR and LMS sampling design, and $\sum_{j=1}^{\infty}\langle \Gamma_1 e_j,e_j\rangle \rightarrow\sum_{j=1}^{\infty}\langle \Delta_5 e_j,e_j\rangle$ under any HE$\pi$PS sampling design as $\nu\rightarrow\infty$ \textit{a.s.} $[\mathbb{P}]$. Further, it can be shown using the same line of argument as in the $2^{nd}$ paragraph of this proof that for RHC sampling design, 
\begin{align}
\begin{split}
&\Gamma_2\rightarrow  \Delta_6=c E_{\mathbb{P}}(X_i)E_{\mathbb{P}}\bigg\{\bigg(Y_i-E_{\mathbb{P}}(Y_i)-C_{zy}\big(( Z_i-E_{\mathbb{P}}(Z_i))C_{zz}^{-1}\big)\bigg)\otimes \\
&\bigg(Y_i-E_{\mathbb{P}}(Y_i)-C_{zy}\big(( Z_i-E_{\mathbb{P}}(Z_i)) C_{zz}^{-1}\big)\bigg)X_i^{-1}\bigg\}
\end{split}
\end{align} 
 with respect to the HS norm, and $\sum_{j=1}^{\infty}\langle \Gamma_2 e_j,e_j\rangle\rightarrow  \sum_{j=1}^{\infty}\langle \Delta_6 e_j,e_j\rangle$ as $\nu\rightarrow\infty$ \textit{a.s.} $[\mathbb{P}]$. Thus $\Sigma_2$=$ \Delta_6$ in this case.
\end{proof}
\begin{proof}[\textbf{Proof of Lemma $4$}]
Note that $(\langle S_1, e_1\rangle,\ldots,\langle S_1, e_r\rangle)$=$\sqrt{n}(\hat{\overline{W}}_1-\overline{W})$. Let us first consider SRSWOR, LMS sampling design and any HE$\pi$PS sampling design. Note that under the above-mentioned sampling designs, $\Gamma_{1,r}\rightarrow \Sigma_{1,r}$ as $\nu\rightarrow\infty$ \textit{a.s.} $[\mathbb{P}]$ because $\Gamma_1\rightarrow \Sigma_1$ under these sampling designs as $\nu\rightarrow\infty$ \textit{a.s.} $[\mathbb{P}]$ in view of  Lemma $3$.  Moreover, $\Sigma_{1,r}$ is a n.n.d. matrix since $\Sigma_1$ is a n.n.d. operator. Now, consider the case, when $\Sigma_{1,r}$ is p.d. Then, under the above sampling designs, $m \Gamma_{1,r} m^T>0$ for any $m\in\mathbb{R}^r$ and $m\neq 0$, and all sufficiently large $\nu$ \textit{a.s.} $[\mathbb{P}]$. For any $\epsilon>0$, define $H(\epsilon,m)$=$\{1\leq i \leq N: |m(W_i-\tilde{T} \pi_i)^T|>\epsilon \pi_i N(n^{-1} m \Gamma_{1,r} m^T)^{1/2}\}$,  $J(\epsilon,m)$=$\sum_{i \in H(\epsilon,m)}(m(W_i-\tilde{T} \pi_i)^T)^2 (\pi^{-1}_i-1) (n^{-1}N^2 m \Gamma_{1,r} m^T)^{-1}$ and   $\tilde{W}_i$= $(n(N\pi_i)^{-1})W_i-(nN^{-1})\tilde{T}$, where $\tilde{T}$=$\sum_{i=1}^N W_i(1-\pi_i)(\sum_{i=1}^N \pi_i(1-\pi_i))^{-1}$. Then, given any $\delta>0$,  
\begin{equation}\label{eq 14}
J(\epsilon,m)\leq (m \Gamma_{1,r} m^T)^{-(1+\delta/2)} n^{-\delta/2}\epsilon^{-\delta}N^{-1}||m||^{2+\delta}\sum_{i=1}^N ||\tilde{W}_i||^{2+\delta}\big(n^{-1}N\pi_i\big) 
\end{equation}
 since $|m\tilde{W}_i^T|(\epsilon (n m \Gamma_{1,r}m^T)^{1/2})^{-1}>1$ for all $i\in H(\epsilon,m)$. It follows from Jensen's inequality and the fact $\sum_{i=1}^N \pi_i$=$n$ that 
\begin{equation}\label{eq 15}
N^{-1}\sum_{i=1}^N||\tilde{W}_i||^{2+\delta} \big(n^{-1}N\pi_i\big)\leq 2^{\delta+1}\bigg(N^{-1} \sum_{i=1}^N \bigg|\bigg|n W_i(N\pi_i)^{-1}\bigg|\bigg|^{2+\delta}\big(n^{-1}N\pi_i\big)+\bigg|\bigg|nN^{-1}\tilde{T}\bigg|\bigg|^{2+\delta}\bigg).
\end{equation} 
 It also follows from SLLN and Assumption  $3$  that  $N^{-1}\sum_{i=1}^N||W_i||^{2+\delta}$=$O(1)$  as $\nu\rightarrow\infty$ \textit{a.s.} $[\mathbb{P}]$ for any $0<\delta\leq 2$. Further, under all of the above sampling designs,  $n^{-1}\sum_{i=1}^N\pi_i(1-\pi_i)$  is bounded away from $0$ as $\nu\rightarrow\infty$ \textit{a.s.} $[\mathbb{P}]$ because  ($10$) in Lemma $1$  holds under these sampling designs, and Assumption  $1$ \color {black} holds. Therefore, 
\begin{equation}\label{eq 16}
N^{-1} \sum_{i=1}^N \bigg|\bigg|n W_i (N\pi_i)^{-1}\bigg|\bigg|^{2+\delta}\big(n^{-1}N\pi_i\big)=O(1)\text{ and }\bigg|\bigg|nN^{-1}\tilde{T}\bigg|\bigg|^{2+\delta}=O(1),
\end{equation}
and hence $N^{-1}\sum_{i=1}^N ||\tilde{W}_i||^{2+\delta}(N\pi_i/n)$=$O(1)$  by \eqref{eq 15} as $\nu\rightarrow\infty$ under the above sampling designs \textit{a.s.} $[\mathbb{P}]$. Then, under these sampling designs, $J(\epsilon,m)\rightarrow 0$ as $\nu\rightarrow\infty$ for any $\epsilon>0$ and $m\neq 0$ \textit{a.s.} $[\mathbb{P}]$. Hence, $\inf\{\epsilon>0:J(\epsilon,m)\leq \epsilon\}\rightarrow 0$ as $\nu\rightarrow\infty$ for any $m\neq 0$ \textit{a.s.} $[\mathbb{P}]$. Thus given any $m\neq 0$, the H\'ajek-Lindeberg condition (see \cite{MR1624693}) holds for $\{m W_i^T\}_{i=1}^N$ \textit{a.s.} $[\mathbb{P}]$. Note that SRSWOR is a high entropy sampling design since it is a rejective sampling design. Further, it follows from  Lemma $1$  that LMS sampling design is a high entropy sampling design. Moreover, under all the sampling designs discussed in this paragraph, $\sum_{i=1}^N \pi_i(1-\pi_i)\rightarrow \infty$ since  $n^{-1}\sum_{i=1}^N\pi_i(1-\pi_i)$  is bounded away from $0$ as $\nu\rightarrow\infty$ \textit{a.s.} $[\mathbb{P}]$. Therefore, by Theorem $5$ in \cite{MR1624693}, $\sqrt{n} m(\hat{\overline{W}}_1-\overline{W})^T\xrightarrow{\mathcal{L}}N(0,m \Sigma_{1,r} m^T)$ as $\nu\rightarrow\infty$ for any $m\neq 0$ under these sampling designs \textit{a.s.} $[\mathbb{P}]$. This implies that under these sampling designs, $\sqrt{n}(\hat{\overline{W}}_1-\overline{W})\xrightarrow{\mathcal{L}}$ $N_r(0,\Sigma_{1,r})$ as $\nu\rightarrow\infty$ \textit{a.s.} $[\mathbb{P}]$.
\par
\vspace{.2cm}

Next, consider the case, when $\Sigma_{1,r}$ is a positive semi definite (p.s.d.) matrix. Let $A_1$=$\{m\neq 0:m\Sigma_{1,r} m^T>0\}$ and $A_2$=$\{m\neq 0:m\Sigma_{1,r} m^T=0\}$. Then, under the sampling designs mentioned in the  preceding  paragraph, $\sqrt{n}m(\hat{\overline{W}}_1-\overline{W})^T\xrightarrow{\mathcal{L}}N(0,m \Sigma_{1,r} m^T)$ as $\nu\rightarrow\infty$ for any $m\in A_1$ \textit{a.s.} $[\mathbb{P}]$ in the same way as argued above. Next, suppose that $P(s,\omega)$ denotes one of these sampling designs, and $Q(s, \omega)$ is a rejective sampling design with inclusion probabilities equal to those of $P(s,\omega)$ (cf. \cite{MR1624693}). Note that under $Q(s,\omega)$, $var(\sqrt{n}m(\hat{\overline{W}}_1-\overline{W})^T)$=$m\Gamma_{1,r}m^T(1+h)$ (see Theorem $6.1$ in \cite{MR0178555}) for any $\omega$ and $m$, where $h\rightarrow 0$ as $\nu\rightarrow\infty$ if $\sum_{i=1}^N \pi_i(1-\pi_i)\rightarrow\infty$ as $\nu\rightarrow\infty$. Recall from the preceding paragraph that $\sum_{i=1}^N \pi_i(1-\pi_i)\rightarrow\infty$ as $\nu\rightarrow\infty$ under $P(s,\omega)$ \textit{a.s.} $[\mathbb{P}]$. Therefore, $\sum_{i=1}^N \pi_i(1-\pi_i)\rightarrow\infty$ as $\nu\rightarrow\infty$ under $Q(s,\omega)$ \textit{a.s.} $[\mathbb{P}]$. Next, note that $\Gamma_{1,r}$ depends on the sampling design only through the inclusion probabilities, and $\Gamma_{1,r}\rightarrow\Sigma_{1,r}$ as $\nu\rightarrow\infty$ under $P(s,\omega)$ \textit{a.s.} $[\mathbb{P}]$ as mentioned in the previous paragraph. Therefore, $m\Gamma_{1,r}m^T\rightarrow 0$ as $\nu\rightarrow\infty$ for any $m\in A_2$ under $Q(s,\omega)$ \textit{a.s.} $[\mathbb{P}]$. Hence, $\sqrt{n}m(\hat{\overline{W}}_1-\overline{W})^T$=$o_p(1)$ as $\nu\rightarrow\infty$ for any $m\in A_2$ under $Q(s,\omega)$ \textit{a.s.} $[\mathbb{P}]$. Now, it follows from Lemmas $2$ and $3$ in \cite{MR1624693} that 
\begin{align}\label{eq 17}
\begin{split}
&\sum_{s\in A}P(s,\omega)\leq\sum_{s\in A}Q(s,\omega)+\sum_{s\in\mathcal{S}}|P(s,\omega)-Q(s,\omega)|\leq \sum_{s\in A}Q(s,\omega)+(2D(P||Q))^{1/2}\\
&\leq \sum_{s\in A}Q(s,\omega)+(2D(P||R))^{1/2}, 
\end{split}
\end{align}
 where $A$=$\{s\in\mathcal{S}:|\sqrt{n}m(\hat{\overline{W}}_1-\overline{W})^T|>\epsilon\}$ for $\epsilon>0$, and $R(s,\omega)$ is any other rejective sampling design. Since $P(s,\omega)$ is a high entropy sampling design as discussed earlier in this proof, there exists a rejective sampling design $R(s,\omega)$ such that $D(P||R)\rightarrow 0$ as $\nu\rightarrow\infty$ \textit{a.s.} $[\mathbb{P}]$. Then, under $P(s,\omega)$, $\sqrt{n}m(\hat{\overline{W}}_1-\overline{W})^T$=$o_p(1)$ as $\nu\rightarrow\infty$ for any $m\in A_2$ \textit{a.s.} $[\mathbb{P}]$. Therefore, under $P(s,\omega)$, as $\nu\rightarrow\infty$, $\sqrt{n}m(\hat{\overline{W}}_1-\overline{W})^T\xrightarrow{\mathcal{L}}N(0,m \Sigma_{1,r} m^T)$ for any $m\neq 0$, and hence $\sqrt{n}(\hat{\overline{W}}_1-\overline{W})\xrightarrow{\mathcal{L}}N_r(0,\Sigma_{1,r})$ \textit{a.s.} $[\mathbb{P}]$. 
\par
\vspace{.2cm}

Next, note that $(\langle S_2, e_1\rangle,\ldots,\langle S_2, e_r\rangle)$=$\sqrt{n}(\hat{\overline{W}}_2-\overline{W})$. Also, note that $\Gamma_{2,r}\rightarrow \Sigma_{2,r}$ as $\nu\rightarrow\infty$ \textit{a.s.} $[\mathbb{P}]$ since $\Gamma_2\rightarrow \Sigma_2$ as $\nu\rightarrow\infty$ \textit{a.s.} $[\mathbb{P}]$ in view of  Lemma $3$.  Moreover, $\Sigma_{2,r}$ is a n.n.d. matrix since $\Sigma_2$ is a n.n.d. operator. Let us consider the case, when $\Sigma_{2,r}$ is p.d. Then, $m \Gamma_{2,r} m^T>0$ for any $m\neq 0$ and all sufficiently large $\nu$  \textit{a.s.} $[\mathbb{P}]$. Let us define 
\begin{align}\label{eq 18}
\begin{split}
&J(m)=n\gamma\bigg(\max_{1\leq i\leq N} X_i\bigg)\bigg(N^{-1}\sum_{i=1}^n N_i^3(N_i-1)\sum_{i=1}^N\big(m (W_i \overline{X}X_i^{-1}-\overline{W})^T\big)^4 X_i\bigg)^{1/2}\times\\
&\bigg(\overline{X}^{3/2}\sum_{i=1}^n N_i(N_i-1) m \Gamma_{2,r} m^T\bigg)^{-1},
\end{split}
\end{align}
where $\gamma$=$\sum_{i=1}^n N_i(N_i-1)(N(N-1))^{-1}$ with $N_i$ being the size of the $i^{th}$ group formed randomly in the first step of the RHC sampling design (see Appendix  A)  for $i$=$1,\ldots,n$.  Note that  $N^{-1}\sum_{i=1}^N||W_i||^4$=$O(1)$  as $\nu\rightarrow\infty$ \textit{a.s.} $[\mathbb{P}]$ by Assumption  $3$  and SLLN. Then, we have 
\begin{equation}\label{eq 19}
\bigg(N^{-1}\sum_{i=1}^N\big(m(W_i \overline{X}X_i^{-1}-\overline{W})^T\big)^4 (X_i\overline{X}^{-1})\bigg)^{1/2}=O(1)\text{ and }\overline{X}^{-1}\max_{1\leq i\leq N} X_i =O(1) 
\end{equation}
 as $\nu\rightarrow\infty$ \textit{a.s.} $[\mathbb{P}]$ because Assumption  $2$  holds. Recall from the paragraph preceding Proposition  $3.2$ in Section $3$ of the main text  that $N_i$'s are taken as in  ($4$) of the main article.  Then, by Assumption  $1$, 
\begin{equation}\label{eq 20}
\bigg(\sum_{i=1}^n N_i^3(N_i-1)\bigg)^{1/2}\bigg(\sum_{i=1}^n N_i(N_i-1)\bigg)^{-1}=O(1/\sqrt{n})\text{ and }n\gamma=O(1) 
\end{equation}
 as $\nu\rightarrow\infty$. Therfore, $J(m)\rightarrow 0$ as $\nu\rightarrow\infty$ by \eqref{eq 18}--\eqref{eq 20}, and hence C$1$ in \cite{MR844032} holds for $\{m W_i^T\}_{i=1}^N$ for any $m\neq 0$ \textit{a.s.} $[\mathbb{P}]$. Then, under RHC sampling design, $\sqrt{n}m(\hat{\overline{W}}_2-\overline{W})^T\xrightarrow{\mathcal{L}}N(0,(m \Sigma_{2,r} m^T)$ as $\nu\rightarrow\infty$ for any $m\neq 0$ \textit{a.s.} $[\mathbb{P}]$ by Theorem $2.1$ in \cite{MR844032}. Therefore, under RHC sampling design, $\sqrt{n}(\hat{\overline{W}}_2-\overline{W})\xrightarrow{\mathcal{L}}N_r(0,\Sigma_{2,r})$ as $\nu\rightarrow\infty$ \textit{a.s.} $[\mathbb{P}]$.
\par
\vspace{.2cm}

Next, consider the case, when $\Sigma_{2,r}$ is p.s.d. Let $A_1$=$\{m\neq 0:m\Sigma_{2,r} m^T>0\}$ and $A_2$=$\{m\neq 0:m\Sigma_{2,r} m^T=0\}$. Then, under RHC sampling design, $\sqrt{n}m(\hat{\overline{W}}_2-\overline{W})^T\xrightarrow{\mathcal{L}}N(0,m \Sigma_{2,r} m^T)$ as $\nu\rightarrow\infty$ for any $m\in A_1$ \textit{a.s.} $[\mathbb{P}]$ in the same way as above. Under RHC sampling design, $var(\sqrt{n}m(\hat{\overline{W}}_2-\overline{W})^T)$=$m\Gamma_{2,r}m^T$ (see \cite{MR844032}) for any $\omega$ and $m$. Note that $m\Gamma_{2,r}m^T\rightarrow 0$ as $\nu\rightarrow\infty$ for any $m\in A_2$ \textit{a.s.} $[\mathbb{P}]$. Then, under RHC sampling design, $\sqrt{n}m(\hat{\overline{W}}_2-\overline{W})^T$=$o_p(1)$ as $\nu\rightarrow\infty$ for any $m\in A_2$ \textit{a.s.} $[\mathbb{P}]$. Therefore, under RHC sampling design, as $\nu\rightarrow\infty$, $\sqrt{n}m(\hat{\overline{W}}_2-\overline{W})^T\xrightarrow{\mathcal{L}}N(0,m \Sigma_{2,r} m^T)$ for any $m\neq 0$, and hence $\sqrt{n}(\hat{\overline{W}}_2-\overline{W})\xrightarrow{\mathcal{L}}N_r(0,\Sigma_{2,r})$ \textit{a.s.} $[\mathbb{P}]$. 
\end{proof}
\begin{proof}[\textbf{Proof of Lemma $5$}]
Let us first consider the case, when $P(s,\omega)$ is one of SRSWOR, LMS sampling design and any HE$\pi$PS sampling design. Suppose that $Q(s,\omega)$ is as described in the $2^{nd}$ paragraph of the proof of Lemma $4$.  Then, following similar arguments as in the proof of Theorem $6.1$ in \cite{MR0178555}, we can show that 
\begin{equation}\label{eq 21}
E\langle S_1,e_j\rangle^2=(nN^{-2})\sum_{i=1}^N\langle V_i-T\pi_i,e_j\rangle^2(\pi_i^{-1}-1)(1+h)=\langle \Gamma_1 e_j,e_j\rangle (1+h) 
\end{equation}
 under $Q(s,\omega)$ for any $\omega$ and $j\geq 1$. Here, $h$ does not depend on $\{e_j\}_{j=1}^{\infty}$, and $h\rightarrow 0$ as $\nu\rightarrow\infty$ whenever $\sum_{i=1}^N \pi_i(1-\pi_i)\rightarrow\infty$ as $\nu\rightarrow\infty$. Recall from the $2^{nd}$ paragraph in the proof of Lemma $4$  that $\sum_{i=1}^N \pi_i(1-\pi_i)\rightarrow\infty$ as $\nu\rightarrow\infty$ under $Q(s,\omega)$ \textit{a.s.} $[\mathbb{P}]$. It follows from Lemma $3$  that under $P(s,\omega)$, $\Gamma_1\rightarrow \Sigma_1$ with respect to the HS norm and $\sum_{j=1}^{\infty}\langle\Gamma_1 e_j,e_j\rangle\rightarrow \sum_{j=1}^{\infty}\langle\Sigma_1 e_j,e_j\rangle$ as $\nu\rightarrow\infty$ \textit{a.s.} $[\mathbb{P}]$. Therefore, $\Gamma_1\rightarrow \Sigma_1$ and $\sum_{j=1}^{\infty}\langle\Gamma_1 e_j,e_j\rangle\rightarrow \sum_{j=1}^{\infty}\langle\Sigma_1 e_j,e_j\rangle$ as $\nu\rightarrow\infty$ under $Q(s,\omega)$ \textit{a.s.} $[\mathbb{P}]$ because $\Gamma_1$ depends on the sampling design only through inclusion probabilities, and $P(s,\omega)$ and $Q(s,\omega)$ have the same inclusion probabilities. Thus as $\nu\rightarrow\infty$, $E\langle S_1,e_j\rangle^2\rightarrow \langle\Sigma_1 e_j,e_j\rangle$ for any $j\geq 1$, and $\sum_{j=1}^{\infty}E\langle S_1,e_j\rangle^2\rightarrow\sum_{j=1}^{\infty}\langle\Sigma_1 e_j,e_j\rangle$ under $Q(s,\omega)$ \textit{a.s.} $[\mathbb{P}]$. Then, following the same line of arguments as in the proof of Theorem $1.1$ in \cite{kundu2000central}, we can say that 
\begin{equation}\label{eq 22}
\overline{\lim}_{\nu\rightarrow\infty}\sum_{s\in B_{1,r}}Q(s,\omega)\leq\sum_{j=r+1}^{\infty}\langle\Sigma_1 e_j,e_j\rangle\epsilon^{-2}
\end{equation} 
 \textit{a.s.} $[\mathbb{P}]$ for any $r\geq 1$. Therefore, $\lim_{r\rightarrow\infty}\overline{\lim}_{\nu\rightarrow\infty}\sum_{s\in B_{1,r}}Q(s,\omega)$=$0$ \textit{a.s.} $[\mathbb{P}]$. Further, it can be shown that $\lim_{r\rightarrow\infty}\overline{\lim}_{\nu\rightarrow\infty}\sum_{s\in B_{1,r}}P(s,\omega)$=$0$ \textit{a.s.} $[\mathbb{P}]$ in the same way as the result $\sqrt{n}m(\hat{\overline{W}}_1-\overline{W})^T$=$o_p(1)$ as $\nu\rightarrow\infty$ under $P(s,\omega)$ \textit{a.s.} $[\mathbb{P}]$ is shown in the $2^{nd}$ paragraph of the proof of Lemma $4$. 
\par
\vspace{.2cm}

Let us next consider the case, when $P(s,\omega)$ is RHC sampling design. Note that 
\begin{equation}\label{eq 23}
E\langle S_2,e_j\rangle^2=(n\gamma)(\overline{X} N^{-1})\sum_{i=1}^N\bigg(\langle V_i,e_j\rangle-\langle \overline{V},e_j\rangle X_i\overline{X}^{-1}\bigg)^2 X_i^{-1}=\langle \Gamma_2 e_j,e_j\rangle 
\end{equation}
 under RHC sampling design for any $j\geq 1$ and $\omega$ (cf. \cite{MR844032}). Also, note that as $\nu\rightarrow\infty$, $\Gamma_2\rightarrow\Sigma_2$ with respect to the HS norm and $\sum_{j=1}^{\infty}\langle\Gamma_2 e_j,e_j\rangle\rightarrow \sum_{j=1}^{\infty}\langle\Sigma_2 e_j,e_j\rangle$ \textit{a.s.} $[\mathbb{P}]$ in view of Lemma $3$.  Then, under RHC sampling design, as $\nu\rightarrow\infty$, $E\langle S_2,e_j\rangle^2\rightarrow \langle\Sigma_2 e_j,e_j\rangle$ for any $j\geq 1$, and $\sum_{j=1}^{\infty}E\langle S_2,e_j\rangle^2\rightarrow\sum_{j=1}^{\infty}\langle\Sigma_2 e_j,e_j\rangle$  \textit{a.s.} $[\mathbb{P}]$. Therefore, $\lim_{r\rightarrow\infty}\overline{\lim}_{\nu\rightarrow\infty}\sum_{s\in B_{2,r}}P(s,\omega)$=$0$ \textit{a.s.} $[\mathbb{P}]$ using similar arguments as in the proof of Theorem $1.1$ in \cite{kundu2000central}.
\end{proof}
\begin{proof}[\textbf{Proof of Lemma $6$}]
Note that $\{V_i^{\sharp}\}_{i=1}^N$ are elements of either an infinite dimensional separable Hilbert space or a finite dimensional Euclidean space. Let $\{e_j^{\sharp}\}$ be an orthonormal basis of that space. Further, note that  $N^{-1}\sum_{i=1}^N||V_i^{\sharp}||_{\mathcal{G}}^2$=$O(1)$  as $\nu\rightarrow\infty$ \textit{a.s.} $[\mathbb{P}]$ by SLLN and Assumption  $3$.  Now, suppose that $P(s,\omega)$ is one of SRSWOR, LMS sampling design and any HE$\pi$PS sampling design, and $Q(s,\omega)$ is the corresponding rejective sampling design as described in the $2^{nd}$ paragraph of the proof of Lemma $4$.  Then, one can show that 
\begin{equation}\label{eq 24}
E||S_1^{\sharp}||_{\mathcal{G}}^2=E\bigg(\sum_{j}\bigg\langle S_1^{\sharp}, e_j^{\sharp}\bigg\rangle^2\bigg)=(n N^{-2})\sum_{j}\sum_{i=1}^N\bigg\langle V_i^{\sharp}-T^{\sharp}\pi_i,e_j^{\sharp}\bigg\rangle^2(\pi_i^{-1}-1)(1+h) 
\end{equation}
 for any $\omega$ under $Q(s,\omega)$ in the same way as the derivation of $E\langle S_1,e_j\rangle^2$=$\langle \Gamma_1 e_j,e_j\rangle (1+h)$ in the proof of Lemma $5$.  Here,  $T^{\sharp}$=$\sum_{i=1}^N V_i^{\sharp}(1-\pi_i)\big(\sum_{i=1}^N \pi_i(1-\pi_i)\big)^{-1}$,  $h$ does not depend on $\{e_j^{\sharp}\}$, and $h\rightarrow 0$ as $\nu\rightarrow\infty$ if $\sum_{i=1}^N \pi_i(1-\pi_i)\rightarrow\infty$ as $\nu\rightarrow\infty$. Note that  ($10$) in Lemma $1$  holds under $Q(s,\omega)$ because  ($10$) in Lemma $1$  holds under $P(s,\omega)$ by Lemma $1$,  and $P(s,\omega)$ and $Q(s,\omega)$ have the same inclusion probabilities. Then, $\sum_{i=1}^N \pi_i(1-\pi_i)\rightarrow\infty$ as $\nu\rightarrow\infty$ under $Q(s,\omega)$ \textit{a.s.} $[\mathbb{P}]$.  Therefore, as $\nu\rightarrow\infty$, 
\begin{align}\label{eq 25}
\begin{split}
&(nN^{-2})\sum_{j}\sum_{i=1}^N\bigg\langle V_i^\sharp-T^{\sharp}\pi_i,e_j^\sharp\bigg\rangle^2(\pi_i^{-1}-1)(1+h)=(nN^{-2})\sum_{i=1}^N || V_i^{\sharp}-T^\sharp\pi_i||_{\mathcal{G}}^2\times\\
&(\pi_i^{-1}-1)(1+h)=(nN^{-2})\bigg[\sum_{i=1}^N||V_i^\sharp||_{\mathcal{G}}^2(\pi_i^{-1}-1)-||T^{\sharp}||_{\mathcal{G}}^2\sum_{i=1}^N\pi_i(1-\pi_i)\bigg](1+h)\leq\\
&(nN^{-2})\sum_{i=1}^N \pi_i^{-1}||V_i^\sharp||_{\mathcal{G}}^2(1+h)=O(1)
\end{split}
\end{align}
 under $Q(s,\omega)$ \textit{a.s.} $[\mathbb{P}]$ since  $N^{-1}\sum_{i=1}^N||V_i^{\sharp}||_{\mathcal{G}}^2$=$O(1)$  as $\nu\rightarrow\infty$ \textit{a.s.} $[\mathbb{P}]$. Hence,  $E||S_1^{\sharp}||_{\mathcal{G}}^2$=$O(1)$  as $\nu\rightarrow\infty$ under $Q(s,\omega)$ \textit{a.s.} $[\mathbb{P}]$. Thus  $||S_1^{\sharp}||_{\mathcal{G}}$=$O_p(1)$  as $\nu\rightarrow\infty$ under $Q(s,\omega)$ \textit{a.s.} $[\mathbb{P}]$. Now, it can be shown that $||S_1^\sharp||_{\mathcal{G}}=O_p(1)$  as $\nu\rightarrow\infty$ under $P(s,\omega)$ \textit{a.s.} $[\mathbb{P}]$ in the same way as the result $\sqrt{n}m(\hat{\overline{W}}_1-\overline{W})^T$=$o_p(1)$ as $\nu\rightarrow\infty$ under $P(s,\omega)$ \textit{a.s.} $[\mathbb{P}]$ is shown in the $2^{nd}$ paragraph of the proof of Lemma $4$. 
\par
\vspace{.2cm}

Next, note that under RHC sampling design, as $\nu\rightarrow\infty$, 
\begin{align}\label{eq 26}
\begin{split}
&E||S_2^{\sharp}||_{\mathcal{G}}^2=E\bigg(\sum_{j}\bigg\langle S_2^{\sharp}, e_j^{\sharp}\bigg\rangle^2\bigg)=(n\gamma)(\overline{X}N^{-1})\sum_j \sum_{i=1}^N\bigg(\langle V_i^{\sharp},e_j^{\sharp}\rangle-\langle \overline{V}^{\sharp},e_j^{\sharp}\rangle X_i\overline{X}^{-1}\bigg)^2X_i^{-1}\\
&\leq (n\gamma)\sum_{j} N^{-1}\sum_{i=1}^N \langle V_i^{\sharp},e_j^{\sharp}\rangle^2\overline{X} X_i^{-1}\leq (n\gamma) N^{-1}\sum_{i=1}^N || V_i^{\sharp}||_{\mathcal{G}}^2\overline{X} X_i^{-1}=O(1) 
\end{split}
\end{align}
 \textit{a.s.} $[\mathbb{P}]$ because  $N^{-1}\sum_{i=1}^N|| V_i^{\sharp}||_{\mathcal{G}}^2$=$O(1)$  as $\nu\rightarrow\infty$ \textit{a.s.} $[\mathbb{P}]$, and Assumption  $2$  holds. Also, note that $n\gamma$=$O(1)$ as $\nu\rightarrow\infty$ since $N_i$'s are taken as in  ($4$) of the main text.  Therefore,  $||S_2^{\sharp}||_{\mathcal{G}}$=$O_p(1)$  as $\nu\rightarrow\infty$ under RHC sampling design \textit{a.s.} $[\mathbb{P}]$.
\end{proof}
\vskip .65cm
\noindent
Anurag Dey\\
{\it Indian Statistical Institute}
\vskip 2pt
\noindent
E-mail: deyanuragsaltlake64@gmail.com
\vskip 2pt

\noindent
Probal Chaudhuri\\
{\it Indian Statistical Institute}
\vskip 2pt
\noindent
E-mail:probalchaudhuri@gmail.com
% \vskip .3cm
%\centerline{(Received ???? 20??; accepted ???? 20??)}\par
%\end{document}
%%%%%%%%%%%%%%%%%%%%%%%%%%%%%%%%%%%%%%%%%%%%%%%%%%%%%%%%%%%%%%%%%%%%%%%%%%%%%%%%%%%%%%%%%%%%%%%%%%%%%%%%%%%%%%%%%%%%%%%%%%%%
%%%%%%%%%%%%%%%%%%%%%%%%%%%%%%%%%%%%%%%%%%%%%%%%%%%%%%%%%%%%%%%%%%%%%%%%%%%%%%%%%%%%%%%%%%%%%%%%%%%%%%%%%%%%%%%%%%%%%%%%%%%%